\documentclass[12pt]{article}
\usepackage{graphicx,psfrag,epsfig, color,float}
\usepackage{amssymb,amsmath,amscd,amsthm}
\usepackage{verbatim}
\usepackage{graphicx,psfrag,epsfig}

\usepackage{graphicx}
\usepackage{xfrac}
\usepackage[active]{srcltx}
\usepackage{mathtools}
\usepackage{relsize}
\usepackage{bbm}

\newtheorem{theorem}{Theorem}[section]
\newtheorem{corollary}[theorem]{Corollary}

\newtheorem{lemma}[theorem]{Lemma}
\newtheorem{proposition}[theorem]{Proposition}
\newtheorem{definition}[theorem]{Definition}

\newtheorem{remark}[theorem]{Remark}

\def\be{\color{black}}
\def\br{\color{red}}

\date{}

\begin{document}

\date{}
\title{Metastable Distributions of Semi-Markov Processes}
\author{
L. Koralov\footnote{Dept of Mathematics, University of Maryland,
College Park, MD 20742, koralov@umd.edu},
I. Mohammed Imtiyas\footnote{Dept of Mathematics and Statistics, Mount Holyoke College,  MA 01075, iami@mtholyoke.edu} 
} \maketitle

\begin{abstract}
In this paper, we consider semi-Markov processes whose transition times and transition probabilities depend on a small parameter $\varepsilon$.
Understanding the asymptotic behavior of such processes is needed in order to study the asymptotics of various randomly perturbed dynamical and stochastic systems. The
long-time behavior of a semi-Markov process $X^\varepsilon_t$ depends on how the point $(1/\varepsilon, t(\varepsilon))$ approaches 
infinity. We introduce the notion of complete asymptotic regularity (a certain asymptotic condition on transition probabilities and transition times), originally developed for parameter-dependent Markov chains, which ensures the existence of the metastable distribution for each initial point and a given time scale $t(\varepsilon)$.  
The result may be viewed as a generalization of the ergodic theorem to the case of parameter-dependent semi-Markov processes.
\end{abstract}

{2020 Mathematics Subject Classification Numbers: 60J27, 60K15.
}
\\

{ Keywords: Semi-Markov Processes, Markov Renewal Processes, Asymptotic Regularity, Metastable Distributions. 
}

\section{Introduction.} \label{intro}

In this paper, we describe the asymptotic behavior of parameter-dependent semi-Markov processes. An important application of the results on semi-Markov processes can be found in the study of metastability in non-equilibrium systems. On an intuitive level, metastability of a system is characterized by long periods of random oscillations near an apparent  equilibrium  followed
by a sudden transition towards another more stable equilibrium.   Metastability is exhibited in genetics, molecular dynamics, population dynamics, etc. (see \cite{EREN}, \cite{EST}, \cite{KAUF} and \cite{SCH}). See also \cite{BOV1} and \cite{OLIV} for monographs on metastability that discuss areas of applications.      

Consider a dynamical system or a diffusion process that has $M$ ergodic invariant measures that are supported on disjoint invariant sets. Each invariant set is assumed to  be stable in a certain sense (e.g., an asymptotically stable attractor in the case of a dynamical system). Each invariant set can be associated with a state (or a subset in the state space) of a semi-Markov process, while the transitions between different states are due to large deviations from the unperturbed dynamics. 

Given a process $\Xi_t$ with finitely many invariant sets, the corresponding semi-Markov process is defined by stopping the perturbed 
process $\Xi^\varepsilon_t$ every time it visits a new invariant set (or, in some cases, its $\varepsilon$-dependent neighborhood). Thus the transition kernel and the transition times of the semi-Markov process inherit the $\varepsilon$-dependence from the original system. We will briefly consider several examples where the semi-Markov process corresponds to an underlying randomly perturbed system, but the main result is formulated in an abstract setting. Namely, under certain conditions on the transition times and transition probabilities,
a semi-Markov process exhibits meta-stable behavior. Before describing the assumptions in detail and stating the result precisely, let us
discuss the simpler case of parameter-dependent Markov chains.

Consider a family $X^\varepsilon_t$ of Markov chains on a state space $S = \{1,...,M\}$, where $\varepsilon$ is a small parameter. The time may be continuous or discrete; we focus on the case when $t \in \mathbb{R}^+$, while the case of discrete time is similar. Let $q_{ij}(\varepsilon)$, $i, j \in S$, $i \neq j$,  be the transition rates, i.e.,
\[
\mathbb{P}(X^\varepsilon_{t + \Delta} = j | X^\varepsilon_t = i) = q_{ij}(\varepsilon) (\Delta + o(\Delta))~~~{\rm as}~\Delta \downarrow 0,~~i \neq j.
\]
Let us assume that $q_{ij}(\varepsilon) \downarrow 0$ as $\varepsilon \downarrow 0$ for all $i \neq j$. (This is the case in most systems of interest, and, in any case, can always be achieved by slowing all the transitions by a function of $\varepsilon$.)
We are interested in the behavior of $X^\varepsilon_t$ as $\varepsilon \downarrow 0$ and, simultaneously, $t = t(\varepsilon) \rightarrow~\infty$.
The results on the asymptotic behavior of $X^\varepsilon_t$ can be viewed as a refinement of the ergodic theorem for Markov chains (which concerns the asymptotics with respect to the time variable only) and are closely related to the spectral properties of the transition matrix  (some of the earliest work on
the eigenvalues of matrices with exponentially small entries can be found in \cite{WE}). 
The double limit at hand depends on how the point $(1/\varepsilon, t(\varepsilon))$ approaches infinity. Roughly speaking, one can divide the neighborhood of infinity into a finite number of domains such that $X^\varepsilon_{t(\varepsilon)}$ has a limiting distribution (which depends on the initial point) when $(1/\varepsilon, t(\varepsilon))$ approaches infinity
without leaving a given domain. For different domains, these limits are different. These are referred to as metastable distributions.

Assume, for simplicity, that $ q_{ij}(\varepsilon) > 0$ for each $\varepsilon > 0$ and all $i \neq j$. The chain $X^\varepsilon_t$ is said
to be {\it completely asymptotically regular} if for each $a > 0$ and each $(i_1,...,i_a)$,
$(j_1,...,j_a)$, $(k_1,...,k_a)$, and $(l_1,...,l_a)$  the following finite or infinite limit exists
$$
\lim_{\varepsilon \downarrow 0} \left( \frac{q_{i_1 j_1}(\varepsilon)}{q_{k_1 l_1}(\varepsilon)} \times ... \times \frac{q_{i_a j_a}(\varepsilon)}{q_{k_a l_a}(\varepsilon)}\right) \in [0, \infty],
$$
provided that $i_1 \neq j_1,...,i_a \neq j_a, k_1 \neq l_1,...,k_a \neq l_a$. If the chain is completely asymptotically regular, we can  associate a finite sequence of time scales 
$  0 = \mathfrak{t}_0^{\varepsilon} (i) \ll \mathfrak{t}_1^{\varepsilon} (i) \ll \mathfrak{t}_2^{\varepsilon} (i) \ll \dots  \ll \mathfrak{t}_{n(i)}^{\varepsilon} (i) = \infty $
with each initial state $i \in S$. Here, 
$\mathfrak{t}_{k-1}^{\varepsilon} (i) \ll \mathfrak{t}_k^{\varepsilon} (i) $ means that $\lim_{\varepsilon \downarrow 0} (\mathfrak{t}_{k-1}^{\varepsilon} (i) /\mathfrak{t}_k^{\varepsilon} (i))  = 0.$

Assuming that the process starts at $i$ and
$\mathfrak{t}_{k-1}^{\varepsilon} (i)  \ll t(\varepsilon) \ll \mathfrak{t}_{k}^{\varepsilon} (i) $, $X^{\varepsilon}_{t(\varepsilon)}$ converges to a distribution $\mu_{i, k}$, 
called the  metastable distribution for the initial state $i$ at the time scale $t(\varepsilon)$. This result (generalizing the construction
of \cite{F5}) was proved 
in  \cite{LX} and  \cite{FK1}. \be In the case of discrete-time Markov chains, a closely related result was proved in \cite{BL1}.

The time scales $\mathfrak{t}_{k}^{\varepsilon} (i)$ can be determined in terms of transition times for auxiliary Markov chains that belong to a hierarchy defined by successively reducing the state space, i.e., combining the elements of the state space into subsets (clusters) 
that will serve as states for the chain of higher rank. This construction of the hierarchy of Markov chains generalizes the hierarchy of cycles notion of \cite{F5}.

For a Markov chain, the transition rates $q_{ij}(\varepsilon)$ determine both the time spent in a state $i$ (which is
exponentially distributed with the expectation $1/\sum_{j \neq i} q_{ij}(\varepsilon)$) and the transition probabilities (the
probability of going from $i$ to $j$ is $q_{ij}(\varepsilon)/\sum_{k \neq i} q_{ik}(\varepsilon)$). 
In the case of semi-Markov processes, the transition times and transition probabilities need not be
in this particular relation. Thus the definition of complete asymptotic regularity needs to be modified in order to
incorporate information about the transition times and transition probabilities. In fact, we will see that the state space can be reduced 
inductively to a collection of clusters (until we have one cluster containing all the state space) based on a condition that involves only transition probabilities. In order to determine the effective transition times between clusters, on the other hand, we need both the transition times and transition probabilities of the original process. One difficulty not present when the process was purely Markov is that now the transition time between two states (or parts of the state space) can be asymptotically large even if the  transition probability between them is asymptotically  small, requiring a more sophisticated approach to define the sequence of time scales $\mathfrak{t}_{k}^{\varepsilon} (i)$. \\

There is another feature that distinguishes the set-up in the current paper from the case of a Markov chain with a finite state space. In many applications, the transition probabilities and transition times may depend not on one of finitely many starting points but on a class of equivalence that a starting point belongs to and (to an asymptotically small extent) on the representative in a class of equivalence. Examples of such systems are discussed in Section~\ref{examples}.
Let $S = S_1 \bigcup ... \bigcup S_M$ be a decomposition of a metric space $S$ into a union of disjoint measurable components. Let 
$Q^\varepsilon: S \times \mathcal{B}(S \times [0,\infty)) \rightarrow [0,1]$ be a Markov transition kernel, i.e.,  $Q^\varepsilon(x, \cdot)$ is a probability
measure on $\mathcal{B}(S \times [0,\infty))$ for each $x \in S$, and $Q^\varepsilon(\cdot, A)$ is $\mathcal{B}(S)$-measurable for 
each $A \in \mathcal{B}(S \times [0,\infty))$. Note that the first component in the transition kernel
depends only on the spatial variable. Given an initial point $x \in S$, we can define a discrete-time Markov process
$\xi^{x,\varepsilon}_n = ({\bf X}^{x,\varepsilon}_n, {\bf T}^{x,\varepsilon}_n)$ on $S \times [0,\infty)$ with the kernel $Q^\varepsilon$,
where $\xi^{x,\varepsilon}_0 = ({\bf X}^{x,\varepsilon}_0, {\bf T}^{x,\varepsilon}_0) = (x, 0)$. The components 
${\bf X}^{x,\varepsilon}_n$,  $n \geq 0$,  and $ {\bf T}^{x,\varepsilon}_n$, $n \geq 1$, of this Markov process can be viewed as the states and the inter-arrival times, respectively, of a continuous-time semi-Markov process $X^{x,\varepsilon}_t$. Namely, we define,  for $n \geq 0$, 
\begin{equation} \label{defofSMP}
X^{x,\varepsilon}_t = {\bf X}^{x,\varepsilon}_n~~~{\rm for}~~ \sum_{m =0}^n {\bf T}^{x,\varepsilon}_m \leq t  < 
\sum_{m =0}^{n+1} {\bf T}^{x,\varepsilon}_m. 
\end{equation}
For $x \in S$ and $B \in \mathcal{B}(S)$,  the transition probability associated to the semi-Markov process~is
\[
P^\varepsilon(x,B) = \mathbb{P}({\bf X}^{x,\varepsilon}_1 \in B) = Q^\varepsilon(x, B \times [0,\infty)).
\]
Assuming that $P^\varepsilon(x,B) \neq 0$, the conditional transition time $T^\varepsilon(x,B)$ is defined as a random variable with the distribution 
$$
\mathbb P \left({ T^{\varepsilon}(x, B) \le t  }\right)  = \mathbb{P}({\bf T}^{x,\varepsilon}_1 \leq t| {\bf X}^{x,\varepsilon}_1 \in B) = Q^\varepsilon(x, B \times [0,t])/Q^\varepsilon(x, B \times [0,\infty)).
$$
We assume that 
$P^\varepsilon(x,S_i) = 0$ for each $1 \leq i \leq M$ and $x \in S_i$ (no ``internal" transitions).  Further assumptions, primarily 
generalizing the concept of complete asymptotic regularity, will be stated below. 

The paper is organized as follows. In Section~\ref{assum}, we state the assumptions on the transition probabilities and the
transition times of the semi-Markov process, including the generalization of the notion of complete asymptotic regularity. We conclude this section by stating our main result.  In Section~\ref{examples}, we give several examples of stochastic systems leading to 
parameter-dependent semi-Markov processes that satisfy our assumptions. In Section~\ref{Reduction to a process with nearly independent transition times and transition probabilities}, we show how the problem reduces to considering a system with transition times that are nearly independent of transition probabilities.  In 
Section~\ref{hier}, we collect some results, later to be used recursively, about abstract Markov chains. In Section \ref{hier_main},
we define a hierarchy of Semi-Markov processes on coarser partitions of our state space that capture the dynamics of our process on different time scales. Section~\ref{ergodic_section} is devoted to a description of limiting distributions for parameter-dependent Semi-Markov processes with a particular structure (roughly speaking, the limiting system has one ergodic class and the time under consideration is insufficient to escape from this ergodic class). In Section~\ref{times}, we apply the results of Section \ref{ergodic_section} inductively to different levels of the hierarchy constructed in Section \ref{hier_main} to prove our main result.

\section{Assumptions on the dynamics of the system and formulation of the main result.} \label{assum}

For two functions $f, g: (0, \infty) \rightarrow [0,\infty)$, we will write $f(\varepsilon) \ll g(\varepsilon)$ (as $\varepsilon \downarrow 0$) if $\lim_{\varepsilon \downarrow 0} (f(\varepsilon)/g(\varepsilon)) = 0$. This includes the case when $f \equiv 0$ and $g >0$. For two functions $f, g$ that also depend on a parameter $x$, we will write $ f(x, \varepsilon) \ll g(x, \varepsilon) $ uniformly in $x$ if $\lim_{\varepsilon \downarrow 0}  (f(x, \varepsilon) / g(x, \varepsilon)) = 0$ uniformly in all $x$ being considered. For two functions $f, g: (0, \infty) \rightarrow (0,\infty)$, we will write $f(\varepsilon) \sim g(\varepsilon)$ (as $\varepsilon \downarrow 0$) if $\lim_{\varepsilon \downarrow 0} (f(\varepsilon)/g(\varepsilon)) = 1$. For two functions $f, g$ that also depend on a parameter $x$, we will write $ f(x, \varepsilon) \sim g(x, \varepsilon) $ uniformly in $x$ if $\lim_{\varepsilon \downarrow 0}  (f(x, \varepsilon) / g(x, \varepsilon)) = 1$ uniformly in all $x$ being considered.     

We specify the following assumptions on the process $X^{x, \varepsilon}_t.$ We require the one-step transition probabilities from a point $x \in S_i$ to $S_j$ and the corresponding transition times to depend negligibly on $x$ (if $i$ and $j$ are fixed).  Namely, we assume that,  for every $ i \neq j,  \ P^{\varepsilon} (x, S_j)$ is either  identically $0$ for every $x \in S_i$ or it is positive for all $x \in S_i $ and $\varepsilon > 0,$ and in the latter case there are functions $P_{ij}^{\varepsilon},  \varepsilon > 0, $ such that 
\begin{equation} \label{p_ij^varepsilon_assumption_first}
    \lim \limits_{\varepsilon  \downarrow  0} \dfrac{P^{\varepsilon}(x, S_j)}{P^{\varepsilon}_{ij}} = 1 
\end{equation} 
uniformly in $x \in S_i.$ If $P^{\varepsilon}(x, S_j) \equiv 0$ for $x \in S_i,$ then  we set $P_{ij}^{\varepsilon} \equiv 0.$ In particular, $P_{ii}^{\varepsilon} \equiv 0.$ The functions $P_{ij}^{\varepsilon}$ can be chosen in such a way that $ \sum_{j \neq i} P_{ij}^{\varepsilon} = 1$ for every $i.$  

A condition concerning communication between the sets $S_i$ will be imposed on the functions $P_{ij}^{\varepsilon}.$  We assume that,  for every $\varepsilon > 0$ and every $i \neq j, $ there exists a sequence $i_1, \dots i_k $ in $\{1, \dots M\}$ such that $ i = i_1, j= i_k, $ and $P_{i_ri_{r+1}}^{\varepsilon} > 0$ for every $\varepsilon > 0.$

We assume that for $i \neq j$ there are functions $\tau_{ij}^{\varepsilon},   \varepsilon > 0, $ satisfying 

\begin{equation} \label{time_assumption_first}
    \lim \limits_{\varepsilon  \downarrow  0} \dfrac{\mathbb E \left[{T^{\varepsilon}(x, S_j)}\right]}{\tau^{\varepsilon}_{ij}} = 1 
\end{equation} 
uniformly in $x \in S_i$ provided that $P^{\varepsilon}_{ij}$ is not identically zero. 

We impose a tightness condition. Namely, we assume that
there exist $C, \varepsilon_0 > 0 $ such that 
\begin{equation}
\label{varianceest_new}
{\rm Var} \left({T^{\varepsilon}(x, S_j)}\right) \le C \cdot (\tau^{\varepsilon}_{ij})^2 \end{equation}
 for every $\varepsilon \le \varepsilon_0, x \in S_i$, and $i \neq j$, provided that $P_{ij}^{\varepsilon}$ is not identically zero. 
For any given bounded, continuous function $f:[0, \infty) \to \mathbb R$ and $i, j \in \{1, \dots, M\},$ we assume that \begin{equation} \label{levimetric}
    \lim \limits_{\varepsilon \downarrow 0} \left({ f \left({\dfrac{ T^{\varepsilon}(x_1, S_j) }{ \tau^{\varepsilon}_{ij} }}\right) - f \left({\dfrac{ T^{\varepsilon}(x_2, S_j) }{ \tau^{\varepsilon}_{ij} }}\right)}\right) = 0 \end{equation} uniformly in $x_1, x_2 \in S_i$ provided that $P_{ij}^{\varepsilon}$ is not identically $0$.

\color{black}The next condition helps us avoid purely cyclical behaviour.
Consider a sequence $x_n \in S_i$ with $P^{\varepsilon}_{ij}$ not identically zero, and a sequence $\varepsilon_n$ such that $\varepsilon_n \downarrow 0$ as $n \rightarrow \infty$. We assume that
$$
{\rm if}~~ \frac{T^{\varepsilon_n}(x_n,S_j)}{\tau^{\varepsilon_n}_{ij}} \rightarrow \xi~~{\rm in}~{\rm distribution}~{\rm as}~{n \rightarrow \infty}~{\rm for}~{\rm some}~\xi,~{\rm then}~\xi~
{\rm has}~{\rm no}~{\rm atoms}.
$$

The main assumption on the structure of the process is the following: We will say that the process $X_t^{x, \varepsilon}$ is  completely asymptotically regular if there exists a limit
\begin{equation} \label{carerenewal} \lim \limits_{\varepsilon \downarrow 0} \dfrac{P_{a_1b_1}^{\varepsilon}}{P_{c_1d_1}^{\varepsilon}} \cdot \dfrac{P_{a_2b_2}^{\varepsilon}}{P_{c_2d_2}^{\varepsilon}} \cdots \dfrac{P_{a_rb_r}^{\varepsilon}}{P_{c_rd_r}^{\varepsilon}} \cdot  \dfrac{\tau_{ab}^{\varepsilon}}{\tau_{cd}^{\varepsilon}} \in [0, \infty]
\end{equation}
 for every positive integer $r$ and every $  a,a_i,b,b_i,c,c_i,d,d_i \in \{1,2, \dots M\}$ with $a_i \neq b_i, c_i \neq d_i, a \neq b$ and $c \neq d$ for which the ratios appearing in the limit are defined. \\ \\
\\
\begin{remark}
    It is sufficient to assume the existence of these limits with $a_i = c_i.$ First notice that $P_{ij}^0 = \lim \limits_{\varepsilon \to 0} P_{ij}^{\varepsilon}$ exists for every $i,j $. Indeed, when $P_{ij}^{\varepsilon}$ is not identically zero,    

$$ \lim \limits_{\varepsilon \to 0} P_{ij}^{\varepsilon} = \lim \limits_{\varepsilon \to 0} \dfrac{P_{ij}^{\varepsilon}}{1} =  \lim \limits_{\varepsilon \to 0} \dfrac{P_{ij}^{\varepsilon}}{\displaystyle \sum_{j' } P_{ij'}^{\varepsilon} } = \lim \limits_{\varepsilon \to 0} \  \dfrac{1}{\displaystyle \sum_{j'} \dfrac{P_{ij'}^{\varepsilon}}{P_{ij}^{\varepsilon}} } \in [0, \infty). $$
The limits $P_{ij}^0$ satisfy $\sum_{i \neq j} P_{ij}^0 = 1.$ Hence, for every $i $ there is $j $ such that $P_{ij}^0 > 0$. Suppose we are given  $a,b,c,d \in \{1,2, \dots M\}$. We can find $a', c' \in \{1,2, \dots M\}$ such that $P_{aa'}^0, P_{cc'}^0 > 0$. Then,   $$ \dfrac{P_{ab}^{\varepsilon}}{P_{cd}^{\varepsilon}} \sim  \dfrac{P_{ab}^{\varepsilon}}{P_{aa'}^{\varepsilon}} \cdot \dfrac{P_{cc'}^{\varepsilon}}{P_{cd}^{\varepsilon}} \cdot \dfrac{P_{aa'}^0}{P_{cc'}^0} \;\; \text{as} \;\; \varepsilon \to 0.
$$ We can do this for every factor appearing in $(2)$ and  obtain the following, seemingly weaker but equivalent, condition for complete asymptotic regularity:  
  
\begin{equation} \label{carea}
    \lim \limits_{\varepsilon \downarrow 0} \dfrac{P_{a_1b_1}^{\varepsilon}}{P_{a_1d_1}^{\varepsilon}} \cdot \dfrac{P_{a_2b_2}^{\varepsilon}}{P_{a_2d_2}^{\varepsilon}} \cdots  \dfrac{P_{a_rb_r}^{\varepsilon}}{P_{a_rd_r}^{\varepsilon}} \cdot  \dfrac{\tau_{ab}^{\varepsilon}}{\tau_{cd}^{\varepsilon}} \in [0, \infty]
\end{equation}  
     for every $a_i, b_i, d_i, a,b,c,d \in \{1,2, \dots , M\}$ with $a_i \neq b_i, a_i \neq d_i, a \neq b,$ and $c \neq d$ for which the quantities appearing in the limit exist. \\
\end{remark} 
The following simple lemma (proved in \cite{FK1}) will be useful in applying the property of complete asymptotic regularity. 
 \begin{lemma} \label{appendix_lemma_1}
 Suppose $f_1(\varepsilon), \dots f_m(\varepsilon)$ and $g_1(\varepsilon), \dots , g_{n}(\varepsilon)$ are positive functions which satisfy $$ \lim \limits_{\varepsilon \to 0} \dfrac{f_i(\varepsilon)}{g_{j}(\varepsilon)} \in [0, \infty]  $$ for $ 1 \le i \le m$ and $ 1 \le j \le n.$ Then, 
     $$ \lim \limits_{\varepsilon \to 0} \dfrac{f_1(\varepsilon) + \dots +  f_m(\varepsilon) }{g_1(\varepsilon) + \dots + g_{n}(\varepsilon) } \in [0, \infty]. $$
 \end{lemma}
Next, we formulate the main result on the metastable behavior of the process $X^{x,\varepsilon}_t$ (i.e., on the limiting distribution of the process at different time scales). However, specifying the time scales is accomplished via a construction of the hierarchy below.

We can keep track not only of the current state of the process $X^{x,\varepsilon}_t$ but also of the next state to be visited. Namely, similarly to (\ref{defofSMP}), we define
\[
\hat{X}^{x,\varepsilon}_t = {\bf X}^{x,\varepsilon}_{n+1}~~~{\rm for}~~ \sum_{m =0}^n {\bf T}^{x,\varepsilon}_m \leq t  < 
\sum_{m =0}^{n+1} {\bf T}^{x,\varepsilon}_m. 
\]

\begin{theorem} 
\label{prelimt1} Suppose that the assumptions on the transition probabilities and the transition times stated above are satisfied. 
For each $1 \le i \le M$, there exists a sequence of time scales (functions of $\varepsilon > 0$) satisfying $ 0 = \mathfrak{t}_0^{\varepsilon}(i) \ll \mathfrak{t}_1^{\varepsilon}(i) \ll \mathfrak{t}_2^{\varepsilon}(i) \ll \dots  \ll \mathfrak{t}_{n(i)}^{\varepsilon}(i) = \infty $ and a sequence of probability measures $ \mu_{i, 1}, \mu_{i, 2}, \dots, \mu_{i, n(i)} $ on $\{(j,j'): j, j' \in \{1,...,M\}, j \neq j'\}$ with the following property: if $t(\varepsilon)$ is a function that satisfies $ \mathfrak{t}_{k-1}^{\varepsilon}(i) \ll t(\varepsilon) \ll  \mathfrak{t}_{k}^{\varepsilon}(i)$ for some $k\ge 1$, then $$ \lim \limits_{\varepsilon \to 0} \mathbb P \left({X_{t(\varepsilon)}^{x, \varepsilon} \in S_j, \hat{X}_{t(\varepsilon)}^{x, \varepsilon} \in S_{j'} }\right) = \mu_{i, k}(j, j')  $$ uniformly in $x \in S_i$ for every $j, j' \in \{1, \dots, M\}$, $j \neq j'$.
\end{theorem} 
\noindent
\begin{remark}
 It can be shown (see the remark after Theorem \ref{prelimt2}) that Theorem \ref{prelimt1} could be formulated without the dependence of the time scales on $i$ (by combining the sets of time scales corresponding to different $i$ into one collection of time scales). \end{remark}

It is also interesting to explore what happens if, instead of assuming that $ \mathfrak{t}_{k-1}^{\varepsilon}(i) \ll t(\varepsilon) \ll  \mathfrak{t}_{k}^{\varepsilon}(i)$, we assume that $t(\varepsilon)$ is commensurate with one of the special time scales, i.e., $ { t(\varepsilon) }/{\mathfrak{t}_{k}^{\varepsilon}(i)} \rightarrow  C \in (0, \infty) $ as $\varepsilon \downarrow 0$ for some $k \geq 1$. This question is briefly discussed after Theorem \ref{prelimt2}.

\section{Examples of systems leading to parameter-dependent semi-Markov processes.} \label{examples}

In this section, we discuss several stochastic systems whose long-time behavior can be analyzed by considering an associated semi-Markov process and applying Theorem~\ref{prelimt1}. In each of these examples, the asymptotics has been studied using the particular properties of this or that system; the main point of this section is to show that all these cases fit in a common framework.

\subsection{Random perturbations of dynamical systems with finitely many stable equilibrium points} \label{ex1}

Consider a dynamical system on $\mathbb{R}^d$
\begin{equation} \label{dyns1}
\dot{Y}^{x}_t =  b(Y^{x}_t), ~~Y^{x}_0 = x \in \mathbb{R}^d,
\end{equation}
and a diffusion process that is a small perturbation of the dynamical system
\begin{equation}
\label{diffpr1}
dY^{x, \varepsilon}_t = b(Y^{x, \varepsilon}_t)dt + \varepsilon \sigma(Y^{x, \varepsilon}_t) d W_t, ~~Y^{x, \varepsilon}_0 = x \in \mathbb{R}^d.
\end{equation}

Here, $W_t$ is a $d$-dimensional Wiener process, while $\sigma$  and $b$ are Lipschitz continuous matrix- and vector-valued functions, respectively; $\sigma$ is assumed to be positive-definite for each value of the argument.  If $\varphi: [0.T] \rightarrow \mathbb{R}^d$ 
is a smooth curve with $\varphi(0) = x$, then the asymptotics, as $\varepsilon \downarrow 0$, of the probability that a trajectory of 
$Y^{x,\varepsilon}_t$ on $[0,T]$ belongs to a small neighborhood of $\varphi$ can be expressed in terms of the action functional
\[
S(\varphi) = \frac{1}{2} \int_0^T \langle a(\varphi_t) (\dot{\varphi}(t) - b(\varphi(t))), \dot{\varphi}(t) - b(\varphi(t) \rangle dt,
\]
where $a = (\sigma \sigma^*)^{-1}$.  (See \cite{FW} for the large deviations theory of randomly perturbed dynamical systems). The role of the action functional is that it can be used to define the quasi-potential $V(x,y)$, which, in turn, can be used to express the probabilities of rare transitions (that are due to the small diffusion) between different equilibria of the unperturbed system. Namely,  the quasipotential is defined as
\[
V(x,y) = \inf\{S(\varphi)| \varphi(0) = x, \varphi(T) = y, T \geq 0\}.
\]
Observe that $V(x,y) = 0$ if, for arbitrarily small neighborhoods $U_x$ and $U_y$ of $x$ and $y$, respectively, there is a trajectory of the unperturbed dynamical system going from $U_x$ to $U_y$. Let us assume that there is a finite collection of compact sets $K_1,...,K_M \subset \mathbb{R}^d$ such that

(a) For each $1 \leq i \leq M$, $V(x,y) = V(y,x) = 0$ whenever $x,y \in K_i$.

(b) If $x \in K_i$ and $y \notin K_i$ for some $i$, then $V(x,y) > 0$.

(c) For each $x \in \mathbb{R}^d$, there are $i$ and $y \in K_i$ such that $V(x,y) = 0$.

(d) For some (and therefore for each) $x \in \mathbb{R}^d$, $\lim_{\|y \| \rightarrow \infty} V(x,y) = \infty$.
\\
Let $S_i$ be a $\delta$-neighborhood of $K_i$, where $\delta > 0$ sufficiently small so that $S_i$ are disjoint and, for $x \in S_i$, there is $y \in K_i$ such that $V(x,y) = 0$. Define 
\[
V_{ij} = \inf\{V(x,y)| x\in K_i, y \in K_j\},~~1 \leq i,j \leq M.
\]
According to the Freidlin-Wentzell theory (\cite{FW2}, \cite{FW}), the quantities $V_{ij}$ determine the logarithmic asymptotics of the expected transition times between the sets $S_i$ and $S_j$. For example, let $A^{x,\varepsilon}_j$ be the event that the process $Y^{x,\varepsilon}_t$ 
 starting at $x \in S_i$ hits $S_j$ prior to hitting any $S_k$ with $k \neq i,j$. For $x \in S_i$,  
let $\tau^{x,\varepsilon}$ be the first time when  the process reaches one of the sets $S_j$, other than $S_i$. Assume that $V_{ij} < V_{ik}$ 
for each $k \neq i,j$. It has been shown (see \cite{FW}) that
\begin{equation} \label{logli}
\mathbb{P} (A^{x,\varepsilon}_j) \rightarrow 1~~~{\rm and}~~~ \varepsilon^2\ln \mathbb{E}(\tau^{x,\varepsilon}|A^{x,\varepsilon}_j) \rightarrow  V_{ij}~~~{\rm as}~~\varepsilon \downarrow 0
\end{equation}
uniformly in $x \in S_i$. 

A semi-Markov process $X^{x, \varepsilon}_t$ on the state space $S = S_1 \bigcup ... \bigcup S_M$ can be naturally associated with the diffusion process $Y^{x, \varepsilon}_t$. Namely, assume that $x \in S_i$ for some $i$, and let $0= \tau^{x,\varepsilon}_0,
\tau^{x,\varepsilon}_1,
\tau^{x,\varepsilon}_2,...$ be the times of the successive visits by $Y^{x, \varepsilon}_t$ to different sets $S_i$, i.e., if $Y^{x, \varepsilon}_{\tau^{x,\varepsilon}_k} \in S_i$ with $k \geq 0$, we define
\[
\tau^{x,\varepsilon}_{k+1} = \inf\{t \geq \tau^{x,\varepsilon}_{k}: Y^{x, \varepsilon}_t \in S\setminus S_i\}.
\]
We then define
\[
X^{x, \varepsilon}_t = 
Y^{x,\varepsilon}_{\tau^{x,\varepsilon}_{k}},~~~{\rm for}~~~\tau^{x,\varepsilon}_{k} \leq t < \tau^{x,\varepsilon}_{k+1}.
\]
Verifying, for the process $X^{x,\varepsilon}_t$, all the assumptions on the transition probabilities and transition times outlined in Section~\ref{assum} is not always an easy task since the Freidlin-Wentzell theory provides only the logarithmic asymptotics for the required quantities. Nevertheless, metastability has been described (\cite{F5}, \cite{FW}) for generic randomly perturbed dynamical systems using a construction of the hierarchy of cycles. The main assumption (satisfied by a generic system) is that the notion of the ``next" attractor or the next cycle is correctly defined at each step of the construction. In particular, the assumption (which led to (\ref{logli})) that, for each $i$, the minimum $\min_j V_{ij}$ is achieved for a single value of $j$  is needed to identify the next set visited by the process after it reaches $S_i$.

The precise asymptotics for the expected transition times 
\begin{equation} \label{ekf}
\mathbb{E}(\tau^{x,\varepsilon}|A^{x,\varepsilon}_j) \sim c_{ij} \exp(V_{ij}/\varepsilon^2)
\end{equation}
has been obtained under additional assumptions on the vector field $b$ and the diffusion matrix $\sigma$ and the assumption that the sets $K_i$ are points. For example, for random perturbations that result in reversible diffusion processes, the asymptotics (\ref{ekf}) (the Eyring–Kramers formula) has been established in \cite{BE}, \cite{He} (see also \cite{Ey}, \cite{Kra}, \cite{Land}, \cite{Lang} for earlier work in the physics literature and \cite{Bou} for a generalization to non-reversible processes under certain assumptions on the  quasipotential). The uniqueness of the ``next" attractor is not always necessary to describe the metastable behavior of a randomly 
perturbed dynamical system. In particular, for a class of systems, metastability, including the behavior of the process at the transitional time scales, was recently described in \cite{Lan1}, \cite{Lan2}. 

The connection between the results quoted above and the current paper is that, once the asymptotics (\ref{ekf}) is verified for the underlying diffusion process, it is usually not difficult to show that the assumptions listed in Section~\ref{assum} are satisfied by the corresponding semi-Markov process. Therefore, the metastable behavior of the original process follows from  Theorem~\ref{prelimt1}. 

\subsection{Random perturbations of diffusions with invariant submanifolds}

Here, we briefly recall the results from \cite{FKnew1}, \cite{FKnew2} on metastability for randomly perturbed degenerate diffusions and see how the associated semi-Markov processes fit in the abstract framework.

Consider a diffusion process $Y_t^x$ that satisfies the stochastic differential equation
\[
d Y^x_t = v_0( Y^x_t) dt + \sum_{i=1}^d v_i( Y^x_t) \circ d W^i_t,~~~~~Y^x_0 = x \in \mathbb{R}^d,
\]
where $W^i_t$ are independent Wiener processes and $v_0,...,v_d$ are  bounded $C^3(\mathbb{R}^d)$ vector fields.  
The Stratonovich form  is convenient here since it allows
one to provide a natural coordinate-independent description of the process. The generator of the process is the operator
\[
Lu = L_0 + \frac{1}{2} \sum_{i=1}^d L_i^2,
\]
where $L_i u = \langle v_i, \nabla u \rangle$ is the operator of differentiation along the vector field $v_i$, $i =0,...,d$.   

Let
$K_1, ... , K_M \subset \mathbb{R}^d$ be $C^4$-smooth non-intersecting compact surfaces with dimensions $d_1,...,d_M$, respectively, where $0 \leq d_i < d$, $1 \leq i \leq M$. Points
are considered to be zero-dimensional surfaces. Assume that each of the surfaces is invariant for the process and that the diffusion restricted to a single surface is an ergodic process (the latter condition is trivially satisfied for $K_i$ if $d_i = 0$). Denote 
$D =  \mathbb{R}^d \setminus \left( K_1 \bigcup ... \bigcup K_M \right)$. We assume that the diffusion matrix is non-degenerate on~$D$.  

For each $i = 1,...,M$  and each $x \in K_i$, we define $T(x)$ to be the tangent space to $K_i$ at $x$ (with $T(x)$ being a zero-dimensional space if $d_i = 0$). We
 assume that: 

(a) ${\rm span} (v_0(x),v_1(x),...,v_d(x)) = {\rm span} (v_1(x),...,v_d(x)) = T(x)$ for $x \in K_1 \bigcup ... \bigcup
 K_M$;

(b) ${\rm span} (v_1(x),...,v_d(x)) = \mathbb{R}^d$ for $x \in D$. 
\\
This is just a more convenient (and slightly stronger) way to express the assumptions that we already made: (a)~implies that the surfaces are invariant, and the process is ergodic on each surface; 
(b)~means that the diffusion matrix is non-degenerate on $D$. Under Assumption (a), the set $K_1 \bigcup ... \bigcup
 K_M$ is inaccessible, i.e., the process $X_t$ starting at $x \in D$
does not  reach it in finite time. In fact,
another assumption, concerning the ergodic properties  of the process not only on
the surfaces $K_i$ but in their vicinities and slightly strengthening Assumption (a), is needed (see \cite{FKnew2} for details). 
Under such assumptions,  
each surface can be classified
as either attracting or repelling for $Y^x_t$, depending (roughly speaking) on whether 
$\mathbb{P} (\lim_{t \rightarrow \infty} {\rm dist}(Y^x_t, K_i) = 0) > 0$  for each $x$ sufficiently close to $K_i$ (attracting surface) or 
$\mathbb{P}(\lim_{t \rightarrow \infty} {\rm dist}(Y^x_t, K_i) = 0) = 0$  for each $x \notin K_i$ (repelling surface). 

 Next, we perturb the process $Y^x_t$ by a small non-degenerate diffusion. 
The resulting process $Y^{x,\varepsilon}_t$ satisfies
the following stochastic differential equation:
\[
d Y^{x,\varepsilon}_t = (v_0 + \varepsilon^2   \tilde{v}_0) (Y^{x,\varepsilon}_t ) dt + \sum_{i=1}^d v_i(Y^{x,\varepsilon}_t ) \circ d W^i_t + \varepsilon 
\sum_{i =1}^d \tilde{v}_i (Y^{x,\varepsilon}_t ) \circ d \tilde{W}^i_t,~~~Y^{x,\varepsilon}_t = x \in \mathbb{R}^d,
\]
where $\tilde{W}^i_t$ are independent Wiener processes (also independent of all $W^i_t$), and $\tilde{v}_0,...,\tilde{v}_d$ are 
bounded 
$C^3(\mathbb{R}^d)$ vector fields.   In order to make our assumption on the 
non-degeneracy of the
perturbation more precise, we state it as follows: 

(c) ${\rm span} (\tilde{v}_1(x),...,\tilde{v}_d(x)) = \mathbb{R}^d$ for $x \in \mathbb{R}^d$. 

We also include a condition that ensures that the processes $Y^x_t$ and $Y^{x,\varepsilon}_t$ (with sufficiently small $\varepsilon$) are positive-recurrent. Namely, we assume that:

(d)  $v_0,...,v_d, \tilde{v}_0,...,\tilde{v}_d$ are bounded and $\langle v_0(x), x \rangle < -c \|x\|$ for some $c > 0$ and all sufficiently large
$\|x\|$.  

 The generator of the process $Y^{x,\varepsilon}_t$ is the operator $L^\varepsilon = L + \varepsilon^2 \tilde{L}$ with
\[
\tilde{L}u = \tilde{L}_0 + \frac{1}{2} \sum_{i=1}^d \tilde{L}_i^2,
\]
where $\tilde{L}_i u = \langle \tilde{v}_i, \nabla u \rangle$ is the operator of differentiation along the vector field $\tilde{v}_i$, $i =0,...,d$. 

The process $Y^x_t$  starting at $x$ is approximated well by $Y^{x,\varepsilon}_t$  (with the same starting point)  on finite time intervals as $\varepsilon \downarrow 0$. However, due to the presence of the small non-degenerate component,  
$Y^{x,\varepsilon}_t$ cannot get forever trapped in a vicinity of an attracting surface. The process $Y^{x,\varepsilon}_t$ exhibits metastable behavior: the limiting distribution of $Y^{x,\varepsilon}_{t(\varepsilon)}$, as $(1/\varepsilon, t(\varepsilon))$ approaches infinity, 
 depends  
%
  on the initial point and on the function $t(\varepsilon)$. 

Assume, momentarily, that all the surfaces $K_i$ are attracting and that $d_i < d-1$ for each $i$. In order to relate the situation at hand to the analysis of semi-Markov processes, we again introduce the open sets that contain $K_i$ and that are reached by the perturbed process in finite time with positive probability. Namely, let $S_i^\varepsilon$ be a $\sqrt{\varepsilon}$-neighborhood of $K_i$. This appears 
slightly different from the setup in Section~\ref{assum} since these sets depend on  $\varepsilon$, but each $S_i^\varepsilon$ can be parametrized by a set $S_i$ that is independent of $\varepsilon$.

As in Section~\ref{ex1},
a semi-Markov process $X^{x, \varepsilon}_t$ on the state space  $S^\varepsilon_1 \bigcup ... \bigcup S^\varepsilon_M$ can be naturally associated with the diffusion $Y^{x, \varepsilon}_t$ by stopping the latter process each time it reaches a new set $S^\varepsilon_i$.
Let $A^{x,\varepsilon}_j$ be the event that the process $Y^{x,\varepsilon}_t$ 
 starting at $x \in S^\varepsilon_i$ hits $S^\varepsilon_j$ prior to hitting any $S^\varepsilon_k$ with $k \neq i,j$. For $x \in S^\varepsilon_i$,let $\tau^{x,\varepsilon}$ be the first time when  the process reaches one of the sets $S^\varepsilon_j$, other than $S^\varepsilon_i$. It has been shown (see \cite{FKnew1}, \cite{FKnew2}) \be that there are positive constants $\gamma_i$, $0 \leq i \leq M$, $a_{ij}$, $V_{ij}$, $i \neq j$, such that
\[
\mathbb{P} (A^{x,\varepsilon}_j) \rightarrow a_{ij}~~~{\rm and}~~~ \varepsilon^{\gamma_i} \mathbb{E}(\tau^{x,\varepsilon}|A^{x,\varepsilon}_j) \rightarrow  V_{ij}~~~{\rm as}~~\varepsilon \downarrow 0
\]
uniformly in $x \in S^\varepsilon_i$. Note that, unlike the case of randomly perturbed dynamical systems considered above, the transition times for a randomly perturbed  diffusion process (and, therefore, for the associated semi-Markov process) scale as powers of $\varepsilon$, rather than exponentially. All the assumptions from Section~\ref{assum} placed on the semi-Markov process are easy to verify in this example. 

The situation is slightly more complicated, but still amenable to the same analysis, in the case when some of the surfaces are attracting or have dimension $d_i = d-1$. For example, it may happen that $S^\varepsilon_j$ cannot be directly reached from $S^\varepsilon_i$ if these two sets are separated by a $(d-1)$-dimensional attracting surface, in which case $a_{ij} = 0$. Or, if two sets are separated by a repelling surface, it takes an additional time (that scales as a negative power of $\varepsilon$) to cross a repelling surface.

\subsection{Random perturbations of flows with heteroclinic networks}

Consider again the dynamical system ${Y}^{x}_t$ and its random perturbation $Y^{x, \varepsilon}_t $ 
that satisfy equations (\ref{dyns1}) and
(\ref{diffpr1}), respectively. For simplicity, assume that the dimension of the space is $d = 2$. The Freidlin-Wentzell theory, discussed in Section~\ref{ex1}, concerns transitions between different stable equilibria or limit cycles of the unperturbed system. Now, instead, we will consider
the case where the relevant stationary points are saddles. Namely, assume that the points $A_1,...,A_K$ are stationary for $Y^x_t$, i.e., $b(A_i) = 0$, $1 \leq i \leq K$. Let us linearlize $b$ near each of these points, i.e., let us write $b(x) = B_i (x - A_i) + o(|x-A_i|)$ as $x \rightarrow A_i$, where $B_i$ are  $2 \times 2$ matrices.
Let us assume that each of the matrices $B_i$ has one positive  eigenvalue, $\lambda^+_i$, and one negative eigenvalue, $\lambda^-_i$.

A curve $\gamma$
in $\mathbb{R}^2$ is said to be a heteroclinic orbit connecting  $A_i$ to $A_j$ if $\gamma$ is the image of the mapping $Y^x_\cdot: \mathbb{R} \rightarrow \mathbb{R}^2$, where $Y^x_t$ solves (\ref{dyns1}) with $x \in \gamma$ and $\lim_{t \rightarrow -\infty} Y^x_t = A_i$, $\lim_{t \rightarrow +\infty} Y^x_t = A_j$. Let $\mathcal{N}$ be the heteroclinic network based on the set $\mathcal{A} = \{A_1,...,A_K\}$, i.e., the union of $\mathcal{A}$ and the set of all the heteroclinic orbits connecting points from $\mathcal{A}$. Let us assume that $\mathcal{N}$ is connected and is attracing in the following sense: for each $x \in \mathbb{R}^2$, except a finite number of points, $\lim_{t \rightarrow \infty} {\rm dist} (Y^x_t, \mathcal{N}) = 0$. Whether or not $\mathcal{N}$ attracts solutions of  (\ref{dyns1}) starting in a small neighborhood of $\mathcal{N}$ can be determined from the topology of $\mathcal{N}$ and the set eigenvalues $\lambda^+_i, \lambda^-_i$, $1 \leq i \leq K$ (see \cite{Bakhtin1}, \cite{Bakhtin2}).

Consider now the transitions of the process $Y^{x,\varepsilon}_t$ between different heteroclinic  orbits $\gamma \subset \mathcal{N}$. The transition times are of order at least $O(|\ln(\varepsilon)|)$ due to the logarithmic delay of the process near the saddle points. However, it may take much longer to reach some of the heteroclinic orbits due to the following effect: the process is repelled from $\mathcal{N}$ while passing by a saddle point with $\lambda^+_i > |\lambda^-_i|$, even though $\mathcal{N}$ is attracting as a whole. We will say that a heteroclinic orbit $\gamma_j$ can be easily reached from $\gamma_i$ if there is $x \in \gamma_i$ and a constant $c > 0$ such that $\mathrm{P}(\tau^{x,\varepsilon} \leq c |\ln(\varepsilon)|) \geq 1/c$ for all sufficiently small $\varepsilon$, where  $\tau^{x,\varepsilon}$ is the first time when $Y^{x, \varepsilon}_t$ reaches $\gamma_j$. Two heteroclinic orbits will be said to belong to the same equivalence class if they can be easily reached from each other. A heteroclinic orbit will be said to be transient if another orbit can be reached from it, but the two do not belong to the same equivalence class. 

For each equivalence class that does not contain transient orbits, we can consider the union of the orbits that belong to it, thus obtaining the sets $\Gamma_1,...,\Gamma_M$ corresponding to different equivalence classes. We choose a closed curve $S_i \subset \Gamma_i$, $1 \leq i \leq M$, inside each set $\Gamma_i$ (this is done to avoid atypically fast transitions when the starting point is very close to a saddle point). 

As in Section~\ref{ex1},
a semi-Markov process $X^{x, \varepsilon}_t$ on the state space  $S_1 \bigcup ... \bigcup S_M$ can be naturally associated with the diffusion $Y^{x, \varepsilon}_t$ by stopping the latter process each time it reaches a new set $S_i$.
Let $A^{x,\varepsilon}_j$ be the event that the process $Y^{x,\varepsilon}_t$ 
 starting at $x \in S_i$ hits $S_j$ prior to hitting any $S_k$ with $k \neq i,j$. For $x \in S_i$, let $\tau^{x,\varepsilon}$ be the first time when  the process reaches one of the sets $S_j$, other than $S_i$. It can be shown using the results of \cite{Bakhtin3}
 that there are nonnegative constants $b_{ij}$ and positive constants $\gamma_{ij}$, $a_{ij}$, $V_{ij}$, $i \neq j$, such that
\[
\mathbb{P} (A^{x,\varepsilon}_j) \sim a_{ij}\varepsilon^{b_{ij}}~~~{\rm and}~~~  \mathbb{E}(\tau^{x,\varepsilon}|A^{x,\varepsilon}_j) \sim V_{ij}  |\ln(\varepsilon)| \varepsilon^{-\gamma_{ij}} ~~~{\rm as}~~\varepsilon \downarrow 0
\]
uniformly in $x \in S_i$.

Metastability in noisy heteroclinic networks has been studied in \cite{Bakhtin3}.
The point we wish to make here is that it is possible to fit the analysis of noisy heteroclinic networks in the framework of parameter-dependent semi-Markov renewal processes. The proof that all the assumptions made in Section~\ref{assum} are met in this case is likely to be nontrivial (as is the analysis in \cite{Bakhtin3}).

\section{Reduction to a process with nearly independent transition times and transition probabilities.} \label{Reduction to a process with nearly independent transition times and transition probabilities}

It can be shown that the transition  times of a Markov Renewal Process on a discrete space can be assumed to depend only on the current state and be independent of the next state visited (see for example, \cite{AS}).  This is achieved by extending the state space and studying the process $( {\bf X}_n, {\bf X}_{n+1} ).$ We describe how to extend this construction to our case, where the state space is not necessarily discrete.  

 Let $\overline S = \{ (x, j) : x \in S_i, 1 \le i,j \le M, P_{ij}^{\varepsilon} \not \equiv 0   \}.$ \be Let $\overline Q^{\varepsilon}: \overline S \times  \mathcal B(\overline S) \times  \mathcal B([0, \infty)) \to [0,1]$ be the Markov transition kernel given by $$ \overline Q^{\varepsilon} \left({(x, j), (A \times \{k\}) ,  I}\right) = \dfrac{1}{P^{\varepsilon}(x, S_j)} \int \limits_{A \cap S_j} P^{\varepsilon}(y, S_k) Q^{\varepsilon}(x, dy \times I) $$ when $P^{\varepsilon }(x, S_j) \not \equiv 0$ and $0$ otherwise. Recall that $P^{\varepsilon}(x, A) = Q^{\varepsilon}(x, A \times [0, \infty)).$   
 
 For a given $\overline x = (x, j) \in \overline S,$ let us define the Markov process $\left({ \overline{\bf X}_n^{\overline x, \varepsilon}, \overline{\bf T}_n^{\overline x , \varepsilon}
 }\right)$ with kernel $\overline Q^{\varepsilon} $ and $\left({ \overline{\bf X}_0^{\overline x ,  \varepsilon},  \overline{\bf T}_0^{\overline x,  \varepsilon}
 }\right) = (\overline x,  0).$  Notice that, in the language of the original process,  
\begin{equation} \label{bar_no_bar_kernel_equation}
\begin{aligned}
  \overline Q^{\varepsilon} \left({\overline x, A \times \{k\} ,  I}\right) =&  
 \dfrac{\mathbb P \left({ {\bf X}_1^{x, \varepsilon} 
 \in A \cap S_j, {\bf T}_1^{x, \varepsilon} 
 \in I, 
 {\bf X}_2^{x, \varepsilon} 
 \in  S_k}\right)
}{\mathbb P \left({  {\bf X}_1^{x, \varepsilon} 
 \in  S_j}\right)} \\
 =& \mathbb P \left({ {\bf X}_1^{x, \varepsilon} 
 \in A , {\bf T}_1^{x, \varepsilon} 
 \in I, 
 {\bf X}_2^{x, \varepsilon} 
 \in  S_k \ | \  {\bf X}_1^{x, \varepsilon} 
 \in  S_j}\right)
\end{aligned}
\end{equation} whenever $P^{\varepsilon}(x, S_j) > 0,$ and, therefore, there is a natural correspondence between the finite dimensional distributions of $ \left({ {\bf \overline X }_n^{\overline x , \varepsilon}, {\bf \overline T}_n^{\overline x, \varepsilon}  
 }\right)  $ and $ \left({ {\bf X}_n^{x, \varepsilon}, {\bf T}_n^{x, \varepsilon}  
 }\right) .$ 

Below we will see how the assumptions on the process $ \left({ {\bf X}_n^{x, \varepsilon}, {\bf T}_n^{x, \varepsilon}  
 }\right) $ translate to the process $\left({ {\bf \overline X }_n^{\overline x , \varepsilon}, {\bf \overline T}_n^{\overline x, \varepsilon} 
 }\right) .$ First note that the state space $\overline S$ admits a partition $\{ S_i \times \{j\} : 1 \le i, j \le M, i \neq j\}.$ Let us re-index this partition so that $\{ S_i \times \{j\}  \ : \ 1 \le i, j \le M,  i \neq j \} = \{ \overline S_i : 1 \le i \le \overline M \}.$ As before, we will use $\overline P^{\varepsilon}(\overline x, A) $ to denote $\overline Q^{\varepsilon}(\overline x, A \times [0, \infty)).  $

Let $ 1 \le i, j \le \overline M$ and let  $\overline x \in \overline S_i.$ Then there are $a, b \in \{1, \dots M\}$ and $x \in S_a$ such that  $\overline S_i = S_a \times \{b\} $ and $\overline x = (x, b).$ There are $c, d \in \{1, \dots M\}$ such that $\overline S_j = S_c \times \{d\}.$ 
Then
\begin{equation} \label{extended_space_asymptotics_of_probabilities_1_11}
\begin{aligned}
\overline P^{\varepsilon} \left({ \overline x, \overline S_j}\right) &= \overline Q^{\varepsilon} \left({ \overline x, \overline S_j \times [0, \infty) }\right) \\  &= \dfrac{1 }{P^{\varepsilon}(x, S_b) }   \int \limits_{ S_b \cap S_c} P^{\varepsilon} (x, dy) P^{\varepsilon}(y, S_d) \end{aligned}     \end{equation}  From \eqref{p_ij^varepsilon_assumption_first}, it follows that, when  $b = c, \ \lim \limits_{\varepsilon \to 0}  \dfrac{\overline P^{\varepsilon} \left({ \overline x, \overline S_j}\right)}{ P^{\varepsilon}_{cd} } = 1$ uniformly in $\overline x \in \overline S_i$, while $\overline P^{\varepsilon} \left({ \overline x, \overline S_j}\right) \equiv 0 $ when $b \neq c.$ We conclude that either $\overline P^{\varepsilon} \left({ \overline x, \overline S_j}\right) = 0 $ for every $\overline x \in \overline S_i,$ or $\overline P^{\varepsilon} \left({ \overline x, \overline S_j}\right) > 0 $ for every $\overline x \in \overline S_i$ and there is a positive function $\overline P^{\varepsilon}_{ij}$ such that $$ \lim \limits_{\varepsilon \to 0}  \dfrac{\overline P^{\varepsilon} \left({ \overline x, \overline S_j}\right)}{ \overline P^{\varepsilon}_{ij} } = 1 $$ uniformly in $\overline x \in \overline S_i,$ in which case we take $\overline P^{\varepsilon}_{ij} \equiv P^{\varepsilon}_{cd}.$ Note that the assumption of communication among the functions $P^{\varepsilon}_{ab}$ carries over to the functions $\overline P^{\varepsilon}_{ij}$ naturally. 

Note also, by (\ref{bar_no_bar_kernel_equation}), that if $\overline x = (x, j)$ where $x \in S_i, 1\le i, j \le M, $ then
$$
\mathbb P \left({ \overline {\bf T}_1^{\overline x, \varepsilon} \in I}\right) =  \mathbb P \left({ {\bf T}_1^{ x, \varepsilon} \in I \mid { \bf X}_1^{x, \varepsilon} \in S_j}\right) = \mathbb P \left({ T^{\varepsilon}(x, S_j) \in I }\right) $$
provided that $ P^{\varepsilon}(x, S_j) >0$. 

Adopting notation similar to Section 1, for  $\overline x \in \overline S$  and $B \subseteq \overline S$ such that $\overline 
 P^{\varepsilon}(\overline x, B) > 0, $ let us define $\overline 
 T^\varepsilon(\overline 
 x, B)$ to be a random variable with the distribution $$ \mathbb P \left({ \overline 
 T^\varepsilon(\overline 
 x, B) \in I }\right) = \mathbb P \left({ \overline  {\bf T}_1^{\overline  x, \varepsilon} \in I | \overline {\bf X}_1^{\overline  x, \varepsilon} \in B }\right). $$ $\overline 
 T^{\varepsilon}(\overline  x, B)$ is to be understood as the time in transition from $\overline  x$ to $ \overline {\bf X}_1^{\overline  x, \varepsilon} $ conditioned on $\overline 
 {\bf X}_1^{\overline  x, \varepsilon} \in B.$ If $\overline S_i = S_a \times \{b\}, \overline S_j = S_b \times \{c\},  \overline x = (x, b)$ where $x \in S_a, 1 \le i, j \le \overline M$ and $1 \le a,b,c \le M,$ then,  $$ \begin{aligned}
     \mathbb P \left({ \overline 
 T^\varepsilon(\overline 
 x, B ) \in I }\right) &= \dfrac{\overline Q^{\varepsilon} \left({\overline x, \overline S_j \times I}\right)}{\overline P^{\varepsilon} \left({\overline x, \overline S_j}\right) } = \dfrac{ \dfrac{1}{P^{\varepsilon} \left({ x,  S_b}\right)  } \displaystyle \int_{S_b} P^{\varepsilon} \left({ y, S_c }\right) Q^{\varepsilon} \left({ x, dy \times I }\right)  }{   \dfrac{1}{P^{\varepsilon} \left({ x,  S_b}\right)  } \displaystyle \int_{S_b} P^{\varepsilon} \left({ y, S_c }\right) P^{\varepsilon} \left({ x, dy }\right)   } \\ &= \dfrac{\displaystyle \int_{S_b} P^{\varepsilon} \left({ y, S_c }\right) Q^{\varepsilon} \left({ x, dy \times I }\right) }{\displaystyle \int_{S_b} P^{\varepsilon} \left({ y, S_c }\right) P^{\varepsilon} \left({ x, dy }\right)  }. \end{aligned} $$

The purpose of constructing such a  process on an extended space in the conventional case (discrete space) is to establish independence between ${\bf \overline X}_{n}^{\overline x , \varepsilon}$ and ${\bf \overline T}_{n}^{\overline x , \varepsilon},$ given ${\bf \overline X}_{n-1}^{\overline x , \varepsilon}.$ In our case, we do not retrieve complete independence but instead we have the following results. 

\begin{lemma} \label{fpo} For each $1 \le i \le \overline M$, there is a function $\overline \tau_i^{\varepsilon}$ such that \begin{equation} \label{constructed_time_1_expectation}\lim \limits_{\varepsilon \to 0 } \dfrac{ \mathbb E \left({ \overline T^{\varepsilon}(\overline x, \overline S_j) }\right) }{ \overline \tau_i^{\varepsilon}}  = 1 \end{equation} uniformly in $\overline x \in \overline S_i$ for every $1 \le j \le \overline M$ such that $\overline {P}^{\varepsilon}_{ij}$ is not identically zero. Moreover there are $C, \varepsilon_0 > 0$ such that $$ {\rm Var} \left({\overline{T}^{\varepsilon}(\overline {x}, \overline{S}_j)}\right) \le C \cdot (\overline{\tau}_i^\varepsilon)^2 $$ for every $\varepsilon \le \varepsilon_0$, provided that $\overline {P}^{\varepsilon}_{ij}$ is not identically zero. 
\end{lemma}
\proof  
With $\overline x = (x, b), \overline S_j = S_b \times \{c\},$ 

$$ \begin{aligned}
     \mathbb E \left({ \overline T^{\varepsilon}(\overline x, \overline S_j) }\right) 
  &=  \dfrac{ \displaystyle\int^{\infty} \limits_0 t \  \overline Q^{\varepsilon}(\overline x, \overline S_j \times dt ) }{\overline P^{\varepsilon}(\overline x, \overline S_j  )} = \dfrac{\displaystyle\int^{\infty} \limits_0 t \dfrac{1}{P^{\varepsilon}(x, S_b)} \int \limits_{S_b}P^{\varepsilon}(y, S_c) Q^{\varepsilon}(x, dy \times dt) }{\dfrac{1}{P^{\varepsilon}(x, S_b)} \displaystyle \int \limits_{S_b}P^{\varepsilon}(y, S_c) P^{\varepsilon}(x, dy )} \\
    &\le \dfrac{\sup \limits_{y \in S_b} P^{\varepsilon}(y, S_c)  }{\inf \limits_{y \in S_b} P^{\varepsilon}(y, S_c)} \cdot \dfrac{\displaystyle\int^{\infty} \limits_0 t 
\ Q^{\varepsilon}(x, S_b \times dt) }{P^{\varepsilon}(x, S_b )}.
\end{aligned} $$ The first ratio has the same asymptotics as $P_{bc}^{\varepsilon}/P_{bc}^{\varepsilon} = 1.$ The second ratio is exactly  $\mathbb E \left({T^{ \varepsilon}(x, S_b)}\right)$ (recall the definition of $T^{ \varepsilon}(x, S_b)$ from the original process). Hence the second ratio is asymptotically  equivalent to $\tau_{ab}^{\varepsilon}$ as $\varepsilon \downarrow 0$ due to \eqref{time_assumption_first}. These limits hold uniformly in $x \in S_a.$ The lower estimate is proved in the same way by exchanging $\sup$ and $\inf.$ Since $a, b$ are determined by $\overline x,$ we obtain \eqref{constructed_time_1_expectation} by setting $\overline \tau_i^{\varepsilon } \equiv \tau_{ab}^{\varepsilon}.$ 

Similarly, $$
\begin{aligned}
     \mathbb E \left({ \overline T^{\varepsilon}(\overline x, \overline S_j)^2 }\right) 
  &=  \dfrac{ \displaystyle\int^{\infty} \limits_0 t^2 \  \overline Q^{\varepsilon}(\overline x, \overline S_j \times dt ) }{\overline P^{\varepsilon}(\overline x, \overline S_j  )} = \dfrac{\displaystyle\int^{\infty} \limits_0 t^2 \dfrac{1}{P^{\varepsilon}(x, S_b)} \int \limits_{S_b}P^{\varepsilon}(y, S_c) Q^{\varepsilon}(x, dy \times dt) }{\dfrac{1}{P^{\varepsilon}(x, S_b)} \displaystyle \int \limits_{S_b}P^{\varepsilon}(y, S_c) P^{\varepsilon}(x, dy )} \\
    &\le \dfrac{\sup \limits_{y \in S_b} P^{\varepsilon}(y, S_c)  }{\inf \limits_{y \in S_b} P^{\varepsilon}(y, S_c)} \cdot \dfrac{\displaystyle\int^{\infty} \limits_0 t^2 
\ Q^{\varepsilon}(x, S_b \times dt) }{P^{\varepsilon}(x, S_b )} \\
    &\le \dfrac{\sup \limits_{y \in S_b} P^{\varepsilon}(y, S_c)  }{\inf \limits_{y \in S_b} P^{\varepsilon}(y, S_c)}   \mathbb E \left({ T^{\varepsilon}(x,  S_b)^2 }\right). 
\end{aligned} $$  Due to \eqref{varianceest_new}, $\mathbb E \left({ T^{\varepsilon}(x,  S_b)^2 }\right) \le C \left({\tau_{ab}^{\varepsilon}} \right)^2 =  C \left({\overline {\tau}_{i}^{\varepsilon}} \right)^2$ for some $C > 0$ and small enough $\varepsilon$. This proves our second claim.     \qed
\\

\begin{lemma} \label{lemma_extended_space_bounded_function_uniform_limit_1}
Let  $\overline S_i = S_a \times \{b\}, \overline S_j = S_b \times \{c\}$ for some $ 1 \le i, j \le \overline M, 1 \le a,b,c \le  M.$ For every bounded continuous function $f:[0, \infty) \to \mathbb R, $ $$ \lim \limits_{\varepsilon \downarrow 0} \mathbb E \left({ f \left({ \dfrac{  \overline{ T}^{ \varepsilon} (\overline x, \overline S_{j})  }{\overline \tau_i} }\right) - f \left({ \dfrac{  { T}^{ \varepsilon} (y,  S_{b})  }{ \tau_{ab}} }\right)   }\right) = 0 $$ uniformly in $\overline x \in \overline S_i, y \in S_a.$ \end{lemma} 

\proof Let $x, y \in S_a$ and suppose that $\overline x = (x, b).$ 
Without loss of generality, we can assume that $f \geq 0$. 
$$ \begin{aligned} & \mathbb E  \left({ f \left({ \dfrac{  \overline{ T}^{ \varepsilon} (\overline x, \overline S_{j})  }{\overline \tau_i} }\right) - f \left({ \dfrac{  { T}^{ \varepsilon} (y,  S_{b})  }{ \tau_{ab}} }\right)   }\right)  \\ & = \int_0^{\infty} f\left({\dfrac{t}{\overline \tau_i^{\varepsilon}}}\right) \dfrac{\overline {Q}^{\varepsilon}(\overline x, \overline S_{j} \times dt )}{\overline P ^{\varepsilon}\left({\overline x, \overline S_{j}}\right)} - \int_0^{\infty} f\left({\dfrac{t}{ \tau_{ab}^{\varepsilon}}}\right) \dfrac{ {Q}^{\varepsilon}( y, S_{b} \times dt )}{P ^{\varepsilon}\left({ y, S_{b}}\right)}  \\ &= \int_0^{\infty} f\left({\dfrac{t}{ \tau_{ab}^{\varepsilon}}}\right) \dfrac{\dfrac{  1}{P^{\varepsilon}(x, S_{b})} \displaystyle \int  \limits_{S_{b}}  P^{\varepsilon}(x_1, S_{c}) {Q}^{\varepsilon}( x, dx_1 \times dt )}{\dfrac{  1}{P^{\varepsilon}(x, S_{b})} \displaystyle \int  \limits_{S_{b}}  P^{\varepsilon}(x_1, S_{c}) {P}^{\varepsilon}( x, dx_1 )} - \int_0^{\infty} f\left({\dfrac{t}{ \tau_{ab}^{\varepsilon}}}\right) \dfrac{ {Q}^{\varepsilon}( y, S_{b} \times dt )}{P ^{\varepsilon}\left({ y, S_{b}}\right)} \\ &\le \dfrac{ \sup \limits_{x_1 \in S_b} P^{\varepsilon}(x_1, S_{c}) }{ \inf \limits_{x_1 \in S_b} P^{\varepsilon}(x_1, S_{c})} \int_0^{\infty} f\left({\dfrac{t}{ \tau_{ab}^{\varepsilon}}}\right) \dfrac{ {Q}^{\varepsilon}( x, S_{b} \times dt )}{P ^{\varepsilon}\left({ x, S_{b}}\right)} - \int_0^{\infty} f\left({\dfrac{t}{ \tau_{ab}^{\varepsilon}}}\right) \dfrac{ {Q}^{\varepsilon}( y, S_{b} \times dt )}{P ^{\varepsilon}\left({ y, S_{b}}\right)}  \to 0 \end{aligned} $$  uniformly in $x, y \in S_a$ due to \eqref{levimetric}. A lower estimate is obtained similarly by exchanging $\sup$ and $\inf$ above. 
\qed \\

\begin{corollary}  \label{lemma_extended_space_bounded_function_uniform_limit_1_corollary}   Let $i \in \{1, \dots, \overline M\}.$ For every bounded continuous function $f:[0, \infty) \to \mathbb R, $ $$ \lim \limits_{\varepsilon \downarrow 0} \mathbb E \left({ f\left({ \dfrac{ \overline{ T}^{ \varepsilon} (\overline x_1, \overline S_{j_1}) }{\overline \tau_i^{\varepsilon}} }\right) - f\left({  \dfrac{ \overline{ T}^{ \varepsilon} (\overline x_2, \overline S_{j_2}) }{\overline \tau_i^{\varepsilon}} }\right) }\right) = 0 $$ uniformly in $\overline x_1, \overline x_2 \in \overline S_i$ and $j_1, j_2 \in \{1, \dots, \overline M\}, $ provided that $\overline P_{ij_1}^{\varepsilon}, \overline P_{ij_2}^{\varepsilon} > 0.$
\end{corollary} \proof The result holds since the limit in Lemma \ref{lemma_extended_space_bounded_function_uniform_limit_1} is independent of $\overline S_j.$ \qed   
\\

From this point on, we 
will work with the process $\left({ \overline{\bf X}_n^{\overline x, \varepsilon}, \overline{\bf T}_n^{\overline x , \varepsilon}
 }\right)$. 
 However, since it is inconvenient to use the notation with the bars throughout the paper, we drop those bars from the notation (in particular, we'll write $S$ instead of $\overline{S}$ and $M$ instead of $\overline{M}$).
 In order to distinguish the resulting process from the original one, we restate the properties of the auxiliary process that follow from the assumptions we made on the original one.

(a) Due to \eqref{extended_space_asymptotics_of_probabilities_1_11}, for every $ i \neq j,  \ P^{\varepsilon} (x, S_j)$ is either  identically $0$ for every $x \in S_i$ or it is positive for all $x \in S_i $ and $\varepsilon > 0,$ and in the latter case there are functions $P_{ij}^{\varepsilon},  \varepsilon > 0, $ such that,  
$$ \lim \limits_{\varepsilon  \downarrow  0} \dfrac{P^{\varepsilon}(x, S_j)}{P^{\varepsilon}_{ij}} = 1 $$ 
uniformly in $x \in S_i.$ If $P^{\varepsilon}(x, S_j) \equiv 0$ for $x \in S_i,$ then  we set $P_{ij}^{\varepsilon} \equiv 0.$ In particular, $P_{ii}^{\varepsilon} \equiv 0.$ The functions $P_{ij}^{\varepsilon}$ can be chosen in such a way that $ \sum_{j \neq i} P_{ij}^{\varepsilon} = 1$ for every $i.$ For every $\varepsilon > 0$ and every $i \neq j, $ there exists a sequence $i_1, \dots i_k $ in $\{1, \dots M\}$ such that $ i = i_1, j= i_k, $ and $P_{i_ri_{r+1}}^{\varepsilon} > 0$ for every $\varepsilon > 0.$
 \color{black} 

(b) (follows from Lemma~\ref{fpo})   
For every $1 \le i \le M,$ there is a function $\tau_i^{\varepsilon}$ such that \begin{equation} \label{time_assumption_first_extended_space} \lim \limits_{\varepsilon \to 0 } \dfrac{ \mathbb E \left({ T^{\varepsilon}(x, S_j) }\right) }{ \tau_i^{\varepsilon}}  = 1 
 \end{equation} uniformly in $ x \in S_i$ for every $j$ such that $\mathbb P \left({ {\bf X}^{ x, \varepsilon}_{1} \in S_j }\right) = P^{\varepsilon} \left({ x, S_j}\right) 
 > 0 $ and for every $ x \in S_i.$ 

(c) (follows from Lemma \ref{fpo}) There are $C, \varepsilon_0 > 0 $ such that 
\begin{equation}
\label{varianceest_extended_space}
{\rm Var} \left({T^{\varepsilon}(x, S_j)}\right) \le C \cdot (\tau_i^\varepsilon)^2 \end{equation} for every $x \in S_i$ and $i \neq j, $ provided that $P_{ij}^{\varepsilon}$ is not identically $0$.

(d) (follows from Corollary \ref{lemma_extended_space_bounded_function_uniform_limit_1_corollary}) For any given bounded continuous function $f: [0, \infty) \to \mathbb R$ and $i, j_1, j_2 \in \{ 1, \dots, M\}, $ 

$$ \lim \limits_{\varepsilon \downarrow 0}
\mathbb{E}\left({f\left({\frac{{T}^{\varepsilon}({x_1}, {S}_{j_1})}{{\tau}^\varepsilon_i}}\right) - 
f\left({\frac{{T}^{\varepsilon}({x_2}, {S}_{j_2})}{{\tau}^\varepsilon_i}}\right)}\right) = 0 
$$ uniformly in $x_1, x_2 \in S_i$ provided that $P^{\varepsilon}_{ij_1}, P^{\varepsilon}_{ij_2}$ are positive functions. \be

(e) (follows from Lemma \ref{lemma_extended_space_bounded_function_uniform_limit_1}) For a sequence $x_n \in S_i$ and a sequence $\varepsilon_n$ such that $\varepsilon_n \downarrow 0$ as $n \rightarrow \infty,$
$$
{\rm if}~~\frac{T^{\varepsilon_n}(x_n, S_j )}{\tau^{\varepsilon_n}_{i}} \rightarrow \xi~~{\rm in}~{\rm distribution}~{\rm as}~{n \rightarrow \infty}~{\rm for}~{\rm some}~\xi,~{\rm then}~\xi~
{\rm has}~{\rm no}~{\rm atoms}.
$$

\be 

(f) Complete asymptotic regularity for the new process (which follows from assumption \eqref{carea}):   
$$ \lim \limits_{\varepsilon \downarrow 0} \dfrac{P_{a_1b_1}^{\varepsilon}}{P_{a_1d_1}^{\varepsilon}} \cdot \dfrac{P_{a_2b_2}^{\varepsilon}}{P_{a_2d_2}^{\varepsilon}} \cdots  \dfrac{P_{a_rb_r}^{\varepsilon}}{P_{a_rd_r}^{\varepsilon}} \cdot  \dfrac{\tau_{a}^{\varepsilon}}{\tau_{c}^{\varepsilon}} \in [0, \infty] $$
     for every $a_i, b_i, d_i, a,c \in \{1,2, \dots , M\}$ with $a_i \neq b_i, a_i \neq d_i $ for which the quantities appearing in the limit exist. 

The new extended process ($\left({ \overline{\bf X}_n^{\overline x, \varepsilon}, \overline{\bf T}_n^{\overline x , \varepsilon}
 }\right)$ in the old notation, but with the bars dropped) already contains information about the current state as well as the next state of the original process. Therefore, Theorem~\ref{prelimt1} is implied by the following.
\begin{theorem} 
\label{prelimt2} Suppose that conditions (a)-(f) are satisfied. Then,
for each $1 \le i \le M$, there exists a sequence of time scales (functions of $\varepsilon > 0$) satisfying $ 0 = \mathfrak{t}_0^{\varepsilon}(i) \ll \mathfrak{t}_1^{\varepsilon}(i) \ll \mathfrak{t}_2^{\varepsilon}(i) \ll \dots  \ll \mathfrak{t}_{n(i)}^{\varepsilon}(i) = \infty $ and a sequence of probability measures $ \mu_{i, 1}, \mu_{i, 2}, \dots, \mu_{i, n(i)} $ on $ \{1,...,M\}$ with the following property: if $t(\varepsilon)$ is a function that satisfies $ \mathfrak{t}_{k-1}^{\varepsilon}(i) \ll t(\varepsilon) \ll  \mathfrak{t}_{k}^{\varepsilon}(i)$ for some $k \ge 1$, then $$ \lim \limits_{\varepsilon \to 0} \mathbb P \left({X_{t(\varepsilon)}^{x, \varepsilon} \in S_j }\right) = \mu_{i, k}(j)  $$ uniformly in $x \in S_i$ for every $j \in \{1, \dots, M\}$.
\end{theorem} 
\noindent
\begin{remark}
 It will shown in the proof of the theorem that the time scales corresponding to different $i$ can be assumed to be comparable, i.e., there are  limits $\lim_{\varepsilon \downarrow 0} (\mathbf{t}^\varepsilon_{k_1}(i_1)/\mathbf{t}^\varepsilon_{k_2}(i_2)) \in [0,\infty]$.
We will call two time scales equivalent if this limit is positive and finite.
  Thus,  by collecting the time scales for each $1 \le i \le M$ and retaining only one from each set of equivalent scales, we get a set of time scales $  0 = \mathfrak{t}_0^{\varepsilon} \ll \mathfrak{t}_1^{\varepsilon} \ll \mathfrak{t}_2^{\varepsilon} \ll \dots  \ll \mathfrak{t}_{n}^{\varepsilon} = \infty $ and a collection of probability measures $\left\{{\mu_{i, k} : 1 \le k \le n, 1 \le i \le M }\right\}$ on $\{1, \dots, M\}$ such that if $t(\varepsilon)$ is a function that satisfies $ \mathfrak{t}_{k-1}^{\varepsilon} \ll t(\varepsilon) \ll  \mathfrak{t}_{k}^{\varepsilon}$ for some $k \ge 1$, then $$ \lim \limits_{\varepsilon \to 0} \mathbb P \left({X_{t(\varepsilon)}^{x, \varepsilon} \in S_j }\right) = \mu_{i, k}(j)  $$ uniformly in $x \in S_i$ for every $i, j \in \{1, \dots, M\}$. Observe that the indexing of the measures $\mu_{i, k}$ here may be different from that in the theorem since it may happen that $n > n(i)$ for a given $i$. \end{remark} 

Theorem~\ref{prelimt2} includes the assumption that the time scale $t(\varepsilon)$ where the system is observed is deep inside an  interval between a pair of special time scales, i.e., the relation 
$ \mathfrak{t}_{k-1}^{\varepsilon}(i) \ll t(\varepsilon) \ll  \mathfrak{t}_{k}^{\varepsilon}(i)$ holds. It can be shown that, under additional assumptions, a similar conclusion holds even for time scales that satisfy $\lim_{\varepsilon \downarrow 0} t(\varepsilon)/ \mathfrak{t}_{k}^{\varepsilon}(i)  \in (0, \infty)$  for some $k \geq 1$. We will see later (in Section \ref{times}) that the time scales $\mathfrak{t}_{k}^{\varepsilon}(i)$ correspond to expected exit times from certain sets that are of the form $\bigcup_{j \in J_k} S_j$ for some index set $J_k \subseteq \{1, 2, \dots, M\}$ to be specified. In particular, we will take $\mathfrak{t}_{1}^{\varepsilon}(i) = \tau_i^{\varepsilon}$. For those situations where $J_k$ consists of more than one element, it is shown in Theorem \ref{time_exponential_rank_r} that the exit time from $\bigcup_{j \in J_k} S_j$ has an almost exponential distribution, and this fact can be used to study the metastable distributions at time scales commensurate with $\mathfrak{t}_{k}^{\varepsilon}(i)$. However, without making further assumptions, we can't say anything similar about the limiting exit time from  $\bigcup_{j \in J_k} S_j$ if $J_k$ consists of only only one element. 

The following proposition is stated here without a proof since it's not the main focus of the current work.

\begin{proposition}
Suppose that conditions (a)-(f) are satisfied.
Let $x \in S_i$ where $1 \le i \le M.$ Consider the associated  sequence of time scales  $0  =  \mathfrak{t}_0^{\varepsilon}(i) \ll \mathfrak{t}_1^{\varepsilon}(i) \ll \mathfrak{t}_2^{\varepsilon}(i) \ll \dots  \ll \mathfrak{t}_{n(i)}^{\varepsilon}(i) = \infty $ introduced in Theorem~\ref{prelimt2}. 
Suppose that either

(1) For each $i$, $T^\varepsilon(x, S_j)/\tau^\varepsilon_i$ converges weakly for some $x \in S_i$ and some $j$ (the limiting distribution does not have atoms due to (e) and does not depend on $x$ or $j$ due to (d)), or

(2)  For each $i$, 
 $\lim \limits_{\varepsilon \to 0 } ({ t(\varepsilon) }/\tau^\varepsilon_i) \in \{0, \infty\} $.
\\
Suppose that the function
$t(\varepsilon)$  satisfies, $$\lim \limits_{\varepsilon \to 0 } \dfrac{ t(\varepsilon) }{\mathfrak{t}_{k}^{\varepsilon}(i)} = C \in (0, \infty) $$ for some $k \geq 1$ (this is only possible with $k \geq 2$ if the second assumption holds). Then there exists a probability measure $\mu_{i,k,C}$ such that 
\begin{equation}
\label{limmeas}
\lim \limits_{\varepsilon \to 0} \mathbb P \left({X_{t(\varepsilon)}^{x, \varepsilon} \in S_j }\right) = \mu_{i,k,C}(j)  
\end{equation}
uniformly in $x \in S_i$ for every $j \in \{1, \dots, M\}.$
\end{proposition} 

Let us note that (\ref{limmeas}) was proved in \cite{Lan1} and \cite{Lan2} for a certain class of  randomly perturbed dynamical systems (which can be associated with parameter-dependent semi-Markov processes, as discussed in Section~\ref{examples}).

\section{Asymptotic results for discrete-time Markov chains} \label{hier}

In this section, we are only interested in transition probabilities rather than transition times. We start by discussing a more general setup, but will retain some of the same notation. 

Consider a family of Markov chains $ {\bf X}_{n}^{x, \varepsilon}$ on a metric space $S$ given by the Markov kernel $P^{\varepsilon}(x, A), \ x \in S, A \in \mathcal{B}(S) $ - note again that this chain may be distinct from the one related to the semi Markov process $({\bf X}^{x,\varepsilon}_n, {\bf T}^{x,\varepsilon}_n)$ defined above. We will make the following assumptions. 
\\

(A) There is a partition of the state space $S = S_1 \cup S_2 \cup \dots S_M.$ We will assume that $P^{\varepsilon}(x, S_i) = 0$ for every $x \in S_i$ and $ \varepsilon > 0.$  
\\

(B) For every $ i \neq  j,  \ P^{\varepsilon} (x, S_j)$ is either  identically $0$ for every $x \in S_i$ or it is positive for all $x \in S_i $ and $\varepsilon > 0,$ and in the latter case there are functions $P_{ij}^{\varepsilon}, \  \varepsilon > 0, $ such that  
\begin{equation}\label{asymptotic_discrete} \lim \limits_{\varepsilon \downarrow 0} \dfrac{P^{\varepsilon}(x, S_j)}{{P}_{ij}^{\varepsilon}} = 1 \end{equation} uniformly in $x \in S_i.$  In the case when $P^{\varepsilon}(x, S_j)$ is identically zero, we set $P^{\varepsilon}_{ij}$ to also be identically zero. For every $\varepsilon > 0$ and every $i \neq j, $ there exists a sequence $i_1, \dots i_k $ in $\{1, \dots M\}$ such that $ i = i_1, j= i_k, $ and $P_{i_ri_{r+1}}^{\varepsilon} > 0$ for every $\varepsilon > 0.$ 
\\

(C)  For every $i \neq j, i' \neq j',$ $$
  \lim \limits_{\varepsilon \to 0} \dfrac{P^{\varepsilon}_{ij}}{P^{\varepsilon}_{i'j'}} \in [0,  \infty] $$ whenever the denominator is a function that is not identically zero.  \\ 
  \\Assumption (C), referred to as asymptotic regularity, which is weaker than the assumption of complete asymptotic regularity made in \eqref{carerenewal},  ensures that the limits $P_{ij}^0 = \lim \limits_{\varepsilon \to 0} P_{ij}^{\varepsilon} $ exist for every $i ,  j.$ 
\\

(D) Suppose that there is a partition of $\{1,2, \dots, M\}$ into $E \cup F.$ We will assume that $E$ forms one ergodic class with respect to the transition probabilities ${P}_{ij}^{0}$ (this implies that $E$ has more than one element since $P_{ii}^{0} = 0$ for each $i$).  The process on $E$ with transition probabilities ${P}_{ij}^{0}$ and initial point $i$ will be denoted by ${\bf X}_{n}^{i, 0}$. We also assume that $  \{1, 2, \dots M\} $ also forms one ergodic class with respect to  $P^{\varepsilon}_{ij}$ for every $\varepsilon > 0.$ \color{black} 
\\

For $i \in E$ and $ x  \in \bigcup \limits_{j \in E} S_j,$ let $$ \theta^{x, \varepsilon}  = \min \left\{{ n :   {\bf X}^{x,\varepsilon}_n  \in \bigcup_{k \in F} S_k} \right\}, \ \ \ {N}^{\varepsilon}(x, i) = \displaystyle \sum_{n = 0}^{\theta^{x, \varepsilon} } \chi \left({ {\bf X}_{n}^{x, \varepsilon} \in S_i }\right). $$ 
Thus 
$N^{\varepsilon}(x, i)$ is the number of visits  to $S_i$ by ${\bf X}_n^{x, \varepsilon}$ before this process leaves $\bigcup \limits_{j \in E} S_j$.  Let $\lambda(\cdot)$ be the invariant measure associated with the chain on $E$ governed by the transition probabilities $P_{ij}^0.$ We will prove Lemmas \ref{abstract_chain_lemma} and \ref{abstract_chain_lemma_2} under this setting.

\begin{lemma} \label{abstract_chain_lemma}
Suppose that the  assumptions (A)-(D) hold and let $i \in E$. Then   $$ \lim \limits_{\varepsilon \downarrow 0} \dfrac{ \mathbb P \left({\displaystyle \bigcap \limits_{n = 1}^{\theta^{x, \varepsilon}} ( {\bf X}_{n}^{x, \varepsilon} \not \in S_i  ) }\right)}{\dfrac{1}{\lambda(i)} \displaystyle \sum \limits_{\substack{j \in E, \\ k \in F}}\lambda(j) P^{\varepsilon}_{jk} } = 1   $$ uniformly in $ x \in S_i $. Moreover,  $$ \lim \limits_{\varepsilon \to 0}  \mathbb E \left({N^{\varepsilon}(x, i)}\right) \cdot\dfrac{1}{\lambda(i)} \displaystyle \sum \limits_{\substack{j \in E, \\ k \in F}}\lambda(j) P^{\varepsilon}_{jk}  = 1  $$  uniformly in $x \in \bigcup \limits_{j \in E} S_j$.

\end{lemma}

\proof Let $D = E \setminus \{i\}.$ Let us introduce the following notation: $S^E = \bigcup \limits_{j \in E} S_j, \ S^F = \bigcup \limits_{j \in F} S_j $ and $S^D = \bigcup \limits_{j \in D} S_j.$ Let $ \sigma_i^{x, \varepsilon}  = \min \left\{{ \theta^{x, \varepsilon}, \min \{n > 0:  {\bf X}^{x, \varepsilon}_{n} \in S_i\}}\right\}.$  Consider the kernel $\Gamma_i^{\varepsilon}(x, A) = \mathbb P \left({{\bf X}^{x, \varepsilon}_{\sigma_{i}^{x, \varepsilon}}\in A}\right) $ defined for $x \in S^E$ and $A \in \mathcal B(S)$. We will sometimes avoid mentioning the dependence of $\sigma_{i}^{x, \varepsilon}$ on $x, \varepsilon$ below to minimize notation. Then  for $x \in S_{i},$
 $$  \mathbb P \left({\displaystyle \bigcap \limits_{n = 1}^{\theta^{x, \varepsilon}} ( {\bf X}_{n}^{x, \varepsilon} \not \in S_i  ) }\right) =  \Gamma_i^{\varepsilon} \left({x, S^F}\right) = \mathbb P \left({ {\bf X}^{x, \varepsilon}_{\sigma_{i}} \in S^F }\right) = \sum_{j \in E } \mathbb P \left({ {\bf X}^{x, \varepsilon}_{\sigma_{i}} \in S^F, {\bf X}^{x, \varepsilon}_{\sigma_{i} - 1 } \in S_j  }\right).$$  Let $ j \neq i .$  Let $\mathcal A(m, n)$ be the set of all paths $(i_1, i_2, \dots , i_n)$ of  length $n$ such that each $i_{\alpha} \in E \setminus \{i\}, i_n = j$, and $i_\alpha = j$ for exactly $m$ different indices $\alpha.$ 

$$ \begin{aligned}
     \mathbb P &\left({ {\bf X}^{x, \varepsilon}_{\sigma_{i}} \in S^F, {\bf X}^{x, \varepsilon}_{\sigma_{i} - 1 } \in S_j  }\right) \\ &= \sum_{m = 1 }^{\infty} \mathbb P \left({ {\bf X}^{x, \varepsilon}_{\sigma_{i}} \in S^F, \ {\bf X}^{x, \varepsilon}_{\sigma_{i} - 1 } \in S_j, S_j~~\text{was}~~\text{visited}~~m~~\text{times}}\right) \\ &= \sum_{m = 1 }^{\infty} \  \sum_{n = 1}^{\infty } 
 \ \sum_{\left({i_{\alpha} }\right) \in \mathcal A (m, n)} \mathbb P \left({ {\bf X}^{x, \varepsilon}_{1} \in S_{i_1}, \dots, {\bf X}^{x, \varepsilon}_{n} \in S_{i_n}, {\bf X}^{x, \varepsilon}_{n +1} \in S^F  }\right) \\ &\le   \sup_{ y \in S_j } P^{\varepsilon} \left({y, S^F}\right) \sum_{m = 1 }^{\infty} \  \sum_{n = 1}^{\infty } 
 \ \sum_{\left({i_{\alpha} }\right) \in \mathcal A (m,n)} \mathbb P \left({ {\bf X}^{x, \varepsilon}_{1} \in S_{i_1}, \dots, {\bf X}^{x, \varepsilon}_{n} \in S_{i_n}  }\right). \end{aligned} $$  Due to \eqref{asymptotic_discrete} and the fact that $  {\bf X}^{i, 0}_{n} $ is ergodic on $E,$ the terms in the above sum decay exponentially fast with respect to $m$ and $n$, uniformly in $x \in S_i$ and small enough $\varepsilon.$ Since $ P^{\varepsilon} \left({y, S^F}\right) \sim \displaystyle \sum_{k \in F} P_{jk}^{\varepsilon}$ uniformly in $y \in S_j,$ we obtain from the above inequality that $$ \begin{aligned}  \limsup \limits_{\varepsilon \downarrow 0}  \ &\dfrac{ \sup \limits_{x \in S_{i}} \mathbb P \left({ {\bf X}^{x, \varepsilon}_{\sigma_{i}} \in S^F, {\bf X}^{x, \varepsilon}_{\sigma_{i} - 1 } \in S_j  }\right) }{ \displaystyle \sum_{k \in F} P_{jk}^{\varepsilon}} \\ &\le  \sum_{m = 1 }^{\infty} \sum_{n = 1 }^{\infty}   \  \sum_{\left({i_{\alpha} }\right) \in \mathcal A (m, n)} \mathbb P \left({ {\bf X}^{i, 0}_{1} = {i_1}, \dots, {\bf X}^{i, 0}_{n} = {i_n}  }\right) \\ &= \sum_{m = 1 }^{\infty} \mathbb P \left({ \left({ {\bf X}^{i, 0}_{n} }\right)~~\text{visited}~~j~~\text{at}~~\text{least}~~m~~\text{times}~~\text{before}~~\text{returning}~~\text{to}~~i}\right) \\ &= \mathbb E \left({ \text{Number}~~\text{of}~~\text{visits}~~\text{to}~~j~~\text{between}~~\text{visits}~~\text{to}~~i~~\text{by}~~\left({ {\bf X}^{i, 0}_{n} }\right) }\right) \\ &= \dfrac{\lambda(j)}{\lambda(i)}.
    \end{aligned} $$ By \eqref{asymptotic_discrete}, $   \mathbb P \left({ {\bf X}^{x, \varepsilon}_{\sigma_{i}} \in S^F, {\bf X}^{x, \varepsilon}_{\sigma_{i} - 1 } \in S_i  }\right) = P^{\varepsilon} \left( {x, S^F} \right) \sim \displaystyle \sum_{k \in F} P_{jk}^{\varepsilon},  $ uniformly in $x \in S_i$, and therefore the inequality above holds even when $j = i.$ This shows that $$ 
 \limsup \limits_{\varepsilon \downarrow 0} \ \dfrac{ \sup \limits_{x \in S_{i} } \Gamma_i^{\varepsilon} \left({x, S^F}\right) }{ \dfrac{1}{\lambda(i)} \displaystyle \sum_{j \in E} \lambda(j) \displaystyle \sum_{k \in F} P_{jk}^{\varepsilon} } \le 1.$$ We can obtain a lower estimate similarly to show that \begin{equation} \label{abstract_chain_asymptote_lemma_initial_formula_1}
     \Gamma_i^{\varepsilon} \left({x, S^F}\right) 
 \sim  \dfrac{1}{\lambda(i)} \displaystyle \sum_{j \in E} \lambda(j) \displaystyle \sum_{k \in F} P_{jk}^{\varepsilon} 
 \end{equation} uniformly in $x \in S_{i}.$ This proves the first claim in the Lemma. 

 Now let $x \in S^E \setminus S_i.$  Then for $ n \ge 1, $   $$ \mathbb P \left({N^{\varepsilon}(x, i) \ge n}\right)=  \int_{S_i} \cdots  \int_{S_i} \Gamma_i^{\varepsilon}(x, dx_1) \Gamma_i^{\varepsilon}(x_1, dx_2) \dots \Gamma_i^{\varepsilon}(x_{n-1}, dx_n ).$$ Therefore, $$\Gamma_i^{\varepsilon}(x, S_i) \left({ \inf \limits_{y \in S_i} \Gamma_i^{\varepsilon}(y, S_i)}\right)^{n-1}  \le  \mathbb P \left({N^{\varepsilon}(x, i) \ge n}\right) \le \Gamma_i^{\varepsilon}(x, S_i) \left({ \sup \limits_{y \in S_i} \Gamma_i^{\varepsilon}(y, S_i)}\right)^{n-1} .
$$ Since $ \mathbb{E}N^{\varepsilon}(x, i)  = \displaystyle \sum_{n=1}^{\infty }  \mathbb P \left({N^{\varepsilon}(x, i) \ge n}\right),   $ we have  

$$
\dfrac{\Gamma_i^{\varepsilon}(x, S_i)}{ \sup \limits_{y \in S_i} \Gamma_i^{\varepsilon}(y, S_i)}  \le   \mathbb{E}N^{\varepsilon}(x, i)  \le \dfrac{\Gamma_i^{\varepsilon}(x, S_i)}{ \inf \limits_{y \in S_i} \Gamma_i^{\varepsilon}(y, S_i)}.  
$$

The result is proved  for $x \in S^E \setminus S_i$ after noting that $\Gamma_i^{\varepsilon}(x, S_i) \to 1$ as $\varepsilon \downarrow 0$ uniformly in $x \in S^E$. For $x \in S_i$, we simply note that the probability of leaving $S^E$ in one step is asymptotically small, use the Markov property on the location of the first transition in $S^E$, and observe that adding one to $\mathbb EN^{\varepsilon}(x, i)$ does not change its asymptotics.    \qed
\\

\begin{lemma}  \label{abstract_chain_lemma_2}  
 Suppose that assumptions (A)-(D) hold. Let $G \subseteq F, G \neq \emptyset.$ Let $S^G = \bigcup \limits_{k \in G} S_k.$   Then  $$  \mathbb P \left({{\bf X}^{x, \varepsilon}_{\theta^{x, \varepsilon}} \in S^G}\right) \sim   \displaystyle \sum_{k \in G} \displaystyle \sum_{j \in E} \lambda(j) P_{jk}^{\varepsilon}\displaystyle  \text{\Large /} \sum_{k \in F} \sum_{j \in E} \lambda(j) P_{jk}^{\varepsilon} ~~as~~\varepsilon \rightarrow 0$$
uniformly in $x \in \bigcup \limits_{j \in E} S_j.$ 
\end{lemma} 
 \proof Fix an arbitrary $i \in E.$ Note that by replacing $F$ with $G$ in all the arguments  in the proof of  \eqref{abstract_chain_asymptote_lemma_initial_formula_1},  it can be shown that $\Gamma_i^{\varepsilon} \left({x, S^G}\right) 
 \sim  \dfrac{1}{\lambda(i)} \displaystyle \sum_{j \in E} \lambda(j) \displaystyle \sum_{k \in G} P_{jk}^{\varepsilon} $ uniformly in $x \in S_{i},$ where the definitions for $\Gamma_i^{\varepsilon}$ and $\sigma^{x, \varepsilon }_i$ are borrowed from the same lemma. For $x \in S_i,$ $$  \mathbb P \left({{\bf X}^{x, \varepsilon}_{\theta^{x, \varepsilon}} \in S^G}\right) =  \sum_{n = 1}^{\infty} \int_{S^G}  \int_{S_i} \cdots  \int_{S_i} \Gamma_i^{\varepsilon}(x, dx_1) \Gamma_i^{\varepsilon}(x_1, dx_2) \dots \Gamma_i^{\varepsilon}(x_{n-1}, dx_n ).  $$ By estimating this using a geometric series (as we did in lemma \ref{abstract_chain_lemma}) we find that $$  \dfrac{ \inf \limits_{x \in S_i} \Gamma_i^{\varepsilon} \left({x, S^G}\right)  }{ \sup \limits_{x \in S_i} \Gamma_i^{\varepsilon} \left({x, S^F}\right) } \le  \mathbb P \left({{\bf X}^{x, \varepsilon}_{\theta^{x, \varepsilon}} \in S^G}\right) \le \dfrac{ \sup \limits_{x \in S_i} \Gamma_i^{\varepsilon} \left({x, S^G}\right)  }{ \inf \limits_{x \in S_i} \Gamma_i^{\varepsilon} \left({x, S^F}\right) }. $$ Therefore, we obtain that $\mathbb P \left({{\bf X}^{x, \varepsilon}_{\theta^{x, \varepsilon}} \in S^G}\right) \sim \displaystyle \sum_{j \in E} \lambda(j) \displaystyle \sum_{k \in G} P_{jk}^{\varepsilon} \ \text{\Large /}  \displaystyle \sum_{j \in E} \lambda(j) \displaystyle \sum_{k \in F} P_{jk}^{\varepsilon}   $ uniformly in $x \in S_{i}.$ Since these asymptotics do not depend on $i,$ they hold uniformly in $x \in \bigcup \limits_{j \in E} S_j.$
\qed 
\\
\\

\section{Hierarchy of clusters} \label{hier_main}

 In this section, we return to considering our original Markov renewal process $\left({ {\bf X}_n^{x, \varepsilon}, {\bf T}_n^{x, \varepsilon} }\right)$ that satisfies assumptions (a) - (f) from Section \ref{Reduction to a process with nearly independent transition times and transition probabilities}. Our goal in this section will be to construct a hierarchy of subsets of $S$ that will prove to  be informative later. The transitions between different higher ranked sets in this hierarchy will be very rare. We will also construct a hierarchy of Markov renewal processes. For each such process, the state space $S$ will be partitioned into a collection of sets that will be unions of some of the sets $S_i$. We will show that assumptions  (a) - (f) from Section \ref{Reduction to a process with nearly independent transition times and transition probabilities} are satisfied for each such process. Let us make things precise.

Below, we will define a collection $\{S_i^r\}$ of subsets of $S$ for $0 \le r \le R$ and $1 \le i \le M_r.$ $ S_i^r$ will be called the $i$th cluster of rank $r.$ For a fixed $r, $ $\{S_1^r, \dots, S_{M_r}^r\}$ will form a partition of $S.$ For $r \ge 1, $ each $S_i^r$ will be equal to a union of certain clusters of rank $r-1.$ Before defining these clusters inductively, let us first fix some notation. For $0 \le r \le R, $ let $\pi^{r}$ be the projection onto the relevant cluster of rank $r$, i.e., if $x \in S_i^r$ (or $A \subseteq S_i^r$), then $\pi^{r} \left({x}\right)  = S_i^r$  ($ \pi^{r} \left({A}\right)= S_i^r$). Consider the following stopping times: $$ \theta_{0}^{x, \varepsilon}(r)  = 0, $$
$$ \theta_{n}^{x, \varepsilon}(r)  = \min \left\{{ m > \theta_{n-1}^{x, \varepsilon}(r) :  \pi^{r} \left({{\bf X}_m^{x, \varepsilon}}  \right)
\neq \pi^{r} \left({{\bf X}_{\theta_{n-1}^{x, \varepsilon}(r)}^{x, \varepsilon}}
\right)  }\right\}, \ n \ge 1. $$ Then, the kernel given by $P^{r, \varepsilon}(x, A) = \mathbb P \left({{\bf X}_{\theta_{1}^{x, \varepsilon}(r)}^{x, \varepsilon} \in A  }\right)$ describes the Markov chain $ {{\bf X}_{\theta_{n}(r)}^{x, \varepsilon}} $ that captures the transition of the process between clusters of rank $r.$  

We define $S_i^0 = S_i$ for $0 \le i \le M$ and $M_0 = M.$ We now describe the clusters $S_i^r$ inductively for $r \ge 1$.

Suppose that the clusters are defined up to rank $r$ for some $r \ge 0$.  Suppose that the number of distinct clusters of  rank $r$ is greater than one, i.e., $M_r \neq 1$. Then,  $\theta^{x, \varepsilon}_n(r),  \ {{\bf X}_{\theta_{n}(r)}^{x, \varepsilon}}$ and $P^{r, \varepsilon}( \cdot \ ,  \ \cdot)$ are all well-defined objects. Consider the following assumptions.  

(a) $P^{r, \varepsilon}(x, S_i^r) = 0 $ for every $x \in S_i^r$ and $ \varepsilon > 0.$ For every $ i \neq j,  \ P^{r, \varepsilon} (x, S^r_j)$ is either  identically $0$ for every $x \in S^r_i$ or it is positive for all $x \in S^r_i $ and $\varepsilon > 0,$ and in the latter case there are functions $P_{ij}^{r, \varepsilon}, \  \varepsilon > 0, $ such that  \begin{equation}
    \label{rank_r_hierarchy_prob_initial} \lim \limits_{\varepsilon \downarrow 0} \dfrac{P^{r, \varepsilon}(x, S^r_j)}{{P}_{ij}^{r, \varepsilon}} = 1  
\end{equation}uniformly in $x \in S^r_i.$

(b) Suppose that \begin{equation}
\label{asymptotic_regularity_rank_r_hierarchy}    \lim \limits_{\varepsilon \downarrow 0} \dfrac{P_{i_1j_1}^{r, \varepsilon}}{P_{i_1'j_1'}^{r, \varepsilon}} \cdots \dfrac{P_{i_nj_n}^{r, \varepsilon}}{P_{i_n'j_n'}^{r, \varepsilon}}   \in [0, \infty]
\end{equation}
 for every positive  integer $n$ and every $i_{\alpha} \neq j_{\alpha}, i_{\alpha}', \neq j_{\alpha}' $ for which the ratios appearing in the limit are defined.   
 
 We will show, inductively, that these assumptions are satisfied for all values of $r$. Notice that these are satisfied for $r = 0$ in line with our assumptions. These assumptions ensure that the limits $P_{ij}^{r, 0} = \lim \limits_{\varepsilon \to 0} P_{ij}^{r, \varepsilon}$ exist and  $\sum_{j \neq i} P_{ij}^{r, 0} = 1. $ Then $(P_{ij}^{r, 0})$ is the transition matrix for a Markov chain on $\{1, 2, \dots, M_r\}.$ Associated with this chain is a partition $\{ E_1, E_2, \dots E_{M_{r+1}} \}$ of $\{1, 2, \dots, M_r\}$ into ergodic classes and singleton sets consisting of transient states. This induces a partition of $S$ given by $ \{ S_1^{r+1}, \dots ,  S^{r+1}_{M_{r+1}} \},$ where $S_k^{r+1} = \bigcup \limits_{i \in E_k} S_i^r$, which explains the construction of the clusters of next rank. 

Let us show that the assumptions that were just outlined for rank $ r$ are also satisfied for rank $r+1.$ Fix $S_k^{r+1}, S_l^{r+1}$, where $k \neq l.$ Let $\lambda_k^{r}(\cdot)$ be the invariant measure on $E_k$ associated with the transition probabilities $P_{ij}^{r, 0}.$ In the future we will drop the subscript $k$ when convenient. Lemma \ref{abstract_chain_lemma_2} goes through for the chain ${\bf X}^{x, \varepsilon}_{\theta^{x, \varepsilon}_n(r)}$ with $E = E_k, \ G = E_l, \  \lambda = \lambda_k^{r} $ and $\theta^{x, \varepsilon} =  \min \left\{ { m  \ : \  \pi^{r+ 1 } \left({ {\bf X}^{x, \varepsilon}_{ \theta_m^{x, \varepsilon}(r) }  }\right) } \neq \pi^{r+1}(x) \right\} .$ Therefore, $$  P^{r + 1, \varepsilon} (x, S_l^{r+1})\sim \displaystyle \sum_{i : S_i^r \subseteq S_k^{r+1}}  \sum_{j : S_j^r \subseteq S_l^{r+1}} \lambda^{r}_k(i)P_{ij}^{r, \varepsilon}\text{\Large /} \displaystyle \sum_{i : S_i^r \subseteq S_k^{r+1}}  \sum_{j : S_j^r \not \subseteq S_k^{r+1}} \lambda^{r}_k(i)P_{ij}^{r, \varepsilon} \ \text{as} \ \varepsilon \to 0 $$ uniformly in $x \in S_k^{r+1}.$ We simply define $$P^{r+1, \varepsilon}_{kl} = \displaystyle \sum_{i : S_i^r \subseteq S_k^{r+1}}  \sum_{j : S_j^r \subseteq S_l^{r+1}} \lambda^{r}_k(i)P_{ij}^{r, \varepsilon}\text{\Large /} \displaystyle \sum_{i : S_i^r \subseteq S_k^{r+1}}  \sum_{j : S_j^r \not \subseteq S_k^{r+1}} \lambda^{r}_k(i)P_{ij}^{r, \varepsilon},$$ which is the quantity that appears in the right hand side above. From this definition it follows that assumption (a) holds for rank $r+1$ if it holds for rank $r.$ Assumption (b) (formula  \eqref{asymptotic_regularity_rank_r_hierarchy}) follows for rank $r+1$ due to Lemma  \ref{appendix_lemma_1}.

This completes the construction of the hierarchy of sets $\{S_i^r:  1 \le i \le M_r, 0 \le  r \le R\}$ and shows that properties (a) and (b) hold for all $0 \leq r <  R$.

\begin{remark}
    \label{remark_one_jump_prob_smaller}

Let $ 0 \le  r \le R - 2,1 \le  i_0, j_0 \le M_r, 1 \le k, l \le M_{r + 1} $ and $ k \neq l$. Suppose $S_{i_0}^r \subseteq S_k^{r + 1}, S_{j_0}^r \subseteq S_l^{r + 1} $ and $P_{i_0j_0}^{r, \varepsilon}$ is not identically zero.  From the definition of $P_{kl}^{r + 1, \varepsilon}$ stated above, it follows that $$ \begin{aligned} \dfrac{ P_{i_0j_0}^{r, \varepsilon} }{ P_{kl}^{r + 1, \varepsilon} }  &\le {P_{i_0j_0}^{r, \varepsilon} } \left({ \dfrac{ \lambda_k^r(i_0) P_{i_0j_0}^{r, \varepsilon} }{ \displaystyle \sum_{i : S_i^r \subseteq S_k^{r+1}}  \sum_{j : S_j^r \not \subseteq S_k^{r+1}} \lambda^{r}_k(i)P_{ij}^{r, \varepsilon}}}\right)^{-1} \\ &= \dfrac{1}{ \lambda_k^r(i_0) }  \displaystyle \sum_{i : S_i^r \subseteq S_k^{r+1}}  \sum_{j : S_j^r \not \subseteq S_k^{r+1}} \lambda^{r}_k(i)P_{ij}^{r, \varepsilon} \to 0 \end{aligned} $$ as $\varepsilon \downarrow 0$. Therefore, by induction, the same limit holds if $r$ is replaced with any $0 \le s \le r$ in the numerator in the left-hand side.   
    
\end{remark}

Recall that $\theta^{x, \varepsilon}_1{(r)}$ is the first time the process ${\bf X}_{n}^{x, \varepsilon}$ exits $\pi^{r}(S_i).$

\begin{definition} Let $1 \le i \le M, 0 \le r < R$ and $ x \in \pi^{r}(S_i)$. Let us define $$ N^{\varepsilon}(x, i, r) = \sum_{n = 0}^{\theta^{x, \varepsilon}_1{(r)}} \chi_{S_i}({\bf X}_{n}^{x, \varepsilon})$$
 to be the random variable that counts the number of visits to $S_i$ before ${\bf X}_{n}^{x, \varepsilon}$ exits $\pi^r(x)$, the cluster of rank $r$ to which $x$ belongs.

\end{definition} 
 
 \begin{lemma}
 \label{rank_r_expectation_lemma}
 Let $1 \le  r < R$  and $ 1 \le i \le M.$ For $1\le s \le r$, let $j_s$ be such that $\pi^{s}(S_i) = S^s_{j_s}$, Then
 $$\lim \limits_{\varepsilon \to 0} \ \mathbb E \left({N^{\varepsilon}(x, i, r)}\right) \times  \prod_{s = 0}^{r-1} \left( ({\lambda^{s}(j_s)})^{-1}
 \displaystyle  \sum_{\substack{j \ : \ S_j^s \subseteq \pi^{s+1}(S_i) \\ k \ : \ S_k^s \not \subseteq \pi^{s+1}(S_i)} }  {\lambda^{s}(j)}P_{jk}^{s, \varepsilon} \right) = 1$$  uniformly in $x \in \pi^{r}(S_i). $\end{lemma}
  \proof  This result is true for $r=1$ due to  Lemma \ref{abstract_chain_lemma}.   Let us assume the statement is true for ranks smaller than $r+1< R.$  
  
 Let $\tilde{\sigma}^{x,\varepsilon}_0 = 0$
  and 
 $\tilde{\sigma}^{x,\varepsilon}_{ n}$, $n \geq 1$,  be the time of the $n$th   return  (or visit, if $x \notin \pi^r(S_i)$) by ${\bf X}^{x,\varepsilon}_n$  to $\pi^r(S_i)$.  Define  
${\sigma}^{x,\varepsilon}_n =
\min(\tilde{\sigma}^{x,\varepsilon}_n, \theta^{x,\varepsilon}_1(r+1))$, $n \geq 0$. Thus $\sigma^{x, \varepsilon}_n$ are constant in $n$ for large $n$ corresponding to the times after the process exits $\pi^{r+1}(S_i)$.

Let $x \in \pi^{r+1}(S_i)$ and suppose that $\pi^{r+1}(S_i) \neq \pi^{r}(S_i)$. Note that $\{{\bf X}_{\sigma_{n }}^{x, \varepsilon} \in \pi^r(S_i)
\}$ is the event that the process visits $\pi^r(S_i)$ at least $n$ times before exiting $\pi^{r+1}(S_i).$ We observe that $\{{\bf X}_{\sigma_{n }}^{x, \varepsilon} \in \pi^r(S_i)
\} \subseteq \mathcal{F}^{x, \varepsilon}_{\sigma_{n}} $ and $$ N^{\varepsilon}(x, i, r+1)  = \displaystyle \sum_{n = 0}^{\infty} \chi_{\pi^r(S_i)} \left({ {\bf X}_{\sigma_n}^{x, \varepsilon} }\right)\displaystyle \sum_{m = \sigma_n^{x, \varepsilon}}^{ \sigma_{n + 1 }^{x, \varepsilon} -  1 } \chi_{S_i}\left({{\bf X}_{m}^{x, \varepsilon}}\right).$$ By the Strong Markov Property, for every $n  \ge 0,$   
$$  \begin{aligned} \inf \limits_{y \in \pi^{r}(S_i) } \mathbb E \left({N^{\varepsilon}(y, i, r)}\right) &\cdot \sum_{n=0}^{\infty}  \mathbb P \left({{\bf X}_{\sigma_{n }}^{x, \varepsilon} \in \pi^r(S_i)
}\right) \\ &\le \mathbb E \left({N^{\varepsilon}(x, i, r+1)}\right) \\ &\le \sup \limits_{y \in \pi^{r}(S_i) } \mathbb E \left({N^{\varepsilon}(y, i, r)}\right) \cdot \sum_{n=0}^{\infty}  \mathbb P \left({{\bf X}_{\sigma_{n }}^{x, \varepsilon} \in \pi^r(S_i)
}\right). \end{aligned} $$ By the induction hypothesis, $$  \inf \limits_{y \in \pi^{r}(S_i) } \mathbb E \left({N^{\varepsilon}(y, i, r)}\right)  \sim \sup \limits_{y \in \pi^{r}(S_i) } \mathbb E \left({N^{\varepsilon}(y, i, r)}\right)  \sim  \left({  \prod_{s = 0}^{r-1}\displaystyle  \sum_{\substack{j \ : \ S_j^s \subseteq \pi^{s+1}(S_i) \\ k \ : \ S_k^s \not \subseteq \pi^{s+1}(S_i)} }  \dfrac{\lambda^{s}(j)}{\lambda^{s}(j_s)}P_{jk}^{s, \varepsilon} } \right)^{-1} .$$ Furthermore, we realize $ \displaystyle  \sum_{n=0}^{\infty}  \mathbb P \left({{\bf X}_{\sigma_{n }}^{x, \varepsilon} \in \pi^r(S_i)
}\right)$ as the expectation of the number of visits of the process ${\bf X}_{\theta_n(r)}^{x, \varepsilon}$ to $\pi^r(S_i)$ before it exits $\pi^{r+1}(S_i).$ Lemma \ref{abstract_chain_lemma} applies directly to this situation establishing that $$ \lim \limits_{\varepsilon \to 0} \  \sum_{n=0}^{\infty}  \mathbb P \left({{\bf X}_{\sigma_{n }}^{x, \varepsilon} \in \pi^r(S_i)
}\right) \ \  \times \displaystyle  \sum_{\substack{j \ : \ S_j^r \subseteq \pi^{r+1}(S_i) \\ k \ : \ S_k^r \not \subseteq \pi^{r+1}(S_i)} }  \dfrac{\lambda^{r}(j)}{\lambda^{r}(j_r)}P_{jk}^{r, \varepsilon} = 1 $$ uniformly in $x \in \pi^{r+1}(S_i)$. This completes the proof.
\qed 
\\

Let us recall the hierarchy of  clusters $\left\{{S_i^r : 1 \le  i \le M_r, 0 \le r \le R}\right\}$ constructed above. Let us also recall the stopping times $\theta_{n}^{x, \varepsilon}(r)$ for $n \ge 1,0 \le r < R$  that record visits by ${\bf X}_n^{x, \varepsilon}$ to each new cluster of rank $r$. Recall also that for $0 \le r < R, P^{r, \varepsilon}( x, \cdot  )$ is a probability kernel on $S \times \mathcal B(S)$ defined by $$ P^{r, \varepsilon}( x, A ) = \mathbb P \left({ {\bf X}_{\theta_1(r)}^{x, \varepsilon} \in A}\right).$$ Finally let us recall the functions $P_{ij}^{r, \varepsilon}$ defined for $1 \le i, j \le M_r$ that are either identically zero or positive for every $\varepsilon$, and, in the latter case, satisfy $$ \lim \limits_{\varepsilon \downarrow 0} \dfrac{ P^{r, \varepsilon}(x, S_j^r) }{ P_{ij}^{r, \varepsilon} } = 1$$ uniformly in $x \in S_i^r$.

\begin{definition} \label{rank_r_definition} For $ \varepsilon > 0, 0 \le r  <  R$, and $x \in S, $ let us define the Markov chain $\left({ {\bf X}_n^{x, r, \varepsilon}, {\bf T}_n^{x, r, \varepsilon} }\right)$ by  $$  {\bf X}_n^{x, r, \varepsilon} = {\bf X}_{ \theta_n(r) }^{x, \varepsilon},   \ {\rm for \ } n \ge 0,  $$ $$ \ {\bf T}_0^{x, r, \varepsilon} = 0, \    {\bf T}_n^{x, r, \varepsilon} = \sum_{ m =  \theta^{x, \varepsilon}_{n - 1}(r) + 1}^{ \theta^{x, \varepsilon}_{n}(r)  } {\bf T}_{ m }^{x, \varepsilon}, \ {\rm for \ } n \ge 1. $$  The semi-Markov process driven by the Markov chain $\left({ {\bf X}_n^{x, r, \varepsilon}, {\bf T}_n^{x, r, \varepsilon} }\right)$ will be denoted by ${ X}_t^{x, r, \varepsilon}$. For completeness, we restate the definition of the probability kernel $P^{r, \varepsilon}(x, \cdot)$. For $ A \in \mathcal B (S)$, $$ P^{r, \varepsilon}(x, A) = \mathbb P \left({ {\bf X}_1^{x, r, \varepsilon } \in A }\right).  $$ For $1 \le  j \le M_r$ and $r < R$, we define $T^{r, \varepsilon}(x, S_j^r)$ to be a random variable with distribution given by $$ \mathbb P \left({T^{r, \varepsilon}(x, S_j^r) \in I }\right) = \mathbb P \left({ {\bf T}_1^{x, r, \varepsilon} \in I \vert {\bf X}_1^{x, r, \varepsilon} \in S_j^r }\right) $$ whenever $ P^{r, \varepsilon}(x, S_j^r) > 0.$ Noting that there is only one cluster of rank $R$, we define $ { X}_t^{x, R, \varepsilon} = x$ for all $x \in S$ and $ t \ge 0$. We also put $ {\bf T}_1^{x, R, \varepsilon} = \infty$ for every $x \in S$.  \end{definition}

$X_t^{x, r, \varepsilon}$ is a semi-Markov process on $S$ (as in  Section 1) associated with the Markov chain $\left({ {\bf X}_{n}^{x, r, \varepsilon} , {\bf T}_{n}^{x, r, \varepsilon} }\right).$ It makes jumps only between sets $\{ S_1^r, \dots, S_{M_r}^r \}.$ Note that $X_t^{x, 0, \varepsilon} $ is the same as our original process $ X_t^{x, \varepsilon}$.

 Recall from  \eqref{time_assumption_first_extended_space}  that for $1 \le i \le M$ we assume the existence of the functions $\tau_i^{\varepsilon}$ such that $\mathbb E \left({ T^{\varepsilon}(x, S_j)  }\right)/\tau_i^{\varepsilon} \to 
 1 $ as $\varepsilon \to 0$ uniformly in $x \in S_i,$ for every $1\le j \le M$, provided that $P_{ij}^{\varepsilon}$ is not identically zero.

\begin{theorem} \label{prelim_things_abt_rank_r_times}

Let $0 \le  r < R, 1 \le  k \le M_r$. There is a function $\tau_k^{r, \varepsilon}$ that satisfies $\tau_k^{r, \varepsilon} \sim   \displaystyle \sum_{i: S_i \subseteq S_k^r } \tau_i^{\varepsilon} \cdot  \mathbb E \left({ N^{\varepsilon}(x, i, r) }\right)$ uniformly in $x \in S_k^r$ such that $$ \lim \limits_{\varepsilon \downarrow 0}\dfrac{\mathbb E T^{r, \varepsilon}(x, S_l^r)}{\tau_k^{r, \varepsilon}} = 1 $$ uniformly in $x \in S_k^r$, for each $1 \le l \le M_r$, provided that $P_{kl}^{r, \varepsilon}$ is not identically zero.

\end{theorem}

\proof Note that the statement of the theorem is true for $r = 0$, in line with our assumptions in Section \ref{Reduction to a process with nearly independent transition times and transition probabilities}. So let us restrict ourselves to $1 \le r < R$. 

For a given $1 \le i \le M$ such that $S_i \subseteq S_k^r$, it was shown in Lemma \ref{rank_r_expectation_lemma} that the asymptotics of $\mathbb E \left({N^{\varepsilon}(x, i, r)}\right)$ does not depend on the choice of $x \in S_k^r$. So let us fix some $x_0 \in S_k^r$ and define the function $ \tau_k^{r, \varepsilon} $ by setting $$  \tau_k^{r, \varepsilon} = \displaystyle \sum_{i: S_i \subseteq S_k^r } \tau_i^{\varepsilon} \cdot  \mathbb E \left({ N^{\varepsilon}(x_0, i, r) }\right).   $$ Then $\tau_k^{r, \varepsilon} \sim  \displaystyle \sum_{i: S_i \subseteq S_k^r } \tau_i^{\varepsilon} \cdot  \mathbb E \left({ N^{\varepsilon}(x, i, r) }\right)$ uniformly in $x \in S_k^r$ as $\varepsilon \downarrow 0$.   

 For $x \in S$ and $ \varepsilon > 0 $, recall that ${\bf X}_1^{x, r, \varepsilon}$ is the location of the process when it first exits the  rank-$r$ cluster to which $x$ belongs. We claim that \begin{equation} \label{time_theorem_formula_1} \lim \limits_{\varepsilon \downarrow 0} \dfrac{\mathbb E \left({ {\bf T}_1^{x, \varepsilon} \chi_{S_l^r} \left({{\bf X}_{1 }^{x,r,  \varepsilon}}\right)    }\right)}{\tau_i^{\varepsilon} P_{kl}^{r, \varepsilon}}
     = 1 
\end{equation}  uniformly in $x \in S_i$ for every $1 \le i \le M, 1 \le k, l \le M_r$, whenever $S_i \subseteq  S_k^r$ and $P_{kl}^{r, \varepsilon} > 0$ for every $\varepsilon > 0$. Recall from \eqref{rank_r_hierarchy_prob_initial} that  $\mathbb P \left({ {\bf X}_{1}^{x, r,  \varepsilon} \in S_{l}^r }\right) \sim P_{kl}^{r, \varepsilon}$ uniformly in $x \in S_k^r$ as $\varepsilon \downarrow 0$. 

Conditioning on the location of ${\bf X}_1^{x, \varepsilon}$, we obtain \begin{equation} \label{time_theorem_formula_7}
\begin{aligned}  &\mathbb E \left({ {\bf T}_1^{x, \varepsilon} \chi_{S_l^r} \left({{\bf X}_{1 }^{x,r,  \varepsilon}}\right)  }\right)  \\  & = \mathbb E \left({ {\bf T}_1^{x, \varepsilon} \chi_{S_l^r} \left({{\bf X}_{1 }^{x,\varepsilon}}\right)   }\right) + \mathbb E \left({ {\bf T}_1^{x, \varepsilon} \chi_{S_l^r} \left({{\bf X}_{1 }^{x,r,  \varepsilon}}\right)  \chi_{S_k^r} \left({{\bf X}_{1 }^{x,\varepsilon}}\right)   }\right).  \end{aligned} \end{equation}  From \eqref{time_assumption_first_extended_space}, the first expectation in the sum has asymptotics given by $$ \mathbb E \left({ {\bf T}_1^{x, \varepsilon} \chi_{S_l^r} \left({{\bf X}_{1 }^{x,\varepsilon}}\right)  }\right) \sim \tau_i^{\varepsilon} \sum \limits_{j: S_{j} \subseteq S_l^r  } P_{ij}^{\varepsilon} $$ uniformly in $x \in S_i$ as $\varepsilon \downarrow 0$. 
 Observe that 
$P^\varepsilon_{ij}/P^{r,\varepsilon}_{kl} \rightarrow 0$ as $\varepsilon \downarrow 0$ when  $\pi^r(S_i) = S^r_k$ and $\pi^r(S_j) = S^r_l$  as a consequence of Remark \ref{remark_one_jump_prob_smaller}. Therefore, the first term on the right-hand side of (\ref{time_theorem_formula_7}) does not contribute to the limit in 
(\ref{time_theorem_formula_1}). 

We condition the second expectation in the sum on $\mathcal F_1, $ the $\sigma$-algebra generated by $\left({ {\bf X}_1^{x, \varepsilon}, {\bf T}_1^{x, \varepsilon}  }\right)$.  Then,
denoting $g(x, \varepsilon) =      \mathbb P \left({ { {\bf X}_{1}^{x, r,  \varepsilon} \in S_{l}^r }   }\right)  $,

$$ \begin{aligned} \mathbb E \left({ {\bf T}_1^{x, \varepsilon} \chi_{S_l^r} \left({{\bf X}_{1 }^{x,r,  \varepsilon}}\right)   \chi_{S_k^r} \left({{\bf X}_{1 }^{x,\varepsilon}}\right)   }\right) &=   \mathbb E \left({\mathbb E \left({ {\bf T}_1^{x, \varepsilon} \chi_{S_l^r} \left({{\bf X}_{1 }^{x,r,  \varepsilon}}\right)   \chi_{S_k^r} \left({{\bf X}_{1 }^{x,\varepsilon}}\right)   \vert   \mathcal F_1 }\right)  }\right) \\ &=  \mathbb E \left({  {\bf T}_1^{x, \varepsilon} \chi_{S_k^r} \left({{\bf X}_{1 }^{x,\varepsilon}}\right)   g( {\bf X}_{1 }^{x,\varepsilon}, \varepsilon) }\right)  \\ & \sim P_{kl}^{r, \varepsilon}    \mathbb E \left({  {\bf T}_1^{x, \varepsilon} \chi_{S_k^r} \left({{\bf X}_{1 }^{x,\varepsilon}}\right)    }
\right) \\ & \sim P_{kl}^{r, \varepsilon} \tau_i^{\varepsilon} \sum \limits_{j: S_{j} \subseteq S_k^r  } P_{ij}^{\varepsilon}  \end{aligned} $$  uniformly in $x \in S_i$ as $\varepsilon \downarrow 0$.
 Formula \eqref{time_theorem_formula_1} now follows by noticing that $ \sum \limits_{j: S_{j} \subseteq S_k^r  } P_{ij}^{\varepsilon} \rightarrow 1 $ as $\varepsilon \downarrow 0$. 

  For $x \in S_k^r$ and $i$ such that $S_i \subseteq S_k^r$,  let us define $\tilde{\sigma}^{x,\varepsilon}_0(i) = 0$
  and 
$\tilde{\sigma}^{x,\varepsilon}_n (i)$, $n \geq 1$, to be the time of the $n$th visit (or return, if $x \in S_i$)   to $S_i$. Define  
${\sigma}^{x,\varepsilon}_n (i) =
\min(\tilde{\sigma}^{x,\varepsilon}_n(i), \theta^{x,\varepsilon}_1(r))$, $n \geq 0$. Then $$ \mathbb E \left({T^{r, \varepsilon}(x, S_l^r)}\right) =   \dfrac{ \mathbb E \left({  \displaystyle  \sum \limits_{i : S_i \subseteq S_k^r} \sum_{n = 0}^{\infty} {\bf T}_{\sigma_n(i) + 1}^{x, \varepsilon} \cdot \chi_{S_{i}} \left({{\bf X}_{\sigma_n(i) }^{x, \varepsilon}}\right) 
 \chi_{S_l^r} \left({{\bf X}_{1 }^{x,r,  \varepsilon}}\right)}\right) }{ \mathbb P \left({ {\bf X}^{x, r  , \varepsilon}_{1} \in S_l^{r} }\right) 
 }.   $$ We wish to examine the expectation of the interior sum appearing in the numerator for a fixed $i$ such that $S_i \subseteq S_k^r$. Let $h(x, \varepsilon) = \mathbb E \left({ {\bf T}_{ 1}^{ x, \varepsilon} \chi_{S_l^r} \left({{\bf X}_{1 }^{x,r,  \varepsilon}}\right)  }\right)$. By conditioning on $ \mathcal F_{\sigma_n^{x, \varepsilon}(i)}$ and using \eqref{time_theorem_formula_1}, $$ \begin{aligned} \mathbb E \left({ \displaystyle \sum_{n = 0}^{\infty} {\bf T}_{\sigma_n(i) + 1}^{x, \varepsilon} \cdot \chi_{S_i} \left({{\bf X}_{\sigma_n(i) }^{x, \varepsilon}}\right) \chi_{S_l^r} \left({{\bf X}_{1 }^{x,r,  \varepsilon}}\right)   }\right) & = \sum_{n = 0}^{\infty} \mathbb E \left({   \chi_{S_i} \left({{\bf X}_{\sigma_n(i) }^{x, \varepsilon}}\right) h \left({ {\bf X}_{\sigma_n(i) }^{x, \varepsilon}, \varepsilon  }\right)}\right) \\ & \sim \tau_i^{\varepsilon} P_{kl}^{r, \varepsilon} \sum_{n = 0}^{\infty} \mathbb P \left({ {\bf X}^{x, \varepsilon}_{\sigma_n(i)} \in S_i }\right) \\ & \sim  \tau_i^{\varepsilon} P_{kl}^{r, \varepsilon} \mathbb E \left({N^{\varepsilon}(x, i, r)}\right) 
\end{aligned} $$  uniformly in $x \in S_k^r$. This proves the Theorem.  \qed \\



    \begin{remark}
        \label{prelim_things_abt_rank_r_times_remark} For  rank $R$, where there is just one cluster $S_1^R $, we can set $\tau_1^{R, \varepsilon} = \infty$.  
    \end{remark}

Consider again the functions $\tau_k^{r, \varepsilon}, 0 \le  r < R, 1 \le  k \le M_r$, introduced in Theorem \ref{prelim_things_abt_rank_r_times}.  

\begin{theorem} \label{prelim_things_abt_rank_r_times_var}

Let $0 \le  r < R, 1 \le  k \le M_r$. Then there are $C, \varepsilon_0 > 0$ such that $$ \mathbb E  \left({T^{r, \varepsilon}(x, S_l^r)}\right)^2 \le C \cdot  \left({\tau_k^{r, \varepsilon}}\right)^2  $$ for each $x \in S_k^r, 1 \le  l \le M_r$, and $\varepsilon \le \varepsilon_0$, provided that $P_{kl}^{r, \varepsilon}$ is not identically zero.  \end{theorem}

\proof Let $1 \le l \le M_r$  be such that $P_{kl}^{r, \varepsilon}$ is a positive function. For $x \in S_k^r$, let us recall the stopping times $\sigma_n^{x, \varepsilon}(i), n \ge 0$ (the time of the $n$th visit to $S_i \subseteq S_k^r$ prior to exiting $S_k^r$) that were  introduced in the proof of Theorem \ref{prelim_things_abt_rank_r_times}. Then $$ \mathbb E \left({T^{r, \varepsilon}(x, S_l^r)}\right)^2 =   \dfrac{ \mathbb E \left({  \displaystyle  \sum \limits_{i : S_i \subseteq S_k^r} \sum_{n = 0}^{\infty} {\bf T}_{\sigma_n(i) + 1}^{x, \varepsilon} \cdot \chi_{S_{i}} \left({{\bf X}_{\sigma_n(i) }^{x, \varepsilon}}\right) 
 \chi_{S_l^r} \left({{\bf X}_{1 }^{x,r,  \varepsilon}}\right)}\right)^2 }{ \mathbb P \left({ {\bf X}^{x, r  , \varepsilon}_{1} \in S_l^{r} }\right) 
 }.   $$ Since the sum in $i$ is finite, it is enough to estimate the following for a fixed $i$.   
\begin{equation}
\begin{aligned}  \label{time_theorem_formula_4}  & \displaystyle    \mathbb E \left({   \sum_{n = 0}^{\infty} {\bf T}_{\sigma_n(i)  + 1}^{x, \varepsilon} \cdot \chi_{S_{i}} \left({{\bf X}_{\sigma_n(i) }^{x, \varepsilon}}\right)   
 \chi_{S_l^r} \left({{\bf X}_{1 }^{x,r,  \varepsilon}}\right)}\right)^2 
 \\   & =    \sum_{n = 0}^{\infty} \mathbb E \left({ \left({ {\bf T}_{\sigma_n(i)  + 1}^{x, \varepsilon} }\right)^2 \chi_{S_{i}} \left({{\bf X}_{\sigma_n(i) }^{x, \varepsilon}}\right)  \chi_{ S_l^{r} } \left({{\bf X}^{x, r,  \varepsilon}_{1}}\right) }\right) \ +  \\ & \;\;\;\;\;\;\;  2 \sum_{n = 0}^{\infty} \sum_{m = n + 1}^{\infty} \mathbb E \left({  {\bf T}_{\sigma_n(i) + 1}^{x, \varepsilon} {\bf T}_{\sigma_m(i) + 1}^{x, \varepsilon}\chi_{S_{i}} \left({{\bf X}_{\sigma_n(i) }^{x, \varepsilon}}\right) \chi_{S_{i}} \left({{\bf X}_{\sigma_m(i) }^{x, \varepsilon}}\right) \chi_{ S_l^{r} } \left({{\bf X}^{x, r, \varepsilon}_{1}}\right)}\right)
 . \end{aligned}  \end{equation}

 We first claim that that there is  $ C_1 > 0 $ such that \begin{equation} \label{time_theorem_formula_2}\mathbb E \left({ \left({ {\bf T}_1^{x, \varepsilon}}\right)^2 \chi_{S_l^r} \left({{\bf X}_{1 }^{x,r,  \varepsilon}}\right) }\right) \le C_1 \cdot (\tau_i^{\varepsilon})^2 P_{kl}^{r, \varepsilon}    \end{equation} for every $x \in S_i$ and small $\varepsilon$. The proof of this claim is similar to that of \eqref{time_theorem_formula_1} in the proof of Theorem \ref{prelim_things_abt_rank_r_times}. Indeed, for $x \in S_i$, $$ \begin{aligned}  &\mathbb E \left({ \left({ {\bf T}_1^{x, \varepsilon}}\right)^2 \chi_{S_l^r} \left({{\bf X}_{1 }^{x,r,  \varepsilon}}\right)  }\right)  \\  & = \mathbb E \left({ \left({ {\bf T}_1^{x, \varepsilon}}\right)^2 \chi_{S_l^r} \left({{\bf X}_{1 }^{x,\varepsilon}}\right)   }\right) + \mathbb E \left({ \left({ {\bf T}_1^{x, \varepsilon}}\right)^2 \chi_{S_l^r} \left({{\bf X}_{1 }^{x,r,  \varepsilon}}\right)  \chi_{S_k^r} \left({{\bf X}_{1 }^{x,\varepsilon}}\right)   }\right).  \end{aligned}  $$  Let us look at the first term on the right hand side. As a consequence of \eqref{varianceest_extended_space}, there is a constant $C_2 > 0$ such that $$ \mathbb E \left({ \left({ {\bf T}_1^{x, \varepsilon}}\right)^2 \chi_{S_l^r} \left({{\bf X}_{1 }^{x,\varepsilon}}\right)   }\right) \le C_2 \cdot \left({\tau_i^{\varepsilon}}\right)^2 \sum \limits_{j: S_{j} \subseteq S_l^r  } P_{ij}^{\varepsilon}   $$ for every $x \in S_i $ and small enough $\varepsilon$. For the second term, as before, we condition on $\mathcal F_1$ and utilize the function $g_1(x, \varepsilon) = \mathbb P \left({ {\bf X}_1^{x, r, \varepsilon} \in S_l^r}\right)$ with the Markov property. We find that there is a constant $C_3 > 0$ such that  $$ \begin{aligned} \mathbb E \left({ \left({{\bf T}_1^{x, \varepsilon}}\right)^2 \chi_{S_l^r} \left({{\bf X}_{1 }^{x,r,  \varepsilon}}\right)   \chi_{S_k^r} \left({{\bf X}_{1 }^{x,\varepsilon}}\right)   }\right) &=   \mathbb E \left({\mathbb E \left({  \left({{\bf T}_1^{x, \varepsilon}}\right)^2\chi_{S_l^r} \left({{\bf X}_{1 }^{x,r,  \varepsilon}}\right)   \chi_{S_k^r} \left({{\bf X}_{1 }^{x,\varepsilon}}\right)   \vert   \mathcal F_1 }\right)  }\right) \\ &=  \mathbb E \left({  \left({{\bf T}_1^{x, \varepsilon}}\right)^2 \chi_{S_k^r} \left({{\bf X}_{1 }^{x,\varepsilon}}\right)   g_1( {\bf X}_{1 }^{x,\varepsilon}, \varepsilon) }\right)  \\ & \sim P_{kl}^{r, \varepsilon}    \mathbb E \left({   \left({{\bf T}_1^{x, \varepsilon}}\right)^2 \chi_{S_k^r} \left({{\bf X}_{1 }^{x,\varepsilon}}\right)    }
\right) \\ &\le C_3 \cdot  \left({\tau_i^{\varepsilon}}\right)^2 \cdot  P_{kl}^{r, \varepsilon}  \sum \limits_{j: S_{j} \subseteq S_k^r  } P_{ij}^{\varepsilon}  \end{aligned} $$ for every $x \in S_i$ and small enough $\varepsilon$. Since $ \left({P_{kl}^{r, \varepsilon} }\right)^{-1} \cdot \sum \limits_{j: S_{j} \subseteq S_l^r  } P_{ij}^{\varepsilon}  \to 0  $ (due to Remark \ref{remark_one_jump_prob_smaller}) and $ \displaystyle \sum \limits_{j: S_{j} \subseteq S_k^r  } P_{ij}^{\varepsilon}  \to 1 $ as $\varepsilon \downarrow  0$, we obtain  \eqref{time_theorem_formula_2}.

We will now estimate the two sums appearing on the right hand side of \eqref{time_theorem_formula_4}. For the first sum, we condition each term on $\mathcal F_{\sigma_n(i)}$ and use the Markov
property. Let $g_2(x, \varepsilon) = \mathbb E \left({ \left({ {\bf T}_1^{x, \varepsilon}}\right)^2 \chi_{S_l^r} \left({{\bf X}_{1 }^{x,r,  \varepsilon}}\right) }\right)$. Note that $g_2(x, \varepsilon) \le C_1 \cdot (\tau_{i}^{\varepsilon})^2 P_{kl}^{r, \varepsilon}$ for every $x \in S_k^r$ and small $\varepsilon$ due to \eqref{time_theorem_formula_2}. Let $\delta(\varepsilon)$ be a function that satisfies $ \delta(\varepsilon) \sim \mathbb E \left({N^{\varepsilon} (x, i, r) }\right)$ uniformly in $x \in S_k^r$. Then $$ \begin{aligned} \sum_{n = 0}^{\infty} \mathbb E \left({ \left({ {\bf T}_{\sigma_n(i) + 1}^{x, \varepsilon} }\right)^2 \chi_{S_{i}} \left({{\bf X}_{\sigma_n^{x, \varepsilon}(i) }}\right) \chi_{ S_l^{r} } \left({{\bf X}^{x, r,  \varepsilon}_{1}}\right) }\right) &= \sum_{n = 0}^{\infty} \mathbb E \left({\chi_{S_{i}} \left({{\bf X}_{\sigma_n^{x, \varepsilon}(i) }}\right) g_2\left({{\bf X}^{x,  \varepsilon}_{n}, \varepsilon}\right)    }\right) \\ &\le C_1 \cdot (\tau_i^{\varepsilon})^2 P_{kl}^{r, \varepsilon} \sum_{n = 0}^{\infty} \mathbb P  \left({\sigma_n^{x, \varepsilon} (i) \in S_i }\right) \\   &\sim C_1 \cdot (\tau_i^{\varepsilon})^2 \delta(\varepsilon) P_{kl}^{r, \varepsilon}  \end{aligned}  $$ for every $x \in S_k^r$ and small $\varepsilon$. 
 
 Now we look at the second sum on the right hand side of \eqref{time_theorem_formula_2}. For $x \in S_i$, let  $$g_3 (x, \varepsilon) = \mathbb E \left({ {\bf T}_1^{x, \varepsilon} \chi_{S_l^r} \left({ {\bf X}_1^{x,r,  \varepsilon}}\right) }\right),$$   $$g_4(x, \varepsilon )=  \mathbb E \left({N^{\varepsilon}(x, i, r)  }\right) = \sum_{n = 0}^{\infty}  \mathbb P \left({\sigma_n^{x, \varepsilon} (i) \in S_i}\right), $$ $$g_5(x, \varepsilon) = \mathbb E \left({{\bf T}_1^{x, \varepsilon}}\right).$$ Recall that there is a constant $C_4$ such that $g_3(x, \varepsilon) \le  C_4 \cdot \tau_i^{\varepsilon}P_{kl}^{r, \varepsilon} $  for every $x \in S_i$ and small $\varepsilon$ (see  \eqref{time_theorem_formula_1}  from the proof of Theorem \ref{prelim_things_abt_rank_r_times}). We also have $g_4(x, \varepsilon) \sim \delta(\varepsilon)$ uniformly in $x \in S_k^r$ as $\varepsilon \downarrow 0$ by definition. There is also a constant $C_5$ such that $g_5(x, \varepsilon) \le C_5 \tau_i^{\varepsilon}$ for every $x \in S_i$ and small $\varepsilon$. Then

 $$\begin{aligned} &  \sum_{n = 0}^{\infty} \sum_{m = n + 1}^{\infty} \mathbb E \left({  {\bf T}_{\sigma_n(i) + 1}^{x, \varepsilon} {\bf T}_{\sigma_m(i) + 1}^{x, \varepsilon}\chi_{S_{i}} \left({{\bf X}_{\sigma_n(i) }^{x, \varepsilon}}\right) \chi_{S_{i}} \left({{\bf X}_{\sigma_m(i) }^{x, \varepsilon}}\right) \chi_{ S_l^{r} } \left({{\bf X}^{x, r, \varepsilon}_{1}}\right)}\right) \\  &= \sum_{n = 0}^{\infty} \sum_{m = n + 1}^{\infty} \mathbb E\left({\mathbb E \left({  {\bf T}_{\sigma_n(i) + 1}^{x, \varepsilon} {\bf T}_{\sigma_m(i) + 1}^{x, \varepsilon} \chi_{S_{i}} \left({{\bf X}_{\sigma_n(i) }^{x, \varepsilon}}\right) \chi_{S_{i}} \left({{\bf X}_{\sigma_m(i) }^{x, \varepsilon}}\right)  \chi_{ S_l^{r} } \left({{\bf X}^{x, r, \varepsilon}_{1}}\right)  \Bigm\vert \mathcal F_{\sigma_m(i)} }\right) }\right) \\ & = \sum_{n = 0}^{\infty} \sum_{m = n + 1}^{\infty}  \mathbb E \left({ {\bf T}_{\sigma_n(i) + 1}^{x, \varepsilon} \chi_{S_{i}} \left({{\bf X}_{\sigma_n(i) }^{x, \varepsilon}}\right) \chi_{S_{i}} \left({{\bf X}_{\sigma_m(i) }^{x, \varepsilon}}\right)  g_3\left({ {\bf X}^{x,  \varepsilon}_{\sigma_m(i)}, \varepsilon }\right) }\right) \\&\le C_4 \cdot \tau_i^{\varepsilon} P_{kl}^{r, \varepsilon} \sum_{n = 0}^{\infty} \sum_{m = n + 1}^{\infty}  \mathbb E \left({ {\bf T}_{\sigma_n(i) + 1}^{x, \varepsilon}   \chi_{S_{i}} \left({{\bf X}_{\sigma_n(i) }^{x, \varepsilon}}\right) \chi_{S_{i}} \left({{\bf X}_{\sigma_m(i) }^{x, \varepsilon}}\right) }\right)\end{aligned} $$ $$ \begin{aligned}  &= C_4 \cdot \tau_i^{\varepsilon} P_{kl}^{r, \varepsilon} \sum_{n = 0}^{\infty} \sum_{m = n + 1}^{\infty} \mathbb E \left({ \mathbb E \left({ {\bf T}_{\sigma_n(i) + 1}^{x, \varepsilon} \chi_{S_{i}} \left({{\bf X}_{\sigma_n(i) }^{x, \varepsilon}}\right) \chi_{S_{i}} \left({{\bf X}_{\sigma_m(i) }^{x, \varepsilon}}\right) \Bigm\vert \mathcal F_{\sigma_n(i) + 1}  }\right) }\right)  \\ & = C_4 \cdot \tau_i^{\varepsilon} P_{kl}^{r, \varepsilon} \sum_{n = 0}^{\infty}  \mathbb E \left({ {\bf T}_{\sigma_n(i) + 1}^{x, \varepsilon} \chi_{S_{i}} \left({{\bf X}_{\sigma_n(i) }^{x, \varepsilon}}\right) g_4\left({ {\bf X}^{x,  \varepsilon}_{\sigma_n(i)+1}, \varepsilon }\right) }\right) \\ &\sim   C_4 \cdot \tau_i^{\varepsilon} \delta(\varepsilon)  P_{kl}^{r, \varepsilon} \sum_{n = 0}^{\infty}  \mathbb E \left({ {\bf T}_{\sigma_n(i) + 1}^{x, \varepsilon} \chi_{S_{i}} \left({{\bf X}_{\sigma_n(i) }^{x, \varepsilon}}\right)  }\right)  \\ &=C_4 \cdot \tau_i^{\varepsilon} \delta(\varepsilon)  P_{kl}^{r, \varepsilon} \sum_{n = 0}^{\infty}  \mathbb E \left({\mathbb E \left({ {\bf T}_{\sigma_n(i)  + 1}^{x, \varepsilon} \chi_{S_{i}} \left({{\bf X}_{\sigma_n(i) }^{x, \varepsilon}}\right)   \Bigm\vert \mathcal F_{\sigma_n(i) } }\right) }\right)\\ &= C_4 \cdot \tau_i^{\varepsilon} \delta(\varepsilon)  P_{kl}^{r, \varepsilon} \sum_{n = 0}^{\infty}  \mathbb E \left({ \chi_{S_{i}} \left({{\bf X}_{\sigma_n(i) }^{x, \varepsilon}}\right)  g_5\left( { {\bf X}_{n}^{x, \varepsilon}, \varepsilon }\right) }\right) \\ &\le C_4 C_5 \cdot (\tau_i^{\varepsilon})^2 \delta(\varepsilon)  P_{kl}^{r, \varepsilon} \sum_{n = 0}^{\infty}  \mathbb P \left({{\bf X}_{\sigma_n(i) }^{x, \varepsilon} \in S_i }\right)  \\ &\sim C_4 C_5 \cdot (\tau_i^{\varepsilon})^2 \delta(\varepsilon)^2 P_{kl}^{r, \varepsilon}, \end{aligned} $$ where the equalities and inequalities above hold for every  $x \in S_k^r$ and small $\varepsilon$, and the asymptotic equivalences hold uniformly $x \in S_k^r$ as $\varepsilon \downarrow 0$.

 We have obtained the following estimate on the left hand side of \eqref{time_theorem_formula_4}: there are $C, \varepsilon_0> 0$ such that $$ \displaystyle    \mathbb E \left({   \sum_{n = 0}^{\infty} {\bf T}_{\sigma_n(i)  + 1}^{x, \varepsilon} \cdot \chi_{S_{i}} \left({{\bf X}_{\sigma_n(i) }^{x, \varepsilon}}\right)   
 \chi_{S_l^r} \left({{\bf X}_{1 }^{x,r,  \varepsilon}}\right)}\right)^2 \le C  \cdot (\tau_i^{\varepsilon})^2 \delta(\varepsilon)^2 P_{kl}^{r, \varepsilon} $$ for every $x \in S_k^r$ and $\varepsilon \le \varepsilon_0$, where we recall that $\delta(\varepsilon) \sim \mathbb E \left({N^{\varepsilon}(x, i, r)}\right)$ uniformly in $x \in S_k^r$ as $\varepsilon \downarrow 0$. This is enough to conclude the statement of the Theorem.  \qed \\

Next, we would like to show that
$ {T^{r, \varepsilon}(x, S_l^r)}/{\tau_k^{r, \varepsilon}}$ converges to an exponential random variable for $x \in S^r_k$.  First, however, we will prove the following more abstract lemma. For the purpose of this lemma, $\mathcal{S}$ is an abstract parameter space. The main point of the lemma is that, under appropriate assumptions, 
a sum of a random number of random variables is almost exponentially distributed, provided that the number of terms is
determined by stopping the summation after each term with a small positive probability.
\begin{lemma} \label{abstrac_exp_lemma} For  $x \in \mathcal S$ and $ \varepsilon > 0$,  let $  \xi_n^{x, \varepsilon}, \zeta_{n}^{x, \varepsilon},  n \ge 0  $, be two sequences of random variables adapted to a filtration $\mathcal {F}^{x, \varepsilon}_n$. Suppose that $\xi_0^{x, \varepsilon} = 0$ and that there is a function $e(\varepsilon)$ such that  $$ \lim \limits_{\varepsilon \downarrow 0 } \dfrac{\mathbb E \left({ \xi_{n+1}^{x, \varepsilon} \vert  \mathcal {F}^{x, \varepsilon}_n }\right) }{ e(\varepsilon) } = 1  $$ uniformly in $x \in \mathcal{S}, n \ge 0$, and a version of the conditional expectation (which is defined up to a set of probability zero) can be taken such that the limit is uniform  with respect to the elements of the probability space. Suppose also that $$  \mathbb E \left({( \xi_{n+1}^{x, \varepsilon})^2 \vert  \mathcal {F}^{x, \varepsilon}_n }\right) \le C (e(\varepsilon))^2 $$ for some constant $C$, small enough $\varepsilon$, every $x \in \mathcal{S}$, and every $n \ge 0$. Suppose that $\zeta_{n}^{x, \varepsilon}$ takes values in $\{0, 1\}$, $\zeta^{x,\varepsilon}_0 = 1$, and that $ \zeta_{n}^{x, \varepsilon} = 0 $ implies $ \zeta_{m}^{x, \varepsilon} = 0 $ for $m \ge n$. Suppose that there is a positive function $r(\varepsilon)$ satisfying $  \lim \limits_{\varepsilon \downarrow 0 } r(\varepsilon)  = 0$ such that 
    $$ \lim \limits_{\varepsilon \downarrow 0}  \dfrac{\mathbb P \left({ \zeta_{n+1}^{x, \varepsilon} = 0  \vert \zeta_{n}^{x, \varepsilon} = 1 }\right) }{ r(\varepsilon) }  = 1$$ uniformly in $x \in \mathcal{S}, n \ge 0$. Let $ K^{x, \varepsilon} = \min \{ n  : \zeta_{n}^{x, \varepsilon} = 0 \}$.  Then for each $t \ge 0$,   $$ \lim \limits_{\varepsilon \downarrow 0} \mathbb P \left({ r(\varepsilon) K^{x, \varepsilon} \ge t  }\right) = e^{-t},  $$ 
        $$ \lim \limits_{\varepsilon \downarrow 0} \mathbb P \left({ {\displaystyle (e(\varepsilon))^{-1}r(\varepsilon) \sum_{n = 0}^{\infty }  \xi_{n+1}^{x, \varepsilon}\zeta_{n}^{x, \varepsilon}}  \ge t  }\right) = e^{-t}   $$
        uniformly in $x \in \mathcal S$. 

\end{lemma}

\proof For $\varepsilon > 0$, let $$ p(\varepsilon) = \sup \limits_{x \in \mathcal S, n \ge 0} \mathbb P \left({ \zeta_{n+1}^{x, \varepsilon} = 0  \vert \zeta_{n}^{x, \varepsilon} = 1 }\right),  q(\varepsilon) = \inf \limits_{x \in \mathcal S, n \ge 0} \mathbb P \left({ \zeta_{n+1}^{x, \varepsilon} = 0  \vert \zeta_{n}^{x, \varepsilon} = 1 }\right).  $$ Then  for each $n \ge 0, x \in \mathcal S$ and $\varepsilon > 0$,   $$ (1 - p(\varepsilon))^{n} \le \mathbb P \left({ K^{x, \varepsilon} >  n  }\right) \le  (1 - q(\varepsilon))^{n}. $$ Furthermore, our assumptions guarantee that  $p(\varepsilon) \sim r(\varepsilon), q(\varepsilon) \sim r(\varepsilon) $ as $\varepsilon \downarrow 0$. Then $$ \limsup \limits_{\varepsilon \downarrow 0} \mathbb P \left({ r(\varepsilon) K^{x, \varepsilon} > t  }\right) \le \limsup \limits_{\varepsilon \downarrow 0} \left({(1 - q(\varepsilon))^{\frac{1}{r(\varepsilon)}} } \right)^{t} = e^{-t} $$ for every $x \in \mathcal S$. The lower limit is obtained using a similar argument by replacing $q(\varepsilon)$ with $p(\varepsilon)$.

Let us now prove the second claim. Let $t > 0$ and $ \delta$ be a small positive quantity. Then  \begin{equation} \label{abstrac_exp_lemma_formula_1}
     \begin{aligned}
   & \mathbb P \left({  {\displaystyle e(\varepsilon)^{-1}r(\varepsilon) \sum_{n = 0}^{\infty }  \xi_{n+1}^{x, \varepsilon}\zeta_{n}^{x, \varepsilon}} > t + \delta \  , \  r(\varepsilon) K^{x, \varepsilon} <  t }\right) \\  & \le  \mathbb P \left({  {\displaystyle  e(\varepsilon)^{-1} r(\varepsilon) \sum_{n = 0}^{ 
 \left\lfloor {\frac {t}{r(\varepsilon)} } \right \rfloor - 1 }  \xi_{n+1}^{x, \varepsilon}}     >   {(t + \delta)}{} }\right) \\ & \le  \mathbb P \left({ \displaystyle  \sum_{n = 0}^{ 
 \left\lfloor {\frac {t}{r(\varepsilon)} } \right \rfloor - 1} \left({\dfrac{  \xi_{n+1}^{x, \varepsilon}} { e(\varepsilon)  } - 1}\right)     >   \dfrac{\delta}{r(\varepsilon)} }\right) \\  
     & \le \delta^{-2} r(\varepsilon)^2   \displaystyle  \sum_{n = 0}^{ 
 \left\lfloor {\frac {t}{r(\varepsilon)} } \right \rfloor - 1 } \mathbb E \left({\dfrac{  \xi_{n+1}^{x, \varepsilon}} { e(\varepsilon)  } - 1}\right)^2 \\&  \;\;\;\;\;\;\;\;\;\;\;\;    +  \ \ {2 \delta^{-2} r(\varepsilon)^2}   \displaystyle  \sum_{n = 0}^{ 
 \left\lfloor {\frac {t}{r(\varepsilon)} } \right \rfloor - 1 }  \sum_{m = n +1}^{ 
 \left\lfloor {\frac {t}{r(\varepsilon)} } \right \rfloor - 1 } \mathbb E \left({\dfrac{  \xi_{n+1}^{x, \varepsilon}} { e(\varepsilon)  } - 1 }\right)  \left({\dfrac{  \xi_{m+1}^{x, \varepsilon}} { e(\varepsilon)  } - 1}\right).   
 \end{aligned} \end{equation}  After conditioning each summand on $\mathcal F_{n}^{x, \varepsilon}$,  the first sum on the right hand side is seen to be bounded by $ Ct\delta^{-2} r(\varepsilon)  $ for some constant $C$, for every $x \in \mathcal S$ and small $\varepsilon$. We next note that there is a positive function $\alpha(\varepsilon)$ (that tends to zero as $\varepsilon \downarrow 0$) such that $$ - \alpha(\varepsilon) \le \mathbb E \left({ \dfrac{\xi_{m + 1}^{x, \varepsilon}}{e(\varepsilon) } - 1 \Bigm \vert \mathcal F_m^{x, \varepsilon}  }\right) \le  \alpha(\varepsilon)   $$ for each $m$. Therefore, the terms in the second sum above are bounded by $$ \alpha(\varepsilon)  \left| {\mathbb E \left({\dfrac{  \xi_{n+1}^{x, \varepsilon}} { e(\varepsilon)  } - 1 }\right) } \right| \le \alpha(\varepsilon) \left({  \mathbb E \left({\mathbb E \left({\dfrac{  \xi_{n+1}^{x, \varepsilon}} { e(\varepsilon)  } - 1 }\right)^2 \Bigm \vert \mathcal F_{n}^{x, \varepsilon} }\right)}\right)^{\frac 1 2} \le C \cdot \alpha(\varepsilon) $$ for some constant $C$, for every $x \in \mathcal S$ and small $\varepsilon$. We conclude that for every $t > 0$ and small $\delta > 0$,  the expression on the left hand side of \eqref{abstrac_exp_lemma_formula_1} tends to zero as $\varepsilon \downarrow 0$ uniformly in $x \in \mathcal S$. 

 Similarly, $$ \begin{aligned} &  \mathbb P \left({  {\displaystyle e(\varepsilon)^{-1} r(\varepsilon) \sum_{n = 0}^{\infty }  \xi_{n+1}^{x, \varepsilon}\zeta_{n}^{x, \varepsilon}} < t 
 \  , \  r(\varepsilon) K^{x, \varepsilon} >   t  + \delta }\right)  \\ & \le  \mathbb P \left({  {\displaystyle  e(\varepsilon)^{-1} \sum_{n = 0}^{ 
 \left\lfloor {\frac {t + \delta }{r(\varepsilon)} } \right \rfloor }  \xi_{n+1}^{x, \varepsilon}}     <    \dfrac{t }{r(\varepsilon)} }\right)  \\ & \le  \mathbb P \left({ \displaystyle  \sum_{n = 0}^{ 
 \left\lfloor {\frac {t + \delta}{r(\varepsilon)} } \right \rfloor } \left({ 1 - \dfrac{  \xi_{n+1}^{x, \varepsilon}} { e(\varepsilon)  } }\right)     >   \dfrac{\delta}{r(\varepsilon)} }\right).  
     \end{aligned} $$ Using arguments similar to those  above, it is easy to show that   the probability on the left hand side of this inequality tends to zero as $\varepsilon \downarrow 0$ uniformly in $x \in \mathcal S$. Combined with the first claim of the lemma, this is enough to conclude that the second claim holds.  \qed \\
 
        \begin{theorem} \label{time_exponential_rank_r}
            Let $1 \le  r < R, 1 \le  k \le M_r$. For each $t \ge 0, $ $$ \lim \limits_{\varepsilon \downarrow 0} \mathbb P \left({ \dfrac{T^{r, \varepsilon}(x, S_l^r)}{\tau_k^{r, \varepsilon}} \ge t  }\right) = e^{-t}$$ uniformly in $x \in S_k^r$, for each $1 \le l \le M_r$, provided that $S_k^r \neq S_i $ for every $1 \le i \le M$ (i.e., $S_k^r$ contains more than one subcluster of rank zero)  and  $P_{kl}^{r, \varepsilon}$ is not identically zero.
        \end{theorem}

        \proof Suppose that $S_k^r$ is the union of at least two distinct clusters of rank $r-1$ (the general case is discussed at the end of the proof).  Let us fix $i$ such that $S_{i}^{r - 1} \subset S_k^r$. For $x \in S_i^{r-1}$, let  $\sigma_0^{x, \varepsilon} = 0$ and $\sigma_n^{x, \varepsilon}$ be the minimum between the time of the $n$th return to $S_i^{r-1}$ and the time of the first exit from $S_k^r$ by the process ${\bf X}_n^{x, r- 1, \varepsilon}$. For $ x \in S_i^{r-1}$, let $ \xi_n^{x, \varepsilon} $ be a sequence of random variables with the following distribution:$$ \begin{aligned} & \mathbb P \left({\xi_0 = 0}\right) = 1, \\
             &\mathbb P \left({\xi_n^{x, \varepsilon}} \in I  \right) \\ &=  {  \mathbb P \left({  \chi_{S_l^{r}}({\bf X}_{\sigma_{n - 1}}^{x, r-1, \varepsilon}) e(\varepsilon) +  \chi_{S_i^{r-1}}({\bf X}_{\sigma_{n - 1}}^{x, r-1, \varepsilon}) \sum_{m = \sigma_{n-1}^{x, \varepsilon}}^{\sigma_{n }^{x, \varepsilon} - 1} {\bf T}_{m + 1}^{x, r-1, \varepsilon}   \in I  \biggm  \vert  {\bf X}_1^{x, r, \varepsilon} \in S_l^r   }\right)  }, n \ge 1\end{aligned}, $$  where $e(\varepsilon)$ is a function to be specified later.
             Thus, $\xi^{x,\varepsilon}_n$ measures the time spent between successive visits to $S^{r-1}_i$ (up to the exit from $S^r_k$), conditioned on exiting to $S^r_l$. However, if ${\bf X}_{n}^{x, r-1, \varepsilon}$ reaches $S^r_l$ by the time $\sigma_{n - 1}$, then $\xi^{x,\varepsilon}_n$ is equal to a constant $e(\varepsilon)$.
             For $ x \in S_i^{r-1}$, let  $\zeta_n^{x, \varepsilon}$ be a sequence of random variables with the following distribution: $$ \begin{aligned} &\mathbb P \left({\zeta_n^{x, \varepsilon}} = 1  \right) = \mathbb P \left( {{\bf X}_{\sigma_{n }}^{x, r-1, \varepsilon} \in S_i^{r-1}  \vert {\bf X}_1^{x, r, \varepsilon} \in S_l^r } \right),  n \ge 0,  \\ & \mathbb P \left({\zeta_n^{x, \varepsilon}} = 0 \right) = \mathbb P \left( {{\bf X}_{\sigma_{n }}^{x, r-1, \varepsilon} \in S_l^r  \vert {\bf X}_1^{x, r, \varepsilon} \in S_l^r } \right), n \ge 0.  \end{aligned}  $$  Then we note that $$ 
            \mathbb P \left({ \dfrac{T^{r, \varepsilon}(x, S_l^r)}{\tau_k^{r, \varepsilon}} \ge t  }\right) \\
        = {  \mathbb P \left({ ( {\tau_k^{r, \varepsilon}} )^{-1} \displaystyle \sum_{n = 0}^{\infty} \xi_{n+1}^{x, \varepsilon} \zeta_n^{x, \varepsilon} \ge t  }\right)  }.  
 $$ We will now verify that $\xi_{n}^{x, \varepsilon},  \zeta_n^{x, \varepsilon} $ satisfy the assumptions of  Lemma \ref{abstrac_exp_lemma} with certain $e(\varepsilon)$ and $r(\varepsilon)$ and also show that  $ {\tau_k^{r, \varepsilon}} \sim (e(\varepsilon))^{-1}r(\varepsilon)$.

From \eqref{abstract_chain_asymptote_lemma_initial_formula_1} and Lemma \ref{abstract_chain_lemma_2}, it follows that $$ \begin{aligned}
     \mathbb P \left({\zeta_1^{x, \varepsilon}} = 0 \right) &= \dfrac{\mathbb P \left( {{\bf X}_{\sigma_{1 }}^{x, r-1, \varepsilon} \in S_l^r   } \right)}{\mathbb P \left({ {\bf X}_1^{x, r, \varepsilon} \in S_l^r}\right)}   \\ &\sim   ({\lambda^{r-1}(i)})^{-1
 } \displaystyle \sum \limits_{\substack{j: S_j^{r- 1} \subset S_k^r , \\ j': S_{j'}^{r- 1} \not \subset S_k^r }}\lambda^{r-1}(j) P^{r- 1, \varepsilon}_{jj'} \end{aligned} $$ uniformly in $x \in S_i^{r-1}$ as $\varepsilon \downarrow 0$, where $\lambda^{r-1}$ is the invariant probability measure on $\{ j : S_j^{r-1} \subset S_k^r\}$ associated with the transition probabilities $P^{r-1, 0}$. After conditioning on $\mathcal F_{\sigma_n}^{x, \varepsilon}$ and using the Markov property, we obtain that $$  \mathbb P \left({\zeta_{n + 1 }^{x, \varepsilon}} = 0 \vert {\zeta_{n}^{x, \varepsilon}} = 1 \right) \sim   ({\lambda^{r-1}(i)})^{-1
 } \displaystyle \sum \limits_{\substack{j: S_j^{r- 1} \subset S_k^r , \\ j': S_{j'}^{r- 1} \not \subset S_k^r }}\lambda^{r-1}(j) P^{r- 1, \varepsilon}_{jj'}  $$ uniformly in $x \in S_i^{r-1}, n \ge 0$ as $\varepsilon \downarrow 0$. Let us denote the function on the right hand side of the asymptotic equivalence above by $r(\varepsilon)$. Observe that $r(\varepsilon) \rightarrow 0$ as $\varepsilon \downarrow 0$ since each term in the sum goes to zero by the construction of the hierarchy.

 Using arguments similar to those that led to \eqref{time_theorem_formula_1}, it is easy to show that $$ \lim \limits_{\varepsilon \downarrow 0} \dfrac{\mathbb E \left({ {\bf T}_1^{x, r-1, \varepsilon} \chi_{S_l^r} \left({{\bf X}_{1 }^{x,r,  \varepsilon}}\right)    }\right)}{\tau_j^{r - 1,\varepsilon} P_{kl}^{r, \varepsilon}} 
     = 1 $$ uniformly in $x \in S_j^{r-1}$ for every $j$ such that $S_j^{r-1} \subset S_k^r$.  (The only difference from \eqref{time_theorem_formula_1} is that here we evaluate the expectation of the exit time from a cluster of rank $r-1$ rather than rank zero.) 
     It is not difficult to show that, since the limiting process ${\bf X}_{{n }}^{x, r-1, 0}$ is ergodic, the expected time spent in $S_j^{r-1}$ between visits to $S_i^{r-1}$ (on the event that the process exits the cycle of rank $r$ to $S^r_l$) is asymptotically equivalent to the product of $\tau_j^{r-1, \varepsilon} P^{r, \varepsilon}_{kl}$ and the limit of the expected number of visits to $S_j^{r-1}$ between visits to $S_i^{r-1}$, i.e., $({\lambda^{r-1}(i)})^{-1}\lambda^{r-1}(j)$. 
     Summing over $j$ and using the Markov property, we obtain that $$ \lim \limits_{\varepsilon \downarrow 0} \dfrac{\mathbb E  \left({\displaystyle \sum_{m =0}^{\sigma_{1 }^{x, \varepsilon} - 1} {\bf T}_{m + 1}^{x, r-1, \varepsilon} \chi_{S_l^r} \left({{\bf X}_{1 }^{x,r,  \varepsilon}}\right)  }\right)}{ P_{kl}^{r, \varepsilon} ({\lambda^{r-1}(i)})^{-1
 }   \displaystyle \sum \limits_{\substack{j: S_j^{r- 1}  \subset S_k^r  }}\lambda^{r-1}(j) \tau^{r- 1, \varepsilon}_{j}  } = 1   $$  uniformly in $x \in S_i^{r-1}$. Let  $$e(\varepsilon) =  ({\lambda^{r-1}(i)})^{-1
 }   \displaystyle \sum \limits_{\substack{j: S_j^{r- 1}  \subset S_k^r  }}\lambda^{r-1}(j) \tau^{r- 1, \varepsilon}_{j}. $$ Then after conditioning on $ \mathcal F_{\sigma_{n}}^{x, \varepsilon}$, we can conclude that $$ \lim \limits_{\varepsilon \downarrow 0 }     \dfrac{\mathbb E \left({ \xi_{n+1}^{x, \varepsilon} \vert  \mathcal {F}^{x, \varepsilon}_{\sigma_n} }\right) }{ e(\varepsilon) } = 1  $$ uniformly in $x \in S_i^{r-1}, n \ge 0$, and the elements of the probability space. 

 Using arguments similar to those that led to \eqref{time_theorem_formula_2}, it is easy to show that $$ \mathbb E \left({ \left({ {\bf T}_1^{x, r- 1,  \varepsilon}}\right)^2 \chi_{S_l^r} \left({{\bf X}_{1 }^{x,r,  \varepsilon}}\right) }\right) \le C \cdot (\tau_j^{r - 1, \varepsilon})^2 P_{kl}^{r, \varepsilon}   $$ for some constant $C$, for each $j $ such that $S_j^{r-1} \subset S_k^r$, and each $ x \in S_j^{r-1}$ and small $\varepsilon$. Again, since the limiting process ${\bf X}_{{n }}^{x, r-1, 0}$ is ergodic, it can be shown that $$  \mathbb E \left({( \xi_{n+1}^{x, \varepsilon})^2 \vert  \mathcal {F}^{x, \varepsilon}_n }\right) \le C (e(\varepsilon))^2 $$ for some constant $C$, small enough $\varepsilon$, and each $x \in S_i^{r-1}, n \ge 0$, and each element of the probability space. 
 
 Finally, let $j$ be such that $S_j^{r-1} \subset S_k$ and $x \in S_j^{r-1}, y \in S_k^r$. Then, by Theorem \ref{prelim_things_abt_rank_r_times}, $ \tau_j^{r - 1, \varepsilon} \sim \displaystyle \sum_{a : S_a \subseteq S_j^{r-1}} \tau_a^{\varepsilon} \mathbb E(N^{\varepsilon}(x, a, r - 1))$. Moreover, by Lemma \ref{rank_r_expectation_lemma}, $$ \mathbb E(N^{\varepsilon}(y, a, r )) \sim   \lambda^{r-1}(j) \left({\displaystyle \sum \limits_{\substack{j: S_j^{r- 1} \subset S_k^r , \\ j': S_{j'}^{r- 1} \not \subset S_k^r }}\lambda^{r-1}(j) P^{r- 1, \varepsilon}_{jj'} } \right)^{-1} \mathbb E(N^{\varepsilon}(x, a, r - 1))  $$ for every $a$ such that $S_a \subset S_k^r$.  Hence, $$ r(\varepsilon)^{-1}e(\varepsilon) \sim \displaystyle \sum_{a : S_a \subset S_k^{r}} \tau_a^{\varepsilon} \mathbb E(N^{\varepsilon}(y, a, r )) \sim \tau_k^{r, \varepsilon}.  $$ The result for the case when $S_k^r$ contains at least two distinct clusters of rank $r-1$ now follows from Lemma \ref{abstrac_exp_lemma}. 
 
 If $S_k^r$ does not contain at least two distinct clusters of rank $r$, then we look at a smaller rank $r'$ such that $S_{k}^r = S_{k'}^{r'} $ for some $k'$ and $S_{k'}^{r'}$ is the union of at least two distinct clusters of rank $r'-1$. To conclude the result we note that in all the arguments used above we can replace $S_l^r$ with  $ \bigcup_{l' \in G} S_{l'}^{r'}$  and $P_{kl}^{r, \varepsilon}$ with $\sum_{l' \in G} P_{k'l'}^{r', \varepsilon}$ where $G$ is any suitable set of indices.  
 \qed \\

We now have all the information we need to claim that assumptions (a)-(f) from Section \ref{Reduction to a process with nearly independent transition times and transition probabilities} hold for every rank $r< R$. Let us collect and summarize these results in the following theorem.

\begin{theorem} \label{main_theorem_section_6}
    If assumptions (a)-(f) from Section \ref{Reduction to a process with nearly independent transition times and transition probabilities} are satisfied by $ \left({ {\bf X}_{n}^{x, \varepsilon} , {\bf T}_{n}^{x, \varepsilon} }\right) $ then they are also satisfied by $\left({ {\bf X}_{n}^{x, r, \varepsilon} , {\bf T}_{n}^{x, r, \varepsilon} }\right)$ for every $0 \le r < R.$
\end{theorem}
\proof The fact that assumption (a) holds for all ranks was discussed in the very beginning of this section since it was necessary to build the hierarchy of clusters. 

Note that we have only defined the functions $\tau_k^{r, \varepsilon}$ up to asymptotic equivalence. So let us fix points $x_i \in S_i$ for each $1 \le i \le M$ and let \begin{equation} \label{main_theo_sec_6_formula_1}
     \tau_k^{r, \varepsilon} = \displaystyle \sum \limits_{ i : S_i \subseteq S_k^r } \tau_i^{\varepsilon} \cdot \mathbb E \left({N^{\varepsilon}(x_i, i, r)}\right) 
\end{equation} for $1 \le  k \le M_r, 0 \le r < R $. Then assumptions (b) and (c) follow for each rank from Theorem \ref{prelim_things_abt_rank_r_times} and Theorem \ref{prelim_things_abt_rank_r_times_var}, respectively. Assumptions (d), (e) are true for all ranks as a consequence of Theorem \ref{time_exponential_rank_r}. 

Complete Asymnptotic Regularity (assumption (f)) follows for every rank due to Lemma \ref{appendix_lemma_1} after using the asymptotics of $  \mathbb{E} \left({N^{\varepsilon}(x_i, i, r)}\right) $ for each $i$ from Lemma \ref{rank_r_expectation_lemma}.
\qed
\\

\begin{remark}
    \label{main_theorem_section_6_remark} From \eqref{main_theo_sec_6_formula_1}  it follows that $\tau_k^{r, \varepsilon}$ and $\tau_l^{s, \varepsilon}$ are asymptotically comparable functions for any $0 \le r, s < R, 1 \le k \le M_r, 1 \le l \le M_s$, i.e., there is a limit $ \lim_{\varepsilon \downarrow 0} (\tau_k^{r, \varepsilon}/ \tau_l^{s, \varepsilon}) \in [0 , \infty]$. Indeed, by Lemma \ref{rank_r_expectation_lemma} it is not difficult to show that $\tau_k^{r, \varepsilon}/\tau_l^{s, \varepsilon} $ is  asymptotically equivalent to a ratio of linear combinations of certain functions of $\varepsilon$, where each such function is of the form $P_{i_1j_1}^{\varepsilon} \cdots P_{i_{\alpha}j_{\alpha}}^{\varepsilon} \tau_i^{\varepsilon}$. If $r \neq s$, then the number of factors of  the form $P_{ij}^{\varepsilon}$ that appear in these products in the numerator of the ratio may be different from that in the denominator. However, we can artificially introduce factors $P_{ij}^{\varepsilon}$ that have positive limits to conclude that  $$ \lim \limits_{\varepsilon \downarrow 0} \dfrac{P_{i_1j_1}^{\varepsilon} \cdots P_{i_{\alpha}j_{\alpha}}^{\varepsilon} \tau_i^{\varepsilon}}{P_{i'_1j'_1}^{\varepsilon} \cdots P_{i'_{\beta}j'_{\beta}}^{\varepsilon} \tau_{i'}^{\varepsilon}} \in [0, \infty] $$ even if $\alpha \neq \beta$, due to complete asymptotic regularity. The claim now follows by  Lemma~\ref{appendix_lemma_1}. 
\end{remark} 

\section{Ergodic Theorems for Semi-Markov Processes} \label{ergodic_section}

Let us consider a general semi-Markov process $X_t^{x, \varepsilon}$ (associated with the Markov chain $\left({ {\bf X}_n^{x, \varepsilon} }, {\bf T}_n^{x, \varepsilon} \right)$)
that satisfies assumptions (a)-(c) and (f) in section \ref{Reduction to a process with nearly independent transition times and transition probabilities}.
This need not be the same as our original process, but rather we would like to get general results that would be applicable under some of the same assumptions. Recall that $P_{ij}^{0}$  are the limits of the quantities $P_{ij}^{\varepsilon}$ as $\varepsilon \downarrow 0.$ Recall also that assumption (f) ensures the existence of the limits $\lim \limits_{\varepsilon \downarrow 0} {\tau_i^{\varepsilon}} / {\tau_j^{\varepsilon}} \in [0, \infty]$  for every $i, j \in \{1, \dots, M\}.$  After rescaling time by dividing by the largest expectation, we may suppose that the limits $\tau_i = \lim \limits_{\varepsilon \downarrow 0 } \tau_i^{\varepsilon}$ exist for every $ 1 \le i \le M, $ and that there is $i_0$ such that $\tau_{i_o} = 1$ and $0 \le \tau_i \le 1$ for every $1 \le i \le M.$ 

Lemma \ref{no_skip_large_interval} and Theorem \ref{ergodic_theorem_first} are proved in a simplified setting where the set $\{1, \dots, M\}$ forms one ergodic class with respect to the transition probabilities $P_{ij}^{0}$. We do this in order to illustrate the main ideas behind these results. Lemma \ref{no_skip_large_interval_small_escape_prob} and Theorem \ref{ergodic_theorem_second} are the extensions of these results to the case that is of more interest to us. 

\begin{lemma} \label{no_skip_large_interval} Suppose there is one ergodic class with respect to the transition probabilities $P_{ij}^0.$ Suppose that assumptions (a)-(c) and (f) from Section \ref{Reduction to a process with nearly independent transition times and transition probabilities} hold. We use the convention that $\tau_{i_0} = 1$ and $0 \leq \tau_i \leq 1$ for all $i$, as above. Then, for every $\delta > 0, $ there are $L,   \varepsilon_0 > 0$ such that 
    $$ \mathbb P \left({ X_s^{x, \varepsilon} = X_t^{x, \varepsilon} \ {\rm for \  every \ } s \in [t, t + L] }\right) < \delta $$ for every $t \ge 0, \varepsilon \le \varepsilon_0$, and $x \in S.$
\end{lemma}
\proof 
Let us first fix $1 \le i \le M$. Let $\boldsymbol{\sigma}_0^{x, \varepsilon} = 0$ and, for each $n \geq 1$, let $\boldsymbol{\sigma}^{x, \varepsilon}_n$  be the time of the first arrival to $S_{i}$ by $X_t^{x, \varepsilon}$ after time $\boldsymbol{\sigma}^{x, \varepsilon}_{n-1}$. 
(We use this bold-face $\boldsymbol{\sigma}$ notation in order to distinguish the stopping times that take real values from integer-valued stopping times that we used earlier.)
Notice that there is a constant  $a > 0$ such that $\mathbb E \boldsymbol{\sigma}_1^{x, \varepsilon} \ge a$ for every $x \in S_{i}$ and small enough $\varepsilon$. Since the set $\{1, \dots, M\}$ forms one ergodic class under $P_{ij}^0$, and since assumption (c) holds, it is not difficult to see that the variance of $\boldsymbol{\sigma}_1^{x, \varepsilon}$ is bounded above uniformly in $x \in S$ and small enough $\varepsilon$. Therefore, there are positive constants $C_1, C_2, C_3$ with $ C_2 \in (0, 1)$ such that $$ \mathbb P \left({ \boldsymbol{\sigma}_n^{x, \varepsilon}  \le  1 }\right) \le C_1 \cdot C_2^{ n } $$ $$ \mathbb P \left({ \boldsymbol{\sigma}_1^{x, \varepsilon} \ge L }\right)  \le \dfrac{C_3}{L^2} $$  for all $ n \ge 0, L > 0, x \in S$ and small enough $\varepsilon.$

For $k \ge 0, $ let $\mathcal A_{k, t}^{x, \varepsilon}$ be the event that the most recent arrival to $S_i$ by the process $X_t^{x, \varepsilon}$ on or before time $t$ occurred in the interval $(t - (k+1), t - k  ].$ For $n \ge 0, $ let $\mathcal B_{k, n, t}^{x, \varepsilon} \subseteq A_{k, t}^{x, \varepsilon}$ be the sub-event that $S_{i}$ was visited on exactly $n$ distinct occasions in the interval $(t - (k+1), t - k  ].$  
 Let $x \in S_i$. Then $$ \begin{aligned}
\mathbb P & \left({ X_s^{x, \varepsilon} = X_t^{x, \varepsilon} \ {\rm for \  every \ } s \in [t, t + L] }\right) \\ &=  \sum_{k = 0}^{\infty}    \mathbb P \left({\left\{{ X_s^{x} = X_t^{x} \ {\rm for \  every \ } s \in [t, t + L] }\right\} \bigcap  \mathcal A^x_{k, t}  }\right)   \\ &= \sum_{k = 0}^{\infty}  \sum_{n = 0}^{\infty}  \mathbb P \left({\left\{{ X_s^{x} = X_t^{x} \ {\rm for \  every \ } s \in [t, t + L] }\right\} \bigcap  \mathcal B^x_{k, n, t}  }\right) \\ &\le \sum_{k = 0}^{\infty} \sum_{n = 0}^{\infty} \sup \limits_{x \in S} \mathbb P \left({ \boldsymbol{\sigma}_n^{x, \varepsilon} \le 1 }\right)  \sup \limits_{x \in S_{i_0}} \mathbb P \left({ \boldsymbol{\sigma}_1^{x, \varepsilon} \ge L + k }\right)   \\ &\le \sum_{k = 0}^{\infty} \sum_{n = 0}^{\infty}   C_1 \cdot C_2^{ n } \cdot \dfrac{C_3}{ (L + k)^2 } \to 0 \ {\rm as} \ L \to \infty. 
\end{aligned} $$  \qed
\\
 Suppose  again that  $\{1, \dots, M\}$ 
 forms one ergodic class with respect 
 to the transition probabilities $P_{ij}^0$. Then there is a unique invariant measure $\lambda$ on $\{1, \dots , M\}$. Consider the measure $\mu $ on $\{ 1, \dots , M \}$ given by $$ \mu(i) =\dfrac{ \lambda(i) \tau_{i}  }{ \displaystyle \sum_{j = 1}^{M} \lambda(j) \tau_{j}   }. $$ 

\begin{theorem} \label{ergodic_theorem_first}
    Suppose that $ t = t(\varepsilon) \rightarrow \infty$ as $\varepsilon \downarrow 0$. Suppose that assumptions (a)-(f) from Section \ref{Reduction to a process with nearly independent transition times and transition probabilities} hold. We use the convention that $\tau_{i_0} = 1$ and $0 \leq \tau_i \leq 1$ for all $i$. If the states $\{1, \dots, M\}$ form one ergodic class with respect to the transition probabilities $P^0_{ij}$, then $$ \lim \limits_{ \substack{ \varepsilon \downarrow 0 } } \mathbb P \left({  X^{x, \varepsilon}_{t(\varepsilon)} \in S_i }\right) = \mu (i) $$ uniformly in $x \in S.$ 
\end{theorem} 

\proof  

If the claim is not true, we can extract sequences $\varepsilon_{n} \to 0, t_n \to \infty$, and $x_n\in S$ such that $ \left|{ \mathbb P \left({ X_{t_n}^{x_n, \varepsilon_n} \in S_{j_0}}\right) - \mu(j_0)}\right| > 3 \delta   $ for some $j_0$ and $\delta > 0.$  By assumption (c), the collection  $ \{T^{\varepsilon}(x, S_j) : x \in S_i,  \varepsilon > 0  \} $ of random variables is tight for each $1\le i, j \le M$ for which $P_{ij}^{\varepsilon} > 0.$ Therefore, we can find a further subsequence (also denoted $(\varepsilon_n, t_n , x_n)$) such that $T^{\varepsilon_n} (y_i, S_j)$ converges in distribution to a random variable $T_{ij}$ as $n \to \infty$ for every $1 \le i, j \le M$, where we note that the distribution of $T_{ij}$ does not depend on $j$ or on the choice of $y_i \in S_i$, as a consequence of assumption (d). Thus we will write $T_i = T_{ij}$. Due to  assumption (c), we have $\mathbb E (T_{i}) = \tau_i.$ 

The family of Markov renewal processes $X_t^{i, 0}$ on $\{1, \dots ,  M\}$ defined by the transition probabilities $P_{ij}^0, $ the sojourn times $T_{i} $ and initial point $X_0^{i, 0} = i,$ admits the limiting distribution $\mu$ defined above, since the distributions $T_{i}$ are non-arithmetic by assumption (e).
Namely, by the ergodic theorem for Markov renewal processes (on a finite state space), 
there exists $t_0 > 0$ such that 
$$
|\mathbb P \left({ X_t^{i, 0} = j} \right)
- \mu(j)| <
\delta 
$$
for every $t \ge t_0$, $1 \leq i,j \leq M$  (\cite{AS}, \cite {CINL}). 
Due to assumptions (a),  (d), and (e), for any fixed $L >0$, there is $N$ such that  $$ \left|{ \mathbb P \left({ X_{t}^{x, \varepsilon_n} \in S_j} \right) - \mathbb P \left({ X_{t}^{i, 0} = j} \right)  }\right| < \delta $$ for every $n \ge N, x \in S_i, 1 \le i, j \le M$ and $t \in [t_0 , t_0 + L].$ If $\boldsymbol{\sigma}_n$ is defined as the first renewal after time $t_n - t_0 - L$  of the process $X_t^{x_n, \varepsilon_n},$ then Lemma \ref{no_skip_large_interval} guarantees that $L$ can be chosen large enough so that this renewal occurs before time $t_n - t_0$ with probability at least $1 - \delta$ for large enough $n.$  Therefore, $$ \begin{aligned}
    &\left|{ \mathbb P \left({ X_{t_n}^{x_n, \varepsilon_n} \in S_{j_0}  }\right) - \mu(j_0) }\right| \\ &<  \delta + \int \limits_{t_n - t_0 - L}^{t_n - t_0} \int \limits_{S} \left|{  \mathbb P \left({X_{t_{n} - t}^{y, \varepsilon_n} \in S_{j_0}}\right)  - \mu(j_0)   }\right| \mathbb P \left({\boldsymbol{\sigma}_n \in dt ,X_{t}^{x_n, \varepsilon_n} \in dy}\right)   \\ &\le 3 \delta . \end{aligned} $$ 
    This brings us to a contradiction and completes the proof of the lemma.
    \qed
    \\

 Now let us consider a slightly more general setup for the process $X_t^{x, \varepsilon}$. Suppose now that there is $1 < M' < M$ such that $\{1, \dots, M'\}$ forms one ergodic class with respect to $P_{ij}^0$. 
 As before, let us suppose that $\tau_{i_0} = 1$ for some $ 
1 \le i_0 \le M'$ and $ \tau_i \le 1 $ for every $1 \le  i \le M'$. Then it can be shown, similarly to Lemma \ref{no_skip_large_interval}, that renewals take place fairly often. We have a slightly stronger result that is valid even in the presence of transient states. Let $S' = S_1 \bigcup \cdots \bigcup S_{M'}$.

 \begin{lemma} \label{no_skip_large_interval_small_escape_prob} Suppose that assumptions (a) - (c) and (f) from section \ref{Reduction to a process with nearly independent transition times and transition probabilities} hold. If the states $\{1, \dots, M'\} $ form one ergodic class with respect to $P_{ij}^0$, then for every $\delta > 0,$ there are $L, \varepsilon_0 > 0$ such that $$ \mathbb P \left({ X_t^{x, \varepsilon} \in S_i, X_s^{x, \varepsilon} \in S_i  \ {\rm for \ every} \ s \in [t , t + L ] }\right)  < \delta,  $$
 $$ \begin{aligned} \mathbb P \left({X_t^{x, \varepsilon} \in S_i, \ {\rm the \ next \ state \ visited \ is \ in \ } S \setminus S' }\right) < \delta  \end{aligned} $$
 for every $1 \le i \le  M', x \in S, \varepsilon \le \varepsilon_0, $ and $ t \ge 0.$ \end{lemma}
 \proof Let us first fix $
 1 \le i \le M'.$ Let  $\boldsymbol{\sigma}_0^{x, \varepsilon}= 0 $  and let $ \boldsymbol{\sigma}_n^{x, \varepsilon}  $ be the time of the first new arrival to $S_i \bigcup (S \setminus S') $ after time $\boldsymbol{\sigma}_{n-1}^{x, \varepsilon}$ by the process $X_t^{x, \varepsilon}$. 
 As in the proof of Lemma \ref{no_skip_large_interval}, for $k \ge 0, $ let $\mathcal A_{k, t}^{x, \varepsilon}$ be the event that the most recent arrival to $S_i$ at or prior to time $t$ occurred in the interval $(t - (k +1), t - k].$ Let  $\mathcal B_{k,n, t}^{x, \varepsilon} \subseteq \mathcal A_{k, t}^{x, \varepsilon}$ be the sub-event that $S_i$ was visited on exactly $n$ distinct occasions in the interval $(t - (k +1), t - k].$ 
 
 It is not difficult to see that $\boldsymbol{\sigma}_1^{x, \varepsilon}$ has expectation bounded away from $0$ uniformly in $x \in S_i$ for small $\varepsilon.$  Therefore, as in Lemma \ref{no_skip_large_interval}, there are constants $C_1 > 0,  C_2  \in (0, 1)$ such that $$ \sup \limits_{x \in S} \mathbb P \left({ \boldsymbol{\sigma}_n^{x, \varepsilon} \le  1 }\right) \le C_1 \cdot C_2^{n}$$ for small enough $\varepsilon$. Due to assumption (c) and the Chebyshev inequality, there is a constant $C_3 >0$ such that $$ \sup \limits_{\substack{x \in S_i \\ 0 \le j \le M}}  \mathbb P \left({T^{\varepsilon}(x, S_j) \ge L }\right) \le \dfrac{C_3}{L^2}  $$ for every $L > 0$ and small enough $\varepsilon$.
 
 Then for every $x \in S$ and small $\varepsilon$, $$ \begin{aligned} & \mathbb P  \left({ X_t^{x, \varepsilon} \in S_i, X_s^{x, \varepsilon} \in S_i  \ {\rm for \ every} \ s \in [t , t + L ]}\right) \\ & =  \sum_{k = 0 }^{\infty} \sum_{n = 1}^{\infty}  \mathbb P \left({ \mathcal B_{k,n, t}^{x, \varepsilon} \bigcap \left\{{ X_t^{x, \varepsilon} \in S_i, X_s^{x, \varepsilon} \in S_i  \ {\rm for \ every} \ s \in [t , t + L ]}\right\} }\right) \\ & \le 
 \sum_{k = 0 }^{\infty} \sum_{n = 1}^{\infty}   \sup \limits_{x \in S} \mathbb P \left({ \boldsymbol{\sigma}_{n}^{x, \varepsilon}\le 1 }\right)  \sum_{1 \le j \le M} \sup \limits_{ \substack{ x \in S_i \\}} P^{\varepsilon} \left({ x,  S_j }\right) \mathbb P \left({  {T^{\varepsilon}(x, S_j) \ge L + k } }\right)  \\  & \le 
 \sum_{k = 0 }^{\infty} \sum_{n = 1}^{\infty}   \sup \limits_{x \in S} \mathbb P \left({ \boldsymbol{\sigma}_{n}^{x, \varepsilon}\le 1 }\right) \sup \limits_{ \substack{ x \in S_i \\0 \le  j \le M}}  \mathbb P \left({  {T^{\varepsilon}(x, S_j) \ge L + k } }\right)       \\ & \le   \sum_{k = 0 }^{\infty} \sum_{n = 1}^{\infty}   C_1 \cdot C_2^{n} \dfrac{C_3}{(L + k)^2}   \\ & < {\delta} \end{aligned}  
$$  for large enough $L.$ 

Furthermore, since the states $\{1, \dots, M'\} $ form one ergodic class with respect to $P_{ij}^0$, for any $\delta' >0, P^{\varepsilon} \left({x, S_j}\right) < \delta'$ for every $x \in S_i, M' + 1 \le  j \le M $ if $\varepsilon $ is small. Therefore, for every $x \in S$ and small $\varepsilon$, $$ \begin{aligned} \mathbb P &\left({X_t^{x, \varepsilon} \in S_i, \ {\rm the \ next \ state \ visited \ is \ in \ } S \setminus S'}\right) \\ & = \sum_{k = 0}^{\infty} \sum_{n = 1}^{\infty} \mathbb P \left({ \mathcal B_{k, n, t}^{x, \varepsilon} \bigcap \left\{{ X_t^{x, \varepsilon} \in S_i, \ {\rm the \ next \ state \ visited \ is \ in \ } S \setminus S')   }\right\} }\right) \\ & \le \sum_{k = 0}^{\infty} \sum_{n = 1}^{\infty} \sup \limits_{x \in S} \mathbb P \left({ \boldsymbol{\sigma}_{n}^{x, \varepsilon}\le 1 }\right) 
 \sup \limits_{\substack{ x \in S_i \\ M' + 1 \le j \le M }} P^{\varepsilon} \left({x, S_j}\right) \mathbb P \left({ T^{\varepsilon}(x, S_j) \ge k }\right)   \\ & \le   \sum_{n = 1}^{\infty}  \delta'  C_1 \cdot C_2^{n} \cdot \  \left(1  + {\sum_{k = 1}^{\infty} \dfrac{C_3}{k^2}}\right)  \\ & < \delta, \end{aligned} $$ where the last inequality hold for an appropriate $\delta'$. \qed
    \\

    Consider the same setting discussed prior to Lemma \ref{no_skip_large_interval_small_escape_prob}. Suppose again that $\{1, \dots, M'\}$ forms one ergodic class with respect to $P_{ij}^0$. After rescaling time, there is $1 \le i_0 \le M'$ such that $ \tau_{i_0} =1$ and $\tau_i \le 1$ for every $1 \le  i \le M'$. Let us redefine the measure $\mu$ as $$ \mu(i)  = \dfrac{\lambda(i) \tau_i}{\displaystyle \sum_{j = 1}^{M'} \lambda(j) \tau_j}$$ for $1 \le  i \le M'$, where $\lambda$ is the invariant measure on $\{1, \dots, M'\}$ associated with $P_{ij}^0$. 
    
    Let $$\eta^{x, \varepsilon} = \inf \left\{{ t \ge 0 : X_t^{x, \varepsilon} \not \in S' }\right\}$$ for $x \in S', \varepsilon > 0.$    

 \begin{theorem} \label{ergodic_theorem_second} Suppose that $t= t(\varepsilon) \to \infty$ as $\varepsilon \downarrow 0.$   Suppose  that assumptions (a) - (f) from Section \ref{Reduction to a process with nearly independent transition times and transition probabilities} hold and the states $\{1, \dots, M'\} $ form one ergodic class with respect to $P_{ij}^0$. Suppose that $ \mathbb P \left({ \eta^{x, \varepsilon} \le t(\varepsilon)}\right) \to 0  $ as $\varepsilon \downarrow 0$ uniformly in $x \in S'$. Then \begin{equation} \label{ergodic_theorem_second_limit_form}
     \lim \limits_{\varepsilon \downarrow 0} \mathbb P \left({ 
 X_{t(\varepsilon)}^{x, \varepsilon} \in S_i}\right) = \mu(i)
 \end{equation} uniformly in $x \in S'.$
 Moreover, suppose that $ t_1(\varepsilon), t_2(\varepsilon)$  are two time scales that satisfy $t_1(\varepsilon) \le t_2(\varepsilon)$ for every $\varepsilon > 0$,  $t_1(\varepsilon) \to \infty $ as $\varepsilon \downarrow 0$ and $ \mathbb P \left({ \eta^{x, \varepsilon} \le t_2(\varepsilon)}\right) \to 0  $ as $\varepsilon \downarrow 0$ uniformly in $x \in S'$. Then the above limit holds uniformly in all time scales $t(\varepsilon)$ that satisfy  $t_1(\varepsilon) \le t(\varepsilon) \le t_2(\varepsilon)$. 
 \end{theorem}
 
  \proof Let's start with the first statement, and consider a particular function $t(\varepsilon)$. Consider an auxiliary process $\Tilde{X}_t^{x, \varepsilon}$ associated with the Markov chain $\left({ \Tilde{\bf X}_n^{x, \varepsilon}, \Tilde{\bf T}_n^{x, \varepsilon} }\right),$ governed by the Markov transition kernel $\Tilde{Q}^{\varepsilon}: \left({S'}\right) \times \mathcal B \left({S'}\right) \times  \mathcal B \left({[0, \infty)}\right) \to [0, 1],$ defined by, $$ \Tilde{Q}^{\varepsilon}\left({x, A \times I}\right) = \dfrac{Q^{\varepsilon}\left({x, A \times I}\right)}{ P^{\varepsilon} \left({x, S'}\right) }  = \mathbb P \left({ {\bf X}_1^{x, \varepsilon} \in A, {\bf T}_1^{x, \varepsilon} \in I \ \vert \ {\bf X}_1^{x, \varepsilon}  \in S' }\right). $$ Then  $\Tilde{X}_t^{x, \varepsilon}$ satisfies assumptions (a) - (f) of Section \ref{Reduction to a process with nearly independent transition times and transition probabilities} with the same functions $P_{ij}^{\varepsilon}$ and $\tau_i^{\varepsilon}$ serving as relevant asymptotics.   Therefore, Theorem \ref{ergodic_theorem_first} goes through for $\Tilde{X}_t^{x, \varepsilon}$ with the new limiting distribution $\mu$ defined above.

   For bounded times, the distributions of $X_t^{x, \varepsilon}$ and $\Tilde X_t^{x, \varepsilon}$ are close. Therefore, for any $\delta > 0, L >0,$ there are $t_0, \varepsilon_0 > 0$ such that $$ \left| { \mathbb P \left({X_t^{x, \varepsilon} \in S_j}\right) - \mu(j) } \right| < \delta$$ for all $\varepsilon \le \varepsilon_0, t \in [t_0, t_0 + L], x \in S_i $ and $ 1\le i, j, \le M'$.

  Suppose there are sequences $\varepsilon_n \to 0, x_n \in S'$ and $1 \le j_0 \le M'$ such that  $$\left|{ \mathbb P \left({X^{x_n, \varepsilon_n}_{t(\varepsilon_n)} \in S_{j_0}}\right) -  \mu(j_0) }\right| > 3 \delta $$ for some $\delta > 0.$ We know that \begin{equation} \label{ergodic_theore_second_formula_1} \mathbb P \left({ \eta^{x_n, \varepsilon_n} > t(\varepsilon_n) }\right)  \ge 1 - \delta \end{equation} for large enough $n. $ Let $\boldsymbol{\sigma}_n$ be the first renewal time after time $t(\varepsilon_n) - t_0 - L$ by the process $X_t^{x_n, \varepsilon_n}$. Consider the event $\mathcal A_n$ $ =  \left\{{\boldsymbol{\sigma}_n \le t(\varepsilon_n) - t_0}\right\} \bigcap $ $ \left\{{X_{\boldsymbol{\sigma}_n}^{x_n, \varepsilon_n} \in S'}\right\}$. Lemma \ref{no_skip_large_interval_small_escape_prob} and \eqref{ergodic_theore_second_formula_1} guarantee that $ \mathbb P \left({ \mathcal A_n  }\right) \ge 1 - 2 \delta$ for large $n$. Therefore, $$ \begin{aligned}
    &\left|{ \mathbb P \left({ X_{t(\varepsilon_n)}^{x_n, \varepsilon_n} \in S_{j_0}  }\right) - \mu(j_0) }\right| \\ & \le 2 \delta  +  \left|{ \mathbb P \left({ \mathcal A_n \bigcap \left\{{X_{t(\varepsilon_n)}^{x_n, \varepsilon_n} \in S_{j_0}}\right\}   }\right)  - \mu(j_0) \mathbb P \left({\mathcal A_n}\right) }\right|  \\ &\le 2 \delta + \int \limits_{t_n - t_0 - L}^{t_n - t_0} \int \limits_{S'} \left|{  \mathbb P \left({X_{t(\varepsilon_n) - t}^{y, \varepsilon_n} \in S_{j_0}}\right)  - \mu(j_0)   }\right| \mathbb P \left({\boldsymbol{\sigma}_n \in dt, X_{t}^{x_n, \varepsilon_n} \in dy}\right)  \\ &\le   3 \delta . \end{aligned}   $$    
   This brings us to a contradiction and establishes the result for a given time scale $t(\varepsilon)$. 
   
   Suppose now that $t_1(\varepsilon),  t_2(\varepsilon)$ are two time scales that satisfy $t_1(\varepsilon) \le t_2(\varepsilon)$ for every $\varepsilon > 0$,  $t_1(\varepsilon) \to \infty $ as $\varepsilon \downarrow 0$ and $ \mathbb P \left({ \eta^{x, \varepsilon} \le t_2(\varepsilon)}\right) \to 0  $ as $\varepsilon \downarrow 0$ uniformly in $x \in S'$. If the limit in (\ref{ergodic_theorem_second_limit_form})  does not hold uniformly in all time scales $ t_1(\varepsilon) \le t(\varepsilon) \le t_2(\varepsilon)$, then, by the definition of uniform convergence, we can identify a time scale that violates the first statement of the theorem, leading to a contradiction. 
 \qed
 \\

Consider the same setting as Theorem \ref{ergodic_theorem_second} without rescaling time, though. For a time scale $t = t(\varepsilon)$, let $E 
 = \{1 \le  i \le M' : \tau_i^{\varepsilon} \ll t(\varepsilon)\}, E' 
 = \{1 \le  i \le M' : \tau_i^{\varepsilon} \gg t(\varepsilon)\} $. 

 \begin{theorem} \label{ergodic_theorem_third} Consider a time scale $t= t(\varepsilon)$. Suppose that assumptions (a) - (f) from Section \ref{Reduction to a process with nearly independent transition times and transition probabilities} hold and that the states $\{1, \dots, M'\}$ form one ergodic class with respect to the transition probabilities $P_{ij}^0$ for some $1 < M' < M$. If $ E, E' \neq \emptyset $ and $ E \bigcup E' =  \{1, \dots, M'\}$, then for each $i \in E,$ there is a probability measure $\mu(i, \cdot)$ concentrated on $E'$  such that $$ \lim \limits_{\varepsilon \downarrow 0} \mathbb P \left({ X^{x, \varepsilon}_{t(\varepsilon) } \in S_j }\right) = \mu(i, j)$$ uniformly in $x \in S_i.$  Moreover, if $t_1(\varepsilon) \le t_2(\varepsilon)$ are two time scales associated with the same sets $E$ and $ E'$, then the limit holds uniformly in all time scales $t(\varepsilon)$ satisfying $t_1(\varepsilon) \le t(\varepsilon) \le t_2(\varepsilon)$.   
     
 \end{theorem} 

\proof 
 As in Theorem \ref{ergodic_theorem_second}, it is enough to prove the result for a particular time scale $t(\varepsilon)$. 
Let us fix an index $ i \in E.$ For $x \in S_i, $ let $$\sigma^{x, \varepsilon} = \min \left \{{n : {\bf X}_n^{x, \varepsilon} \in S_j \ {\rm for \ some} \ j \not  \in E }\right\}. $$  Since $\{1, \dots, M'\}$ forms one ergodic class with respect to $P_{ij}^0$, $\mathbb P \left({ \sigma^{x, \varepsilon} \ge n }\right)$  decays exponentially in $n$. Also, since $t(\varepsilon) \gg \tau_j^{\varepsilon}$ for every $j \in E$,  $$
    \lim \limits_{\varepsilon \downarrow 0} \mathbb P \left({ \displaystyle \sum_{n = 1}^{\sigma^{x, \varepsilon} } {\bf T}_n^{x, \varepsilon} \le t(\varepsilon), {\bf X}_{\sigma}^{x, \varepsilon}  \in S'   }\right)  = 1 
$$   uniformly in $x \in S_i.$ If $\nu(i, j)$ is the probability that the Markov chain on $\{1, \dots, M'\}$ governed by $P_{ij}^0$ and starting at $i$ hits $j$ when it first enters $E'$, then $$ \lim \limits_{\varepsilon \downarrow 0 } \mathbb P \left({ {\bf X}_{\sigma}^{x, \varepsilon} \in S_j }\right) = \mu(i, j)$$ uniformly in $x \in S_i$. Moreover, due to assumption (e), $\mathbb P \left({T^{\varepsilon}(x, S) \le t(\varepsilon) }\right) \to 0$ as $\varepsilon \downarrow 0$ uniformly in $x \in S_j, j \in E'$. Therefore, $$ \lim \limits_{\varepsilon \downarrow 0} \mathbb P \left({ X_{t(\varepsilon)}^{x, \varepsilon} \in S_j }\right) = \lim \limits_{\varepsilon \downarrow 0 } \mathbb P \left({ {\bf X}_{\sigma}^{x, \varepsilon} \in S_j }\right) = \mu(i, j) $$ uniformly in $x \in S_i$.
\qed \\

\section{Proof of Theorem \ref{prelimt2}} \label{times} 

Recall the time scales 
$\tau_k^{r, \varepsilon},1  \le k \le M_{r}$, $0 \le  r \le R$, introduced in Theorem \ref{prelim_things_abt_rank_r_times}.  These functions are comparable to each other due to Remark \ref{main_theorem_section_6_remark}.
Retaining one representative from each set of equivalent time scales, we get a collection that will serve as the time scales $\mathbf{t}_{1}^{\varepsilon}, \dots, \mathbf{t}_{n-1}^{\varepsilon}$, with $\mathbf{t}^\varepsilon_0 = 0$, $\mathbf{t}^\varepsilon_n = \infty$. 
Observe (see remark after Theorem \ref{prelimt2})  that these functions do not depend on $i$. However, only some of those time scales will be relevant when describing the behavior of the process starting at $x \in S_i$ with a given $i$.

Let $1 \le  i \le M.$ Suppose that $t = t(\varepsilon)$ is a positive function that, for every $1 \le  k \le M_r$, $0 \le r < R$, satisfies $ t(\varepsilon) \ll \tau_k^{r, \varepsilon} $ or $t(\varepsilon) \gg \tau_k^{r, \varepsilon}$. If $t(\varepsilon) \ll \tau_i^{0, \varepsilon}$, then define $r = -1$. Otherwise, let $r$ be the largest rank such that $t(\varepsilon) \gg \tau_{k}^{r, \varepsilon}$ and $S_{i} \subseteq S_{k}^{r}$ for some $1 \le k \le M_r$. In other words, $r$ is the first rank such that $t(\varepsilon)$ is asymptotically smaller than the expected time to exit the cluster of rank $r+1$ to which $S_{i}$ belongs. 

If $r = -1$, then  \begin{equation} \label{main_proof_formula_3}
    \lim \limits_{\varepsilon \downarrow 0} \mathbb P \left({X_{t(\varepsilon)}^{x, \varepsilon} \in S_i}\right) = 1
\end{equation}uniformly in $x \in S_i, $ as a simple consequence of assumptions (c) and (e).  In other words, the limiting measure corresponding to $0 \ll t(\varepsilon) \ll \tau^{0,\varepsilon}_i$ is concentrated at $i$. Note that $\mathbf{t}^\varepsilon_1(i) = \tau^{0,\varepsilon}_i $ and $\mu_{i,1}$ is concentrated at $i$, in the notation of Theorem~\ref{prelimt2}, but
the interval $(0, \tau^{0,\varepsilon}_i)$ may contain multiple intervals $( \mathbf{t}^\varepsilon_0, \mathbf{t}^\varepsilon_1)$, $( \mathbf{t}^\varepsilon_1, \mathbf{t}^\varepsilon_2)$(see remark after Theorem \ref{prelimt2}).   Furthermore, it can be shown that if $t_1(\varepsilon) \ll \tau_i^{0, \varepsilon}$, then \eqref{main_proof_formula_3} holds uniformly in all time scales $0 \le t(\varepsilon) \le t_1(\varepsilon)$.   

Suppose that $r = 0$ and $t(\varepsilon) \gg \tau_j^{0, \varepsilon}$ for every $j$ such that $S_j^0 \subseteq \pi^1(S_i)$. In this case, Theorem \ref{ergodic_theorem_second} applies to ensure the existence of a measure, which we will call $\nu^0(i,j)$, that satisfies \begin{equation} \label{main_proof_formula_4}
     \lim \limits_{\varepsilon \downarrow 0} \mathbb P \left({X_{t(\varepsilon)}^{x, \varepsilon} \in S_j^0}\right) = \nu^0(i,j)
\end{equation} uniformly in $x \in S_i.$ 

If $r = 0$ and $t(\varepsilon) \ll \tau_j^{0, \varepsilon}$ for some $S_j^0 \subseteq \pi^1(S_i)$, then Theorem \ref{ergodic_theorem_third} ensures the existence of a probability measure, again called $\nu^0(i, \cdot)$ concentrated on the set of $j$'s satisfying $S_j^0 \subseteq \pi^1(S_i)$, such that \begin{equation} \label{main_proof_formula_5}
    \lim \limits_{\varepsilon \downarrow 0} \mathbb P \left({X_{t(\varepsilon)}^{x, \varepsilon} \in S_j^0}\right) = \nu^0(i,j) 
\end{equation} uniformly in $x \in S_i$. Depending on the number and the magnitude of different time scales $\tau_j^{0, \varepsilon}$ that appear for those $j$ satisfying $S_j^0 \subseteq \pi^1(S_i)$, there may be several such distinct metastable distributions (these will play the role of $\mu_{i, 2}, \mu_{i, 3}$ etc. according to the notation in the statement of Theorem \ref{prelimt2}). The existence of all these distinct measures is a direct consequence of Theorem \ref{ergodic_theorem_third}.  Furthermore, if  $t_1(\varepsilon) \le  t_2(\varepsilon)$ are two time scales whose ratio is bounded away from zero and infinity that satisfy the assumptions imposed on $t(\varepsilon)$, then \eqref{main_proof_formula_4} (or \eqref{main_proof_formula_5}, as applies) holds uniformly in all time scales $\overline {t}(\varepsilon)  $ satisfying $t_1(\varepsilon) \le  \overline {t}(\varepsilon)  \le {t}_2(\varepsilon)$.

Suppose that $r \ge 1.$ We will use induction on the quantity $r$, which was defined as the smallest rank such that $t(\varepsilon)$ is asymptotically smaller than the expected exit time from the cluster of rank $r + 1$ to which $S_i$ belongs.  Therefore, we assume that  $\mathbb P \left({X^{x, \varepsilon}_{t(\varepsilon)} \in S_j}\right)$ has a limit for each $1 \le j \le M$ as $\varepsilon \downarrow 0$, uniformly in $x \in S_i$ for every $1 \le i,j \le M$ if $t(\varepsilon)$ is asymptotically small enough to produce an $r$ that is strictly smaller than the one being considered. We also assume that if $t_1(\varepsilon) \le t_2(\varepsilon)$ are two time scales whose ratio is bounded away from zero and infinity that satisfy the assumptions imposed on $t(\varepsilon)$, then all such limits exist uniformly in time scales $\overline t (\varepsilon)$ satisfying $t_1(\varepsilon) \le \overline t(\varepsilon) \le t_2(\varepsilon)$.   

\subsubsection*{ Case 1: $t(\varepsilon) \gg \tau_k^{r, \varepsilon}$ for every $S_k^r \subseteq \pi^{r + 1}(S_i)$.}

 Recall the process $X_t^{x, r, \varepsilon}$ that tracks our original process along clusters of rank $r$, introduced in Definition \ref{rank_r_definition}. 
 Theorem \ref{ergodic_theorem_second} applies to this process, therefore,  there is a measure $\nu^{r}(i, \cdot)$ concentrated on those indices $k$ for which $ S_k^r \subseteq \pi^{r+1}(S_i) $ such that

\begin{equation} \label{main_proof_formula_1} \lim  
    \limits_{\varepsilon \downarrow 0} \mathbb P \left({ X_{t(\varepsilon)}^{x, \varepsilon} \in S_k^{r}}\right) = \lim \limits_{\varepsilon \downarrow 0} \mathbb P \left({ X_{t(\varepsilon)}^{x,r,  \varepsilon} \in S_k^{r}}\right) = \nu^{r} (i, k) 
\end{equation} uniformly in $x \in \pi^{r + 1}(S_{i}).$  Consider a time scale $s(\varepsilon)$ such that $s(\varepsilon) \ll \tau_k^{r, \varepsilon}$ for every $k $ such that $S_k^r \subseteq \pi^{r+1}(S_i), \nu^r(i, k) > 0, $ and $s(\varepsilon) \gg \tau_{k'}^{r', \varepsilon}$ for every $k', r'$ such that $0 \le r' < r$ and  $S_{k'}^{r'} \subsetneq S_k^r.$ 
In other words, $s(\varepsilon)$ is sufficiently small so that the process doesn't get out of any cluster $S^r_k$ that has positive measure, yet large enough to make many transitions between proper subclusters of $S^r_k$ if the process starts in $S^r_k$.

Let us briefly discuss the existence of such a time scale $s(\varepsilon)$. Note that the collection of functions $\tau_k^{r, \varepsilon}$ for those $k$ satisfying $S_k^r \subseteq \pi^{r + 1}(S_i), \nu^r(i,k) > 0$ are all asymptotically commensurate (the ratio of any two such functions has a positive limit as $\varepsilon \downarrow 0$) as a result of the definition of $\nu^{r}(i, k)$ obtained from Theorem   \ref{ergodic_theorem_second}. If a given cluster $S_k^r$ satisfying $\nu^r(i, k) >0$ is comprised of at least two distinct subclusters of rank $r - 1$, then, by the calculations in Theorem \ref{prelim_things_abt_rank_r_times} and Lemma \ref{rank_r_expectation_lemma}, the expected exit times are asymptotically smaller than $\tau_k^{r, \varepsilon}$, and we can always choose an intermediate time scale $s(\varepsilon)$. However, it may be the case that $S_k^r$ consists of just one cluster of rank $r - 1$. In this case, we investigate whether there is a rank $r' < r$ such that $S_k^r$ is the union of at least two distinct clusters of rank $r'$. If no such $r'$ exists, then $S_k^r = S_{j}$ for some $1 \le j  \le M$, and $\nu^r(i, k)$ will serve as the limit of $\mathbb P \left({X_{t(\varepsilon)}^{x, \varepsilon} \in S_j}\right)$ as $\varepsilon \downarrow 0$ uniformly in $x \in S_i$. When such a rank $r'$ does exist, then, due to the same explanation as above, the expected exit times of such subclusters of rank $r'$ are asymptotically smaller than $\tau_k^{r, \varepsilon}$, and a time scale $s(\varepsilon)$ can be chosen.

Let us fix an index $k$ such that $ \nu^r(i, k) >0.$ Suppose that there are at least two distinct clusters of rank $r-1$ contained within $S_k^r$ (the case when $r ' < r - 1$ is similar, and in the case when $S_k^r = S_j$ for some $j$, we have $ \lim \limits_{\varepsilon \downarrow 0} \mathbb P \left({X_{t(\varepsilon)}^{x, \varepsilon} \in S_j }\right) = \nu^{r}(i, k)$  uniformly in $x \in S_i$).  By the induction hypothesis, there is a measure $\overline {\nu}^r(k, \cdot )$ on the set $\{ j : 1 \le j \le M, S_j \subseteq S_k^r \}$ such that \begin{equation} \label{main_proof_formula_2}\lim \limits_{\varepsilon \downarrow 0 } \mathbb P \left({X^{y, \varepsilon}_{\overline t (\varepsilon)}\in S_j}\right) = \overline {\nu}^r(k, j ) \end{equation}uniformly in $y \in S_k^r$ and all time scales $\overline t (\varepsilon)$ satisfying $ s (\varepsilon) \le \overline t (\varepsilon) \le 2 s(\varepsilon) $.

 Let $ \boldsymbol{\sigma}^{x, \varepsilon} $ be the time of the first renewal of the process $X_t^{x, r-1, \varepsilon}$ after time $t(\varepsilon) - 2 s(\varepsilon).$ Due to Lemma \ref{no_skip_large_interval_small_escape_prob} and \eqref{main_proof_formula_1} applied to the time $t(\varepsilon) - 2s(\varepsilon)$,  $$ \lim \limits_{\varepsilon 
 \downarrow 0 }\mathbb P \left({ X_{t(\varepsilon) - 2 s(\varepsilon)}^{x,  \varepsilon} \in S_k^r, \boldsymbol{\sigma}^{x, \varepsilon} \le t(\varepsilon) - s(\varepsilon), X_{\boldsymbol{\sigma}}^{x, \varepsilon} \in S_k^r }\right)  = \lim \limits_{\varepsilon 
 \downarrow 0 } \mathbb P \left({ X_{\boldsymbol{\sigma}}^{x,  \varepsilon} \in S_k^r  }\right) = \nu^{r}(i, k)$$ uniformly in $x \in \pi^{r + 1}(S_i)$. Let $j$ be such that $S_j \subseteq S_k^r.$ By applying the Markov property to the stopping time $\boldsymbol{\sigma}^{x, \varepsilon }$ and using \eqref{main_proof_formula_2}, we conclude that $$ \lim \limits_{\varepsilon \downarrow 0 } 
    \mathbb P \left({ X^{x, \varepsilon}_{\boldsymbol{\sigma}} \in S_k^r, X^{x, \varepsilon}_{ t(\varepsilon)} \in S_j  }\right) = \nu^r(i, k)\overline {\nu}^r(k, j )  
$$  uniformly in $x \in \pi^{r + 1}(S_i)$.  Therefore,
$$ \liminf \limits_{\varepsilon \downarrow 0 } 
    \mathbb P \left(X^{x, \varepsilon}_{ t(\varepsilon)} \in S_j\right) \geq \nu^r(i, k)\overline {\nu}^r(k, j )  
$$
if $S_j \subseteq S_k^r$ and  $ S_k^r \subseteq \pi^{r+1}(S_i) $. Since the sum of the expressions in the right-hand side over such $j$ and $k$ is equal to one, we obtain

$$ \lim \limits_{\varepsilon \downarrow 0 } 
    \mathbb P \left({ X^{x, \varepsilon}_{ t(\varepsilon)} \in S_j  }\right) = \nu^r(i, k)\overline {\nu}^r(k, j )  
$$  uniformly in $x \in \pi^{r + 1}(S_i.)$   

If $t_1(\varepsilon) \le t_2(\varepsilon)$ are two time scales whose ratio is bounded away from zero and infinity, and which satisfy the assumptions imposed on $t(\varepsilon)$, then this same limit holds uniformly in all time scales $t_1(\varepsilon) \le \overline t (\varepsilon) \le t_2(\varepsilon)$, since, otherwise, we would be able to select a time scale satisfying our assumptions for which this limit would not hold. 

\subsubsection*{ Case 2: $t(\varepsilon) \ll \tau_k^{r, \varepsilon}$ for some $S_k^r \subseteq \pi^{r + 1}(S_i)$.}

Let $E = \left\{{ k : S_k^r \subseteq \pi^{r + 1  }(S_i), t(\varepsilon) \gg \tau_k^{r, \varepsilon} }\right\}, E' = \left\{{ k : S_k^r \subseteq \pi^{r + 1  }(S_i), t(\varepsilon) \ll \tau_k^{r, \varepsilon} }\right\} $. With an abuse of notation, for $x \in \bigcup_{k \in E} S_k^r$, let $\boldsymbol{\sigma}^{x, \varepsilon}$ be now defined as the first hitting time of $\bigcup_{k \in E'} S_k^r$. 
Then $\lim_{\varepsilon \downarrow 0} \mathbb{P} (\boldsymbol{\sigma}^{x,\varepsilon} \leq t(\varepsilon)/2) =1$ uniformly in $x \in \bigcup_{k \in E}S_k^r$.

We claim that there is a  measure $\nu^{r}(i, \cdot)$ on $ \{ j : S_j \subseteq { \bigcup_{k \in E'} S_k^r   } \}$ such that  \begin{equation} \label{main_proof_formula_6}
   \lim \limits_{\varepsilon \downarrow 0 } \mathbb P \left({{ X}_{\boldsymbol{\sigma}}^{x, \varepsilon} \in S_j}\right) = \nu^r(i, j) 
\end{equation} uniformly in $x \in S_i$. This  is similar to Theorem \ref{ergodic_theorem_third} with the difference being that  Theorem \ref{ergodic_theorem_third} would give the result for $S^r_k$ instead of individual sets $S_j$. To get this slightly finer result, we can consider an auxiliary Markov renewal process obtained from $X^{x,\varepsilon}_t$ by 
changing the transition times and transition probabilities out of $S_j$, for $j$ such that $  S_j \subseteq { \bigcup_{k \in E'} S_k^r   }$, in such a way that all such $S_j$ become clusters of rank $r$. Then (\ref{main_proof_formula_6}) follows from Theorem \ref{ergodic_theorem_third} applied to the process on clusters of rank $r$.


By the induction hypothesis, for each $j$ such that $S_j \subseteq \bigcup_{k \in E'} S_k^r$,  there is a measure $\overline{\nu}^r(j, \cdot)$ on $ \{ j' : S_{j'} \subseteq \pi^r(S_j) \}$ such that $$ \lim \limits_{\varepsilon \downarrow 0} \mathbb P \left({X_{\overline{t}(\varepsilon)}^{x, \varepsilon} \in S_{j'}}\right) = \overline {\nu}^r(j, j') $$ uniformly in $x \in S_j$ and every time scale $ t(\varepsilon)/2 \le \overline {t}(\varepsilon) \le t(\varepsilon)$.

Let us fix $j'$ such that $S_{j'} \subseteq \bigcup_{k \in E'} S_k^r $. By applying the Markov property to the stopping time $\boldsymbol{\sigma}^{x, \varepsilon}$ and summing over the location of $X_{\boldsymbol{\sigma}}^{x, \varepsilon}$ inside $\pi^r(S_{j'})$,   we obtain that $$ \liminf \limits_{\varepsilon \downarrow 0 }  \mathbb P \left({ X_{t(\varepsilon)}^{x, \varepsilon} \in S_{j'}  }\right) \ge \sum_{j : S_j \subseteq \pi^r(S_{j'}) } \nu^r(i, j) \overline {\nu}^r(j, j')  $$ uniformly in $x \in S_i$. Since the terms on the the right hand side add up to $1$,  we conclude that   $$ \lim \limits_{\varepsilon \downarrow 0 }  \mathbb P \left({ X_{t(\varepsilon)}^{x, \varepsilon} \in S_{j'}  }\right) = \sum_{j : S_j \subseteq \pi^r(S_{j'}) } \nu^r(i, j) \overline {\nu}^r(j, j')  $$ uniformly in $x \in S_i$.

As before,  if $t_1(\varepsilon) \le t_2(\varepsilon)$ are two time scales whose ratio is bounded away from zero and infinity, and which satisfy the assumptions imposed on $t(\varepsilon)$, then this limit holds uniformly in all time scales $t_1(\varepsilon) \le \overline t (\varepsilon) \le t_2(\varepsilon)$.
This justifies the inductive step and completes the proof.
\qed 
\\
\\
\noindent {\bf \large Acknowledgments}:  Leonid Koralov was supported by the NSF grant DMS-2307377 and by the Simons Foundation Grant MP-TSM-00002743. 
\\
\\


\begin{thebibliography}{999999}




\bibitem{AS} Asmussen S., {\it Applied probability and queues}, Springer-Verlag, (2003). 

\bibitem{Bakhtin1}
Bakhtin Y., {\it
Small noise limit for diffusions near heteroclinic networks},
Dynamical Systems 25 (2010), no. 3, 413--431.

\bibitem{Bakhtin2}
Bakhtin Y., {\it Noisy heteroclinic networks}, 
Probab. Theory Related Fields 150 (2011), no. 1--2, 1--42.

\bibitem{Bakhtin3}
Bakhtin Y., Chen H., Pajor-Gyulai Z., {\it Rare transitions in noisy heteroclinic networks},  to appear in the Memoirs of the American Mathematical Society.

\bibitem{BL1} Betz V., Le Roux S., 
{\it Multi-scale metastable dynamics and the asymptotic stationary distribution of perturbed Markov chains}, 
Stochastic Processes and their Applications, Volume 126, Issue 11, November 2016, Pages 3499--3526.

\bibitem{Bou}  Bouchet, F..  Reygner, J., {\it 
Generalisation of the Eyring–Kramers
Transition Rate Formula to Irreversible
Diffusion Processes}, Ann. Henri Poincare, 17 (2016), 3499–3532.

\bibitem{BOV1} Bovier, A., den Hollander, F., {\it Metastability: A Potential-Theoretic Approach}, Springer International Publishing, (2015). 

\bibitem{BE} Bovier, A., Eckhoff, M., Gayrard, V., Klein, M., {\it Metastability in reversible diffusion processes I: Sharp asymptotics for capacities and exit times}, J. Eur. Math.
Soc. 6(4), 399–424 (2004).



\bibitem{CINL} \c Cinlar E., {\it Introduction to Stochastic Processes}, Dover Publications, Inc., (2013). 

\bibitem{EREN} E W., Ren W., Vanden-Eijnden E., {\it Energy landscapes and rare events}, Proceedings of
the International Congress of Mathematicians, Vol. I (Beijing, 2002), pp. 621–630. Higher Ed. Press, Beijing (2002). 

\bibitem{EST} Eston V. R., Galves A., Jacobi C. M., Langevin R., {\it Dominance switch between
two interacting species and metastability}, Atas do II Simposio dos Ecossistemas da costa sul Brasileira, Sao Paulo: CACIESP,  (1988).

\bibitem{Ey} Eyring, H., {\it  The activated complex in chemical reactions}, J. Chem.
Phys. 3(2), 107–115 (1935).

\bibitem{F5} Freidlin M. I., {\it Sublimiting distributions and stabilization of solutions of parabolic equations
with a small parameter}, Soviet Math. Dokl., 18, No 4 (1977), 1114--1118.

\bibitem{FK1} Freidlin M. I., Koralov L., {\it  Metastable Distributions of Markov Chains with Rare Transitions}, 
Journal of Statistical Physics,  Volume 167, pages 1355 -- 1375 (2017).

\bibitem{FKnew1}  Freidlin M. I., Koralov L., {\it  Perturbations of parabolic equations and diffusion processes with degeneration: boundary problems, metastability, and homogenization},
Ann. Probab. 51 (2023), no. 5, 1752–1784.

\bibitem{FKnew2}  Freidlin M. I., Koralov L., {\it Metastability in Parabolic Equations and Diffusion
Processes with a Small Parameter}, arXiv:2403.12333 

\bibitem{FW2}  Freidlin M. I., Wentzell A. D.,  {\it On small random perturbations of dynamical systems}, Russian Math. Surveys, 25 No. 1 (1970) 1 - 55. 

\bibitem{FW} Freidlin M. I., Wentzell A. D., {\it Random
Perturbations of Dynamical Systems}, third edition, Springer 2012.

\bibitem{He} Helffer, B., Klein, M., Nier, F., {\it Quantitative analysis of metastability in
reversible diffusion processes via a Witten complex approach}, Math. Contemp. 26, 41–86 (2004).


\bibitem{KAUF} Kauffman, S., {\it The Origins of Order}, Oxford University Press, Oxford (1993).

\bibitem{Kra}  Kramers, H.-A., {\it Brownian motion in a field of force and the diffusion model of
chemical reactions},  Physica 7(4), 284–304 (1940).

\bibitem{Land}
Landauer, R., Swanson, J.A., {\it Frequency factors in the thermally activated
process}, Phys. Rev. 121, 1668 (1961).

\bibitem{Lang} Langer, J.S., {\it Statistical theory of the decay of metastable states}, Ann.
Phys. 54(2), 258–275 (1969).


\bibitem{Lan1} Landim C., Lee J., Seo I., {\it Metastability and time scales for parabolic equations with drift 1: the first time scale
}, Digital Preprint,  	arXiv:2309.05546 (2023). 

\bibitem{Lan2} Landim C., Lee J., Seo I., {\it Metastability and time scales for parabolic equations with drift 2: the general time scale
}, Digital Preprint,  		arXiv:2402.07695 (2024).

\bibitem{LX} Landim C.,  Xu T., {\it Metastability of finite state Markov chains: a recursive procedure
to identify slow variables for model reduction}, ALEA Lat. Am. J. Probab.
Math. Stat. 13 (2016), no. 2, 725--751.





 \bibitem{OLIV} Olivieri E., Vares M.E., {\it Large Deviations and Metastability}, Encyclopedia of Mathematics
and Its Applications, vol. 100. Cambridge University Press, Cambridge (2005).



\bibitem{SCH} Schütte, C., Sarich, M., {\it  Metastability and Markov State Models in Molecular Dynamics}, 
Courant Lecture Notes in Mathematics, vol. 24, Courant Institute of Mathematical Sciences/
American Mathematical Society, Providence/New York (2013). 

\bibitem{WE} Wentzell A. D., {\it On the asymptotics of eigenvalues of matrices with elements of order $\exp(-V_{ij}/(2 \varepsilon^2))$}, Soviet Math. Dokl., 13, No 1 (1972), 65-68.  



%
















































%





\end{thebibliography}
\end{document}